\documentclass[runningheads]{llncs}

\usepackage{subfigure}
\usepackage{csvsimple}
\usepackage{color}
\usepackage{graphicx}
\usepackage{braket}
\usepackage[final]{fixme}
\usepackage{amsmath,amsfonts,amssymb}
\usepackage{multirow}

\newcommand{\half}{\frac 12}

\newcommand{\dt}{\;\textnormal{d}t}
\newcommand{\ds}{\;\textnormal{d}s}
\newcommand{\dx}{\;\textnormal{d}x}

\newcommand{\bbf}{\mathbf{b}}
\newcommand{\fp}[2]{\frac{\partial #1}{\partial {#2}}}
\newcommand{\ltwoomega}{{0,\Omega}}

\newcommand{\ltwospace}{{L^2(\Omega)}}
\newcommand{\ltwoV}{{L^2([0,1], V)}}
\newcommand{\ltwoVh}{{L^2([0,1], V_h)}}

\newcommand{\ltwonorm}[1]{\| #1 \|_{0, \Omega}}

\newcommand{\diff}{\,\textnormal{d}}

\newcommand{\pp}[2]{\frac{\partial #1}{\partial #2}}

\newcommand{\DiffV}{{\textnormal{Diff}_V}}

\newcommand{\bnormi}[1]{\| #1 \|_{2, \Gamma_i}}
\newcommand{\bnorm}[1]{\| #1 \|_{2, \Gamma}}

\newcommand{\bspace}{{L^2(\Gamma)}}
\newcommand{\Ltwogamma}{{L^2(\Gamma)}}
\newcommand{\ltwogamma}{{0,\Gamma}}
\newcommand{\ltwogammai}{{0,\Gamma_i}}
\newcommand{\dxhat}{\;\text{d}\hat{x}}

\newcommand{\idop}{\textnormal{\textbf{id}}}

\newcommand{\mdiv}{\,\textnormal{div}_t}
\newcommand{\diver}{\textnormal{div}}

\newcommand{\ddv}{\Phi_t}
\newcommand{\ugrad}{u \cdot \nabla}

\newcommand{\pt}{\frac{\partial}{\partial t}}

\newcommand{\Gamp}{{\omega}}

\newcommand{\Gbo}{{G(\bbf;\Omega)}}
\newcommand{\Ebo}{{E(\bbf;\Omega)}}

\newcommand{\Gboz}{{\stackrel{\circ}{G}(\bbf;\Omega)}}
\newcommand{\Gbohz}{{\stackrel{\circ}{G}_h}}

\newcommand{\traceopg}{\gamma_\Gbo}
\newcommand{\traceoph}{\gamma_{H^1(\Omega)}}

\newcommand{\bgrad}{\bbf \cdot \nabla_t}

\newcommand{\bdiv}[1]{\mdiv ( \bbf #1)}

\newcommand{\gnorm}[1]{\| #1 \|_{\Gbo}}
\newcommand{\enorm}[1]{\| #1 \|_{\Ebo}}

\newcommand{\inv}{^{-1}}

\newcommand{\polyone}{\mathcal{P}^1(K)}

\newcommand{\pone}{\mathcal{P}^1(\Omega_h)}

\newcommand{\transp}{^\top}

\newcommand{\vertiii}[1]{{\left\vert\kern-0.25ex\left\vert\kern-0.25ex\left\vert #1
    \right\vert\kern-0.25ex\right\vert\kern-0.25ex\right\vert}}

\newtheorem{assumption}[theorem]{\rm{\sc Assumption}}

\begin{document}
\title{Space-time metamorphosis}
\titlerunning{Space-time Metamorphosis}
\author{Andreas Bock \and Colin Cotter}
\authorrunning{Bock and Cotter}
\institute{Imperial College London}


\maketitle

\begin{abstract}
{We study the problem of registering images. The framework we use is
metamorphosis and we construct a variational Eulerian space-time setting and
pose the registration problem as an infinite-dimensional optimisation problem. The
geodesic equations correspond to a system of advection and continuity equations
and are solved analytically. Well-posedness of a primal conforming finite
element method is established and its convergence is investigated numerically.
This provides a discrete forward operator for the matching parameterized by a
space-time velocity field. We propose a gradient descent method on this control
variable and show several promising numerical results for this approach.}
{Shape analysis; Metamorphosis; Finite element method.}
\end{abstract}

\section{Introduction}

In shape analysis, \emph{metamorphosis}
\cite{trouve2005metamorphoses,holm2009euler} is a metric framework for shape
matching between potentially topologically different shapes. This image
registration framework is a special case of Grenander's \emph{group action
model} using the theory of deformable templates
\cite{grenandergeneral,grenander2007pattern}. In this setting, registration (or
\emph{matchings}) between different shapes (or \emph{deformable templates}) are
found by identifying the action in some group of transformations so that the
application of the group action maps between the two shapes. A particular
strength of the group action framework is its ability to handle the non-linear
nature of shapes.  Equipping a shape space with a metric not only provides a
quantitative measure of closeness but also generates a meaningful way of
interpolating, in sense of the prescribed metric, between elements of the space.
In this paper, shapes are images given as functions with compact support on some
polygonal Lipschitz domain where the goal is to match a \emph{template} image
$I_0$ with a \emph{target} $I_1$ using metamorphosis. We aim to solve the
geodesic equations associated with the energy of the matching in a space-time
Eulerian frame and discretize these partial differential equations (PDEs) using
the finite element method. Defining solution spaces over space-time allows us to
leverage computational platforms featuring parallelism in both space and time at
the expense of dealing with a larger system. This then defines a forward
operator taking as input a time-dependent family of velocity fields and provides
the associated geodesics for the images, motivating a gradient descent method to
solve the inverse problem for the space-time velocity. Applications are manifold
in everything from medical imaging to music videos 
\cite{maurer1993review,beg2005computing,beg2006computing,wang2007large,du2011whole,fischer2008ill,beier1992feature,bajcsy1989multiresolution}.\\

Diffeomorphic matching frameworks such as the popular \emph{large deformation
diffeomorphic metric mapping} (LDDMM) approach \cite{beg2005computing} preserve
the underlying topology of the shapes on which they act. This class of methods
relies on an observation in \cite{arnold1966geometrie} that spaces of
smooth vector fields $V$ generate a subgroup of diffeomorphisms $\DiffV$ of some
shape space via the following ODE:
\[
\dot{\varphi}_t = u_t \circ \varphi_t,\qquad t \in [0,1],
\]
where $u_t \in V$, $\varphi_t\in\DiffV$ and $\dot{g} = \frac{\partial
g}{\partial t}$.  A match between two shapes $I_0$ and $I_1$ is therefore given
by the curve of diffeomorphisms that minimize the quantity:
\[
\int_0^1 \|u_t\|_V^2 \diff{t}+  c \|I_1 - I_0 \circ \varphi_1\inv\|^2,
\]
for some penalty parameter $c>0$ where the first term is referred to as the
\emph{kinetic energy}. This approach was, to the best of the authors' knowledge,
first used for matching problems in
\cite{christensen1996deformable,trouve1995infinite} with further early
applications in \cite{trouve1998diffeomorphisms}. $\|\cdot\|_V$ is typically
some higher-order Sobolev norm or induced by a smooth kernel.  See also
\cite{younesshapes} for a textbook on shapes and diffeomorphisms and
\cite{bauer2014overview,bruveris2018riemannian} for an overview of the
diffeomorphism group and Riemannian geometry for shape analysis.
When a purely diffeomorphic matching is no longer possible or necessary,
metamorphosis allows for shape matching with topological changes. An early
example of this is in \cite{grenander1994representations} where the growth 
represents a tumour appearing on a medical image over the course of time. For
the technical development of metamorphosis we refer to
\cite{miller2001group,miller2002metrics}. \cite{trouve2005local} takes a more
geometric point of view and describes in great detail a Riemannian
construction.\\

First, we introduce some basic notation in order to sketch the main ideas of
this paper and our contribution to the literature. $\Gamp = [0,1]^d$ denotes the
spatial $d$-dimensional unit domain representing the spatial dimension of the
images we are matching and we let $\Omega = [0,1]\times\Gamp$, where the first
coordinate denotes time. $\Gamp$ will in certain cases be periodic in its $d$
dimensions. We let $\Gamma_0$ and $\Gamma_1$ denote the $d$-dimensional
submanifolds of $\Omega$ given by the restriction of $\Omega$ in the first
coordinate to the boundaries i.e. $t=0$ and $t=1$ specifying where we impose the
images $I_j$, $j=0,1$ that we aim to match. When $d=1$, $\Omega$
is a simple unit square and $\Gamma_j$, $j=1,0$ are simply opposite boundaries;
for $d=2$ these are opposing facets of the unit cube.  Further we let $\Gamma =
\Gamma_0\cup\Gamma_1$.  Similarly, we shall sometimes use the notation $\Omega_t
=\Gamp$, $t\in[0,1]$.  For ease of notation we define the operator $\pp{}{t} +
\ugrad \equiv \bgrad$ taking values in a real Hilbert space $G$ that we define
later on. Here we use the $d$-dimensional row vector $\bbf = (u,\,
1)\transp$, $\nabla_t= (\nabla, \frac{\partial}{\partial t})$ and $-\mdiv$ its adjoint.
In the following $\ltwonorm{\cdot}$ is the $L^2(\Omega)$ norm, and note that
this is over space-time. Furthermore, let $\ltwoV$ denote functions $v$ such
that for $t\in [0,1]$, $t\mapsto \| v(\cdot, t)\|_V$ is $L^2$ and 
for $x\in \omega$, $x\mapsto v(x, t)$ is in $V$. We can compute a
\emph{metamorphic matching} between $I_0$ and $I_1$ be solving the following
infinite-dimensional optimisation problem:
\begin{subequations}\label{abstract_metamorphosis_problem}
\begin{align}
\inf_{u\in \ltwoV, I\in G}\quad & \frac 12 \int_0^1 \half\| u
\|_V^2 \dt + \sigma^{-2} \ltwonorm{\bgrad I}^2\label{abstract_metamorphosis_problem:fnl}\\ \text{subject to} \quad &
I|_{\Gamma_j} = I_j, \quad j=0,1,\label{abstract_metamorphosis_problem:cst}
\end{align}
\end{subequations}
where $u$ and $I$ depend both on space and time and $\sigma^{-2}>0$ is a penalty
parameter. \eqref{abstract_metamorphosis_problem:cst} is understood in an
$L^2(\Omega)$ sense. We aim to solve a type of relaxation of
\eqref{abstract_metamorphosis_problem} and implicitly define $I$ by the velocity
via a forward operator in the form of a PDE. Indeed for fixed $u$, the second
term in \eqref{abstract_metamorphosis_problem} is a least-squares advection
problem for $I$. With the previous idea in mind we write
\eqref{abstract_metamorphosis_problem} as:
\begin{subequations}\label{metamorphosis_problem_relaxed}
\begin{align}
\inf_{(u,z)\in \ltwoV\times\ltwospace}\quad & \half\int_0^1\| u \|_V^2 \dt +
\sigma^{-2} 
\ltwonorm{z}^2\label{metamorphosis_problem_relaxed:fnl}\\
\text{subject to} \quad&
z = \bgrad I^*[u]\quad \textnormal{in } L^2(\Omega)\\
& I^*[u] \triangleq \arg \inf_{\substack{I\in G,\\ I|_{\Gamma_0}=I_0,\\
I|_{\Gamma_1}=I_1}} \| \bgrad I \|^2_{0,\Omega}.\label{metamorphis_problem_relaxed:inner}
\end{align}
\end{subequations}
We aim to solve \eqref{abstract_metamorphosis_problem} in section
\ref{sec:mm} and determine suitable function spaces $V$ and $G$ defined over
\emph{space-time} so that $\bgrad : G\rightarrow \ltwospace$.  Since the
optimisation problem in the constraint is convex in $I$ we can differentiate the
functional in \eqref{metamorphis_problem_relaxed:inner} with respect to a
variation in $I$ and find necessary and sufficient optimality conditions for an
'optimal' image.  Section \ref{inner} is devoted to studying this \emph{inner
problem}, meaning we show that each image metamorphosis solution $I\in G$ can be
be parameterized by a certain vector field $u$. Next, section
\ref{outer_problem} shows that since $I$ satisfies a system of linear equations
\emph{implicitly defined} for any candidate $u\in\ltwoV$, we can cast
\eqref{metamorphosis_problem_relaxed} as an inverse problem in the $u$ variable
to which a gradient method can be applied.  We call this the \emph{outer
problem} and show that on certain subspaces, minimizers of
\eqref{metamorphosis_problem_relaxed:fnl} are also minimizers of
\eqref{abstract_metamorphosis_problem:fnl}. Section \ref{outer:conclusion}
concludes this paper and highlights our accomplishments as well as the main
challenges going forward. Our contributions are as follows. We derive several
theoretical results for the inner problem leading to a primal finite element
method. We study its convergence properties numerically, and based on the
relaxation in \eqref{metamorphosis_problem_relaxed} propose a gradient descent
scheme on the velocity. Several numerical results are shown for one and
two-dimensional images.\\

We highlight some differences between the setting we propose and the literature.
Neither the space-time nor the finite element approaches are novel here.
Generally speaking, space-time methods are advantageous when a \emph{one-shot}
(i.e.  solving a PDE for all time and space simultaneously) solution approach is
preferred to a time-marching scheme. 
Early works applying space-time finite elements include
\cite{hughes1988space,hulbert1990space} in the context of hyperbolic problems
and elastodynamics, motivated by the desire to resolve discontinuities in the
solution via adaptivity. This approach permits a so-called \emph{unstructured}
mesh which means that the elements form an irregular pattern as opposed to e.g.
a uniform quadrilateral discretization in space and time. We are therefore
adopting this framework to not only develop a more computationally expedient
method but also to pave the way towards an adaptive strategy, where local mesh
resolutions in space-time can help accurate resolve sharp image gradients. See
the recent textbook \cite{olaf} for more general space-time methods for PDEs.
Finite element methods have appeared before in the computational anatomy
community \cite{boyd2009image,gunther2011direct}. Much of this work uses
anatomical models i.e. applying finite element methods to model tissues from
volumetric data. To the best of our knowledge using finite elements to directly
model both the image and velocity in an Eulerian setting has not yet been
investigated and this forms our contribution to the literature. Finally, very
recently, \cite{effland} develops a shooting method for metamorphosis of weakly
differentiable images via a time-discrete forward map. The regularity of
deformation is ensured through cubic splines on a coarse mesh, while the image
is discretized using bilinear elements on a fine mesh.\\

\section{Variational Space-time Method}\label{sec:mm}

Before describing the inner problem \eqref{metamorphis_problem_relaxed:inner} in
the next section we begin with a standard assumption on the velocity. Throughout
this paper, $u\in\ltwoV$ will be a velocity field and we state the following
assumptions:
\begin{assumption}\label{ass1}
\hspace{5cm}
    \begin{enumerate}
    \item $V$ is a Lipschitz function space ($W^{1,\infty}(\Gamp)$) such that for $t
        \in [0,1]$, $x\mapsto u(x, t)$ is a Lipschitz mapping in $x$ and $t\mapsto \| u(\cdot, t)\|_V$ is $L^2([0,1])$. When $\Gamp$ is not
        periodic we restrict further to $u(\cdot, t) \in H_0^1(\Gamp) \cap
        W^{1,\infty}(\Gamp)$, where $H_0^1(\Gamp)$ is the kernel of the trace
        map on $\partial \Gamp$.
    \item Denote by $\varphi: [0,1]\rightarrow \textnormal{Diff}(\Gamp)$ the \emph{flow map} associated with each
        $u \in \ltwoV$ such that for any point
        $\hat{x}\in\Gamma_0$, $\varphi$ solves the following system
        uniquely, see \cite{arnold1966geometrie}:
        \begin{equation}\label{inner:varphi}
        \begin{cases}
          \dot{\varphi}_t(\hat{x}) = u(\varphi_t(\hat{x}), t), \quad t \in
          [0,1],\\
          \varphi_0 = \idop.
    \end{cases}
    \end{equation}
\end{enumerate}
\end{assumption}

The work \cite{dupuis1998variational} shows that $H^3(\Gamp)$ is sufficient to ensure
the spatial Lipschitz property in assumption \ref{ass1}, and we let $V$ be
comprised of functions $v$ that are finite in the following norm, where $\langle
\cdot,\cdot\rangle_{0, \Omega}$ denotes the inner product on $L^2(\Omega)$:
\begin{equation}\label{Vnorm}
\| v \|_V^2 = \langle L v,v\rangle_{0, \Omega}.
\end{equation}
Here, $L=(\idop - \alpha^2 \Delta)^3$ and $\alpha^2$ is a
length-scale. A candidate velocity field with the
regularity above is Lipschitz in space. The assumption above also establishes a
bijection between $\Gamma_0$ and $\Gamma_1$ in the sense that for any $\hat
x_0 \in \Gamma_0$ there exists a unique $\hat x_1\in\Gamma_1$ such that $\hat
x_0 = \varphi_1\inv(\hat x_1)$, and that for all times $t\in [0,1]$, $\varphi_t$
is a diffeomorphism of the domain $\Omega_t$.\\

Based on this assumption we now determine the image space $G$ induced by the
space-time vector field $\bbf$.  We define an \emph{energy semi-norm}:
\begin{equation}\label{energy_norm}
\enorm{f} = \ltwonorm{\bgrad f},
\end{equation}
and the \emph{graph norm}:
\begin{equation}\label{graph_norm}
\gnorm{f}^2=\ltwonorm{f}^2+\enorm{f}^2.
\end{equation}
The solution space for the images is then given by:
\begin{equation}\label{graph_space}
\Gbo = \{ f\in L^2(\Omega) \;|\; \gnorm{f} < \infty\},
\end{equation}
and $G=\Gbo$ when it is clear from the context. It is clear that $H^1(\Omega)
\subset \Gbo$. The spaces above admit the inner products:
\begin{align*}
& \langle i,j \rangle_\Ebo = \langle \bgrad i,\bgrad j \rangle_{0,\Omega},\\
& \langle i,j \rangle_\Gbo = \langle i,j \rangle_{0,\Omega} + \langle i, j \rangle_\Ebo.
\end{align*}
We also define the boundary norm:
\begin{equation}\label{boundary_norm}
\bnorm{f}^2 = \sum_{i=0,1} \bnormi{f}^2, \qquad
\bnormi{f}^2 = \int_{\Gamma_i} f^2 \diff \hat x,\; i=0,1.
\end{equation}

\subsection{Inner Problem}\label{inner}

In this section we solve the least-squared advection problem in
\eqref{metamorphis_problem_relaxed:inner} for a fixed $u$ satisfying assumption
\ref{ass1}:
\begin{subequations}\label{inner_problem}
\begin{align}
& \text{Find } I^* = \arg\inf_{I\in G} \half \ltwonorm{\bgrad
I}^2\label{inner_problem:fnl}\\
& \text{subject to}\quad I|_{\Gamma_j} = I_j, \quad
j=0,1.\label{inner_problem:cst}
\end{align}
\end{subequations}

A minimizer of this convex problem can be shown to satisfy a (necessary and
sufficient) coupled system between advection and continuity.  In section
\ref{subsec:analytical_solution} we solve these strong form equations
analytically and formulate a regularity theorem. Next, section
\ref{inner:lsq} contains the main results of this paper wherein we develop a variational
setting in which well-posedness is obtained and from which finite element
methods are derived. 

\subsubsection{Analytical Solution}\label{subsec:analytical_solution}

To derive the analytical solution to the convex problem described by problem
\ref{inner_problem} we assume for the moment that everything is smooth and that
a unique minimizer to this problem exists, then retrace our steps to impose the
minimum required regularity of our solutions. Substituting $\bgrad I = z$ in
\eqref{inner_problem} leads to the equivalent problem:
\begin{subequations}\label{inner_problem2}
\begin{align}
& \inf_{(I, z)\in G\times\ltwospace} \half \ltwonorm{z}^2\\
& \text{subject to } \;\bgrad I = z,\\
& \text{\phantom{subject to }}\quad I|_{\Gamma_j} = I_j, \quad j=0,1.
\end{align}
\end{subequations}
Now using space-time Lagrange multipliers $\phi\in G$ and $\mu\in L^2(\Gamma)$
for the respective constraints in \eqref{inner_problem2} we obtain the following
functional:
\[
F=\half \ltwonorm{z}^2 + \langle \phi, \bgrad I - z\rangle_{0, \Omega} 
+ \langle\mu,I^*-I\rangle_\ltwogamma,
\]
where $I^*$ represents the boundary data in the sense that:
\[ \langle \mu,I^*-I\rangle_\ltwogamma =\sum_{i=0,1}\langle
\mu,I_i-I\rangle_\ltwogammai.
\]
Taking variations in $F$ in the four variables i.e. $\delta F = 0$ and
inspecting the resulting equations:
\begin{subequations}\label{dual:fullsystem}
\begin{align}
& \langle z, \delta z\rangle_\ltwoomega = \langle \phi, \delta
z\rangle_\ltwoomega, \quad\forall \delta z \in G,\\
& - \langle z, \delta \phi\rangle_\ltwoomega 
+ \langle \bgrad I, \delta \phi\rangle_\ltwoomega= 0,\quad\forall \delta \phi \in G,\\
& - \langle \phi, \bgrad \delta I\rangle_\ltwoomega  -\langle \mu,\delta I
\rangle_\ltwogamma = 0, \quad\forall \delta I\in G,\\
& \langle\delta\mu,I^*-I\rangle_\ltwogamma = 0,\quad\forall \delta \mu \in \Ltwogamma.
\end{align}
\end{subequations}
The first equation states an equivalence between $z$ and $\phi$ so we can reduce
the system by simple substitution. Further, choosing $\mu = \bbf\cdot\eta
z|_\Gamma$ as ansatz we can integrate by parts in the third equation to write it
as:
\[
- \langle \mdiv(\bbf z), \delta I\rangle_\ltwoomega = 0, \quad\forall \delta
I\in G.
\]
In summary we arrive at: 
\begin{subequations}\label{eq:advcontsystem}
\begin{align}
& \bdiv{z} = 0,\label{eq:advcontsystem:z}\\
& \bgrad I = z,\\
& I|_{\Gamma_j} = I_j, \quad j=0,1,\\
& \mu = \bbf\cdot\eta z|_\Gamma.
\end{align}
\end{subequations}

Now any solution to \eqref{eq:advcontsystem} is also a solution to
\eqref{dual:fullsystem}. To solve this we need a property of the Jacobian
determinant of the flow map $\varphi_t$ generated by $u$. Let $J_t =
\frac{\partial\varphi_t}{\partial\hat x}$ denote the Jacobian of $\varphi_t$ at
$\hat x$ and the Jacobian determinant by:
\begin{equation}\label{eq:Phit}
\ddv = \det J_t.
\end{equation}
Since $\varphi_t\circ \varphi_t\inv(x)=x$, $\Phi_t$ is not degenerate and so the
inverse function theorem holds:
\begin{align*}
& \Phi_t (\hat x)  = (\Phi_t\inv\circ\varphi_t(\hat x))\inv,\\
& \Phi_t\inv (x)  = (\Phi_t\circ\varphi_t\inv (x))\inv.
\end{align*}
This quantity can be seen to satisfy the continuity equation since for $t\in
[0,1]$ a time slice of the space-time domain $\Omega_t = \{ (x,t) \;|\; x \in
\Gamp \}$ we have:
\begin{equation}\label{eq:phiconst}
\int_{\Omega_t} \Phi_t\inv \dx = \int_{\Gamma_0} \Phi_t\inv\circ\varphi_t \Phi_t \dxhat =
\int_{\Gamma_0}\dxhat = \text{constant},
\end{equation}
implying:
\[
\frac{\textnormal{d}}{\textnormal{d}t} \int_{\Omega_t} \Phi_t\inv \dx = 0,
\]
so the integral of $\Phi_t\inv$ is conserved. By Reynold's transport theorem
\cite{chorin1990mathematical} (differentiation under the integral sign and using
the divergence theorem) we obtain:
\begin{subequations}\label{phieq}
\begin{align}
& \bdiv{\Phi_t\inv} = 0,\\
& \Phi_0\inv = 1.
\end{align}
\end{subequations}
Therefore, since the governing equation \eqref{eq:advcontsystem:z} for $z$ is
the same it must be of the form $z = \lambda \Phi_t\inv$, where $\lambda$ solves
$\pp{\lambda}{t} + u \cdot \nabla \lambda = 0$ as shown by the following
computation:
\begin{align*}
\pp{z}{t} & = \pp{\lambda}{t}\Phi_t\inv + \lambda \pp{\Phi_t\inv}{t}\\
& = - u\cdot\nabla\lambda\Phi_t\inv - \lambda\diver(u \Phi_t\inv)\\
& = - \diver(u \lambda\Phi_t\inv)\\
& = - \diver(u z).
\end{align*}
so clearly $\bdiv{z} = \pp{z}{t} + \diver(u z) = 0$. $\lambda$ is constant along
streamlines which we denote by $\lambda_{\hat x}$ for each selection of initial
$\hat x \in \Gamma_0$. As a result:
\begin{equation}\label{strong:I:constant}
\pp{I}{t} + u \cdot \nabla I = \lambda_{\hat x} \Phi_t\inv\,.
\end{equation}
Letting $\gamma_t(\hat x) = \int_0^t \Phi_s\inv\circ \varphi_s(\hat x)\ds$ we
define the ansatz solution as a kind of interpolation between the boundary
conditions along the streamlines of $u$:
\begin{align}\label{I:ansatz}
I(\varphi_t(\hat x), t) := \tilde I(\hat x, t) & := \frac{\gamma_t(\hat x)}{\gamma_1(\hat
x)}\Big[I_1\circ\varphi_1(\hat x) - I_0(\hat x)\Big] + I_0(\hat x)\,.
\end{align}

Now we differentiate \eqref{I:ansatz} with respect to time yields. On the
left-hand side we get:
\begin{align*}
\frac{\textnormal{d}}{\textnormal{d}t} I(\varphi_t(\hat x), t) & = \nabla
I(\varphi_t(\hat x), t) \cdot \dot \varphi_t\circ\varphi_t(\hat x) +
\fp{I}{t}(\varphi_t(\hat x), t)\\
& = \nabla I(x, t) \cdot u(x, t) + \fp{I}{t}(x, t) \,,
\end{align*}
and on the right:
\begin{align*}
\frac{\textnormal{d}}{\textnormal{d}t} \frac{\gamma_t(\hat x)}{\gamma_1(\hat
x)}\Big[I_1\circ\varphi_1(\hat x) - I_0(\hat x)\Big] + I_0(\hat x) =
\lambda_{\hat x} \fp{\gamma_t}{t}(\hat x) = \lambda_{\hat x} \Phi_t\inv\circ \varphi_t(\hat x)\,,
\end{align*}
where $\lambda_{\hat x} = \frac{1}{\gamma_1(\hat x)}\Big[I_1\circ \varphi_1(\hat x) - I_0(\hat
x)\Big]$. Reverting to $I$ we see that \eqref{I:ansatz} satisfies
\eqref{strong:I:constant}, so the solution to \eqref{eq:advcontsystem} is given
by
\begin{subequations}\label{eq:Iz:analytical}
\begin{align}
& I(x, t) = \tilde I(\varphi_t\inv(x), t)\\
& z(x, t) = \lambda_{\varphi_t\inv(x)}\Phi_t\inv(x)\,.
\end{align}
\end{subequations}

\eqref{eq:Iz:analytical} is valid for $I_j \in L^2(\Gamma_j)$, $j=0,1$.  Next,
we state conditions that imply higher-order regularity of the solutions obtained
above.  Prescribing regularity of Jacobian determinants is in and of itself a
field of study, and is particularly non-trivial when restricted to functions
$W^{1,\infty}(\Gamp)$, see
\cite{dong1994prescribing,sickel1999characterisation}. This poses a challenge
for the study of regularity of the solutions obtained here. We provide one
example below of what we can expect when restricting $V$ to be a sufficiently
smooth subspace.

\begin{proposition}[Streamline regularity]\label{streamsmooth}
Suppose $\varphi_t \in \mathsf{C}^2(\Gamp)^d$ and assume $I_j
\in H^1(\Gamma_j)$, $j=0,1$. Then for any $\hat x\in\Gamma_0$ the map $t \mapsto
\tilde I(\hat x, t)$ has a continuous second derivative.\\
\end{proposition}
\begin{proof}
By using $\pp{\tilde I}{t} = \lambda_{\varphi_t\inv(x)}\Phi_t\inv\circ\varphi_t$ we
simply apply direct computation to arrive at the following expression for the
second derivative:
\begin{align*}
\frac{\partial^2\tilde I}{\partial t^2} & =
 \lambda_{\varphi_t\inv(x)} \frac{\text{d}}{\text{d}t}\Phi_t\inv\circ\varphi_t\\
& = \lambda_{\varphi_t\inv(x)} (\nabla \Phi_t\inv)\circ\varphi_t \cdot \pp{\varphi_t}{t}\\
& = \lambda_{\varphi_t\inv(x)} (\nabla \Phi_t\inv)\circ\varphi_t \cdot u.
\end{align*}
Now the $i$\textsuperscript{th} component of $\nabla \Phi_t\inv$, $i=1,...,d$, can
be expressed using Jacobi's formula (see e.g. \cite{chorin1990mathematical}) as:
\[
\big[\nabla \Phi_t\inv\big]_i = \Phi_t\inv\textnormal{trace}\Big[J_t \pp{J_t\inv}{x_i}\Big].
\]
The right-hand side is well-defined if the second derivative of the flow map is
continuous i.e. $\frac{\partial^2\varphi_t^{k,j}}{\partial x_j\partial x_k}\in
\mathsf{C}^0(\Gamp)$, for all $k,j=1,...,d$. Since it is a diffeomorphism the same
regularity holds for its inverse.
\end{proof}

\subsubsection{Primal Formulation}\label{inner:lsq}

This section studies the least-square advection problem described in problem
\eqref{inner_problem}. We examine the stability of the resulting least-squares
boundary value problem for the image by differentiating the functional in this
problem with respect to a variation in the image. What we present here is not
the first study of hyperbolic equations in the least-squares finite element
literature, so we first provide a brief historical review. The early work of
\cite{nguyen84} presents a least-squares space-time method for the
advection-diffusion equation. The direction we take here is mostly similar to
the work of Perrochet in 1995 \cite{perrochet1995space}. See the textbook
reference \cite{bochev2009least} for an overview of least-square FEM.  Our
problem differs through the boundary conditions we use which in the present
context are the images we are registering.\\ 

We start by stating the following proposition which is central to the analysis
carried out in this section. This follows from the more general results stated
in \cite[Section 10.2.1]{bochev2009least} for general graph norms.
\begin{proposition}\label{prop:Gbodense}
$\Gbo$ is a Hilbert space in which $H^1(\Omega)$ is dense.
\end{proposition}

We now state and prove results for these spaces that are considered standard for
Sobolev spaces: a trace inequality similar to that shown for linear first-order
PDEs in \cite{houston1999posteriori}, Poincar\'e-type inequality, a trace
theorem and well-posedness. We
highlight the reference \cite{de2004least} as a foundation on which our work
relies as we can adapt to the setting described therein with relative ease,
itself using similar techniques to \cite{manteuffel1999boundary} where a neutron
transport problem was interpreted as a least-squares problem.
\newcommand{\vpxs}{\varphi_s(\hat x)}
\newcommand{\cvpxs}{\circ\vpxs}
\newcommand{\vpx}{\varphi_t(\hat x)}
\newcommand{\cvpx}{\circ\vpx}

\begin{theorem}[Trace inequality]\label{def:inner_trace}
Given assumption \ref{ass1}, $\exists C > 0$ depending only on $\Omega$ and $u$
such that for any $I \in \Gbo$,
\[
\bnormi{I}^2 \leq C \gnorm{I}^2,\quad i=0,1\,.
\]
\end{theorem}
\begin{proof}
Let $\hat x \in \Gamma_0$ and let $f\in\mathsf{C}^\infty(\bar\Omega)$  such that the following
integration by parts formula holds:
\[
f(x,t) = f(\hat x, 0) + \int_0^t \pt{f}(\vpxs, s) + u_s\cvpxs\cdot\nabla f(\vpxs,
s)\diff s\,.
\]
In particular, for $f = p^2$, where $p$ is some smooth function we obtain:
\[
p^2(\hat x, 0) \leq p^2(x,t) + \int_0^t \left| \pt{p^2}(\vpxs, s) + u_s\cvpxs\cdot\nabla
(p^2)(\vpxs, s)\right| \diff s\,.
\]
Now integrating in time we get:
\begin{equation}\label{trace1}
1 \cdot p^2(\hat x, 0) \leq \int_0^1 p^2(x,t) + \int_0^1 \left|\pt{p^2}(\vpx, t)
+ u_t\cvpx\cdot\nabla (p^2)(\vpx, t)\right|\diff t\,.
\end{equation}
We integrate over $\Gamma_0$ on both sides of \eqref{trace1} and use
$p^2(x,t) = p^2(\varphi_t(\hat x), t)$:
\begin{align*}
\|p\|_{2,\Gamma_0}^2 &= \int_{\Gamma_0} \int_0^1
p^2(\varphi_t(\hat x), t)\diff t\diff\hat x\\
& + \int_{\Gamma_0}\int_0^1
\left|\pt{p^2}(\vpx, t) + u_t\cvpx\cdot\nabla (p^2)(\vpx, t)\right|\diff t\diff\hat
x\,.
\end{align*}
Examining the first term on the right-hand side:
\begin{align*}
\int_{\Gamma_0} \int_0^1 p^2(\varphi_t(\hat x), t)\diff t\diff\hat x & =
\int_{\varphi_t\circ \Gamma_0} \int_0^1 p^2(x, t)
\left|\frac{\partial\varphi_t\inv}{\partial x}\right|\diff t\diff x
\leq C_u \ltwonorm{p}^2,
\end{align*}
where we have used the fact that $\varphi_t\circ\Gamma_0 \times [0,1] = \Omega$,
for all $t\in [0,1]$ due to the spatial boundary conditions for $u$ and where $C_u =
\sup_{(x,t)\in\Omega}\left|\frac{\partial\varphi_t\inv}{\partial x}\right|$, which by
assumption \ref{ass1} holds.  Now the second term can also be bounded by the
same reasoning and an application
of H\"older's inequality:
\begin{align*}
& \int_{\Gamma_0}\int_0^1 \left|\pt{p^2}(\vpx, t) + u_t\cvpx\cdot\nabla (p^2)(\vpx,
t)\right|\diff t\diff\hat x\\
\leq & \int_{\varphi_t\circ\Gamma_0}\int_0^1 \left|\pt{p^2}(x, t) + u_t(x)\cdot\nabla
(p^2)(x, t)\right|\left|\frac{\partial\varphi_t\inv}{\partial x}\right|\diff t\diff x\\
\leq & C_u \ltwonorm{\bgrad{p^2}}^2\\
\leq & 2 C_u \ltwonorm{p}\ltwonorm{\bgrad{p}}.
\end{align*}
Now using Young's $\epsilon$-inequality on this and collecting the terms above
we get
\begin{align*}
\|p\|_{2,\Gamma_0}^2 & \leq \ltwonorm{p}^2 + \frac{C_u^2}2(\ltwonorm{p}^2 +
\ltwonorm{\bgrad{p}}^2)\\
& \leq (1+\frac{C_u^2}2) \gnorm{p}^2
\end{align*}
showing the result for $\Gamma_0$. Since the ODE \eqref{inner:varphi} is
time-reversible (compose $\varphi_t$ with $g(t)=-t$ in time and take the
derivative: the same equations are obtained), the same procedure can be carried
out if we follow the streamlines backwards in time from $\Gamma_1$ so adding
these proves the claim for smooth functions. Now we simply extend this from the
dense set $\mathsf{C}^\infty(\bar\Omega)$ to $\Gbo$ by a limiting procedure.
\end{proof}


\begin{theorem}[Poincar\'e inequality]\label{def:inner_poincare}
Given assumption \ref{ass1}, $\exists C > 0$ depending only on $\Omega$ and $u$
such that for any $I \in \Gbo$,
\[
\ltwonorm{I}^2 \leq C\Big(\bnorm{I}^2 + \ltwonorm{\bgrad I}^2 \Big)
\]
\end{theorem}
\begin{proof}
As before, let $\hat x \in \Gamma_0$ and let $f\in\mathsf{C}^\infty(\bar\Omega)$
be such that the streamline integration by parts formula from theorem
\ref{def:inner_trace} holds:
\[
f(x,t) = f(\hat x, 0) + \int_0^1\pt{f}(\vpxs, s) + u_s\cvpxs\cdot\nabla f(\vpxs,
s)\diff s\,.
\]
Square both sides, using Jensen's inequality:
\[
f^2(x,t) \leq 2 f^2(\hat x, 0) + \int_0^t \left|\pt{f}(\vpxs, s) + u_s\cvpxs\cdot\nabla f(\vpxs,
s)\right|^2\diff s\,.
\]
Now we multiply both sides of the equation by $\ddv\inv$ defined in
\eqref{eq:Phit} and integrate over time and $\Gamma_0$:
\begin{align*}
\int_{\Gamma_0} \int_0^1 f^2(x,t) \ddv\inv \diff t \diff\hat x &\leq \int_{\Gamma_0} f^2(\hat
x, 0) \ddv\inv\diff\hat x\\
& + \int_{\Gamma_0} \int_0^1 \left|\pt{f}(\vpx, t) + u_t\cvpx\cdot\nabla f(\vpx,
t)\right|^2 \ddv\inv\diff t \diff\hat x\,.
\end{align*}
Using the results from the theorem \ref{def:inner_trace} this inequality
simplifies to:
\begin{align*}
\ltwonorm{f}^2 & \leq \int_{\Gamma_0} f^2(\hat x, 0) \ddv\inv\diff\hat x +
\gnorm{f}^2\,,
\end{align*}
so we can state:
\begin{align*}
\ltwonorm{f}^2 & \leq (1+ C_u) ( \| f \|_{2, \Gamma_0}^2 + \gnorm{f}^2)
\end{align*}
as required, where $C_u$ is defined as in theorem \ref{def:inner_trace}. The
proof follows by density.
\end{proof}

We also have a trace theorem in the setting presented thus far.
\begin{corollary}[Trace theorem ]\label{def:inner_tracethm}
Given assumption \ref{ass1}, there exists a bounded surjective operator:
$\traceopg: \Gbo \rightarrow H^{1/2}(\Gamma)$.
\end{corollary}

\begin{proof}
Proposition \ref{prop:Gbodense} states that $H^1(\Omega)$ is dense in $\Gbo$ so
the standard trace operator $\traceoph: H^1(\Omega) \rightarrow H^{1/2}(\Gamma)$
admits a continuous extension which we denote by $\traceopg$. Boundedness
follows by theorem \ref{def:inner_trace}.
\end{proof}

This theorem does \emph{not} follow from the analytical solution
\eqref{I:ansatz}. Equation \eqref{I:ansatz} tells us that for any data $g_\Gamma
\in\bspace$, the analytical solution \eqref{I:ansatz} provides a function $g \in
\Gbo$ such that $g|_\Gamma = g_\Gamma$. This does not on its own guarantee that
we can evaluate the trace of an arbitrary function in $\Gbo$.\\

Equipped with corollary \ref{def:inner_tracethm} we can define the space $\Gboz$
of functions that vanish on $\Gamma$ as follows:
\begin{equation}\label{def:Gboz}
\Gboz = \{ f \in \Gbo \,|\, \traceopg f = 0 \}.
\end{equation}

\begin{remark}[Historical comments]
Theorem \ref{def:inner_poincare} was first called a \emph{curved Poincar\'e
result} in the early work \cite{perrochet1995space} but was not proved there.
We also mention that the theory of \emph{Friedrich systems}
\cite{houston1999posteriori} forms some of the basis for the analysis of
least-squared advection systems, or can at the very least be applied in certain
circumstances e.g. in the construction of trace theorems.
\end{remark}

We return to show well-posedness of \eqref{inner_problem}.

\begin{theorem}[Existence of a minimizer]\label{theorem:inner_problem_eu}
Given assumption \ref{ass1} there exists a unique minimizer to
\eqref{inner_problem} and the associated necessary first-order optimality
conditions are necessary and sufficient.
\end{theorem}
\begin{proof}
The functional \eqref{inner_problem:fnl} is clearly bounded from below and
admits a minimizing sequence $\mathcal{I} = \{ I_k \}_{k\geq 0}$.  Then by the
Poincar\'e estimate in theorem \ref{def:inner_poincare}, $\mathcal{I}$ is
bounded in $\Gbo$. Because this space is reflexive we find a subsequence of
$\mathcal{I}$ that converges weakly to a minimizer of \eqref{inner_problem:fnl}.
Convexity implies lower semi-continuity of the functional, and the unique
minimizer is obtained.
\end{proof}

The optimality conditions of \eqref{inner_problem} are described by
differentiating \eqref{inner_problem:fnl} with respect to an arbitrary variation
$\delta I \in \Gbo$ in the image $I$ and obtain what we call the \emph{primal
inner problem}. In the following, define the bilinear form $a(I,\delta I) = \langle I, \delta I\rangle_\Ebo$:
\begin{subequations}\label{leastsquareadvection}
\begin{align}
& \text{Find } I \in \Gbo\; \text{ such that:}\\
& a(I, \delta I) = 0, \qquad
\forall \delta I \in \Gbo,\label{pip}\\
& I|_{\Gamma_i} = I_i, \,\quad i=0,\,1.
\end{align}
\end{subequations}

\begin{theorem}[Existence, uniqueness and stability for the inner
problem]\label{def:inner_uniqueness}
Given assumption \ref{ass1}, for any boundary data $I_i \in H^{1/2}(\Gamma_i)$,
$i=0,1$, there exists a unique $I\in\Gbo$ solving
\eqref{leastsquareadvection} such that:
\[
\gnorm{I} \lesssim \| I \|_{2,\Gamma}.
\]
\end{theorem}
\begin{proof}
Let $g$ be any extension of the boundary data $I_i$, $i=0,1$ to the domain
via corollary \ref{def:inner_tracethm}. Thanks to theorem
\ref{def:inner_poincare} the inner product associated to $\enorm{\cdot}$ is
equivalent to $\gnorm{\cdot}$ on the subspace $\Gboz$. By the Riesz
representation theorem we know that we can find a unique $\xi \in \Gboz$ such
that:
\begin{equation}\label{eq:xi}
a(\xi, v) = -a(g, v),\qquad \forall v \in
\Gboz\,.
\end{equation}
By construction we have $I = \xi + g$. The bound is provided by the trace
operator cf. theorem \ref{def:inner_tracethm}.

\end{proof}

We now turn to discretizing problem \ref{leastsquareadvection} using a
conforming finite element space. Based on the results above the rest of the
analysis here is trivial. The Galerkin finite element space \cite{brennerscott}
$G_h \subset \Gbo$ denotes a finite-dimensional
$\mathsf{C}^0(\Omega)$ conforming subset of $\Gbo$ spanned by the standard
$\polyone$ shape functions associated with a shape-regular quasi-uniform
triangulation $\Omega_h$ of $\Omega$ \cite{ErnGuermond2013}.
To describe a discrete analogue of \ref{leastsquareadvection} we first project
the boundary conditions into the lower-dimensional space. Let $g$ be an
extension of the data introduced in theorem \ref{def:inner_uniqueness} and
definte $g_h$ by:
\[
a( g_h - g, \delta I) =0, \quad \forall \delta I \in G_h,
\]
and let $I_i^h = g_h|_{\Gamma_i}$, $i=0,1$. Then a conforming discretisation of
\ref{leastsquareadvection} is as follows:
\begin{subequations}\label{leastsquareadvection:disc}
\begin{align}
& \text{Find } I_h \in G_h\; \text{ such that:}\\
& a(I_h, \delta I) = 0, \qquad
\forall \delta I \in G_h,\\
& I_h|_{\Gamma_i} = I_i^h, \,\quad i=0,1,\label{leastsquareadvection:disc:bcs}
\end{align}
\end{subequations}
where the boundary conditions are understood in a weak sense. The bilinear form
in problem \ref{leastsquareadvection} remains coercive by conformity and hence
the results above remain valid and this finite-dimensional problem
\eqref{leastsquareadvection:disc} is well-posed. Similar to what was done in
theorem \ref{def:inner_uniqueness} we can write the solution $I_h$ of
\eqref{leastsquareadvection:disc} in terms of its homogeneous part $\xi_h$ as
follows:
\begin{equation}\label{eq:xih}
a(\xi_h, v_h) = -a(g_h, v_h),\qquad \forall
v_h \in
\Gbohz\,.
\end{equation}
This implies the bound described by the following proposition.
\begin{proposition}
Suppose assumption \ref{ass1} holds.  Let $I\in\Gbo$ solve
\eqref{leastsquareadvection} and $I_h\in G_h$ solve
\eqref{leastsquareadvection:disc}. Then $I_h$ is optimal in the energy norm in
the following sense:
\[
\enorm{\xi - \xi_h} \lesssim \inf_{\delta I\in G_h} \enorm{\xi - \delta I}.
\]
\end{proposition}
\begin{proof}
Observe that for a test function $\delta I \in G_h$,
\begin{align*}
a(\xi-\xi_h, \xi-\xi_h) & = a(\xi-\xi_h + \delta I - \delta I, \xi-\xi_h)\\
                & = a(\xi - \delta I, \xi -\xi_h) + a(\delta I-\xi_h,\xi -\xi_h)\\
                & = a(\xi - \delta I, \xi -\xi_h) + a(\delta I-\xi_h,g_h-g)\\
                & = a(\xi - \delta I, \xi -\xi_h)\\
                & \lesssim \enorm{\xi - \delta I} \enorm{\xi-\xi_h},
\end{align*}
and the result follows.
\end{proof}

This essentially recovers some Galerkin orthogonality whenever the boundary
conditions for the discrete problem is obtained by projection via the conforming
finite element space.

\subsubsection{Examples for fixed $u$}

\begin{figure}[h!]
\centering
\includegraphics[width=0.24\textwidth]{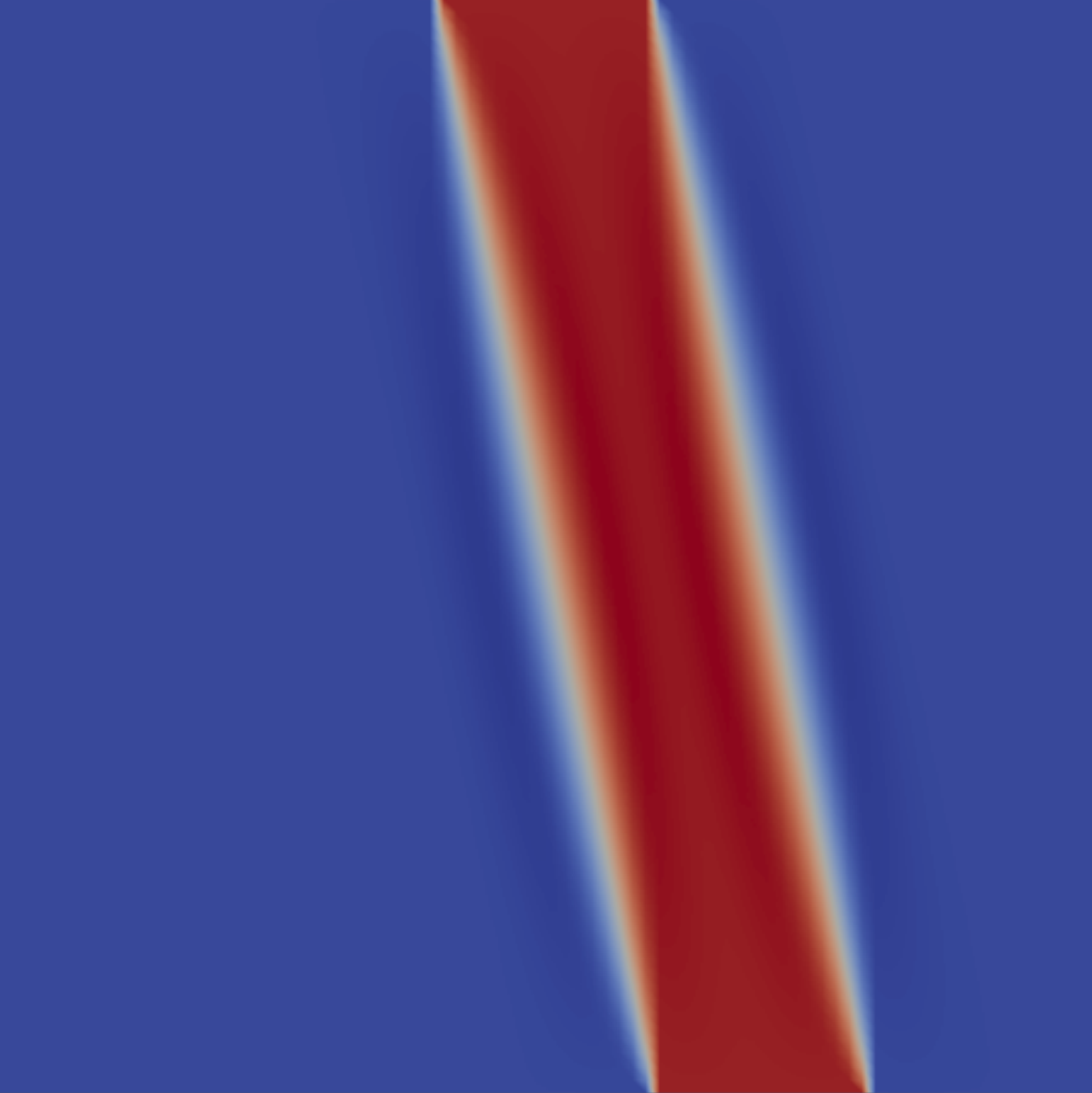}
\includegraphics[width=0.24\textwidth]{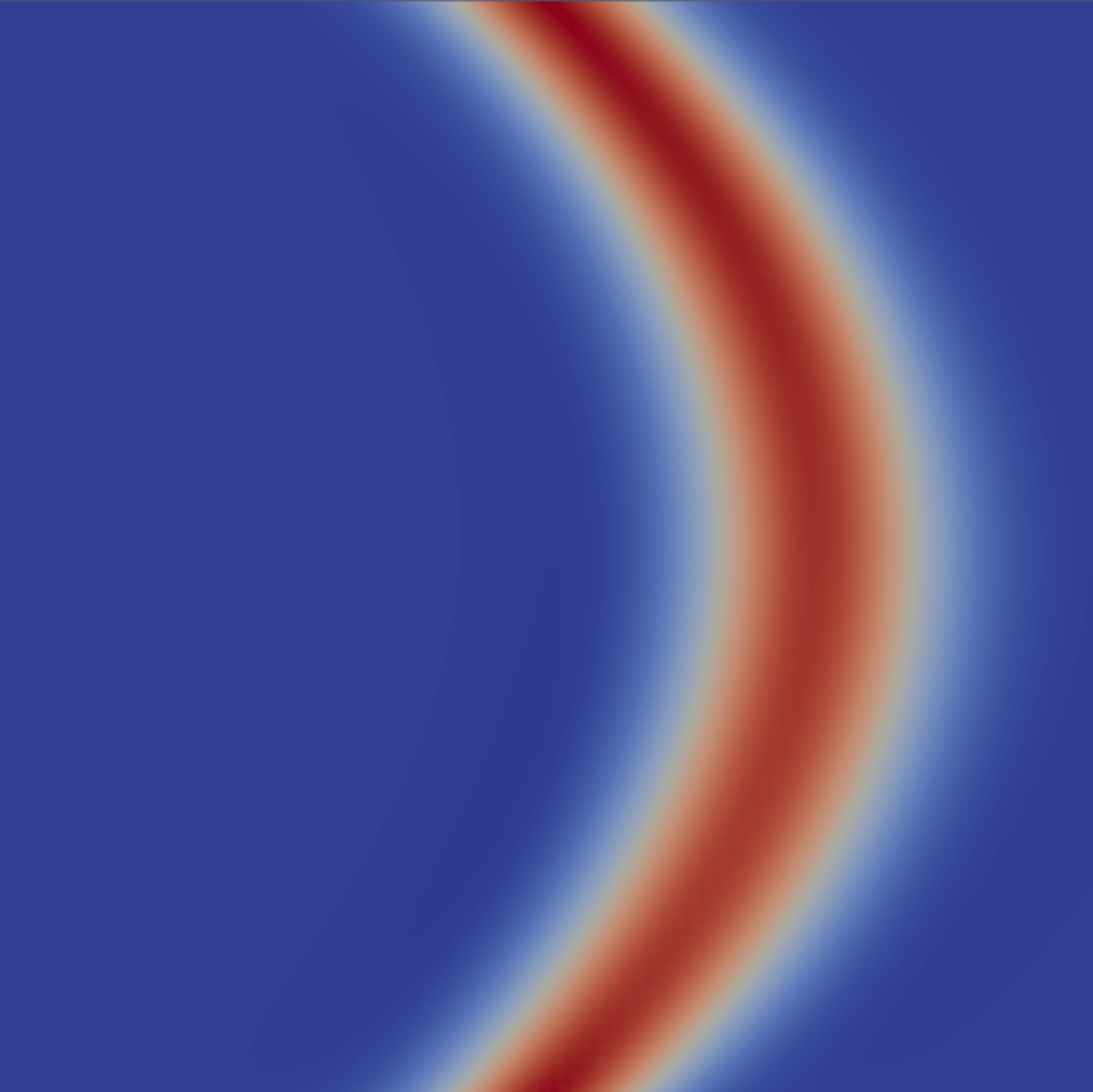}
\includegraphics[width=0.24\textwidth]{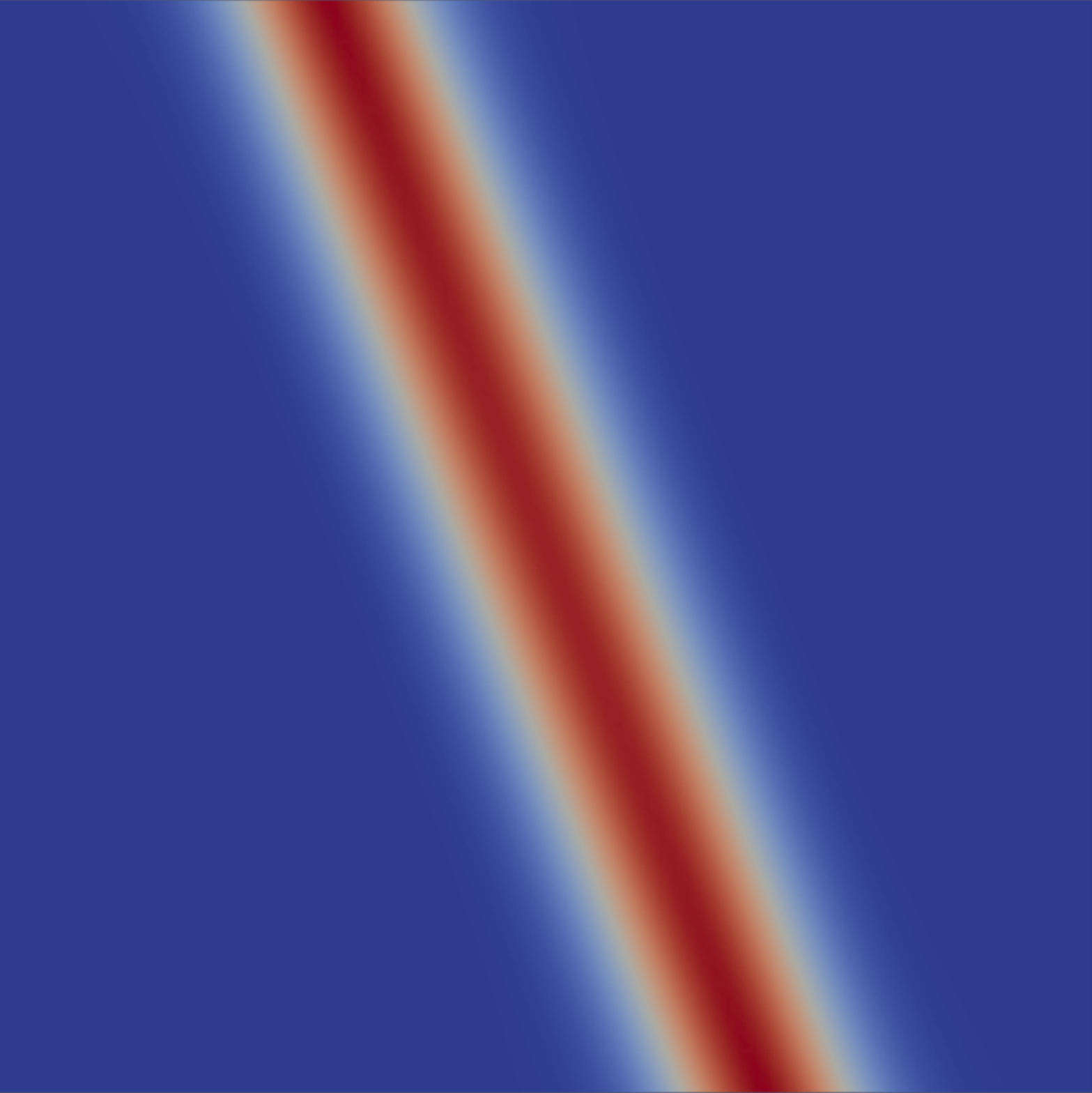}
\includegraphics[width=0.24\textwidth]{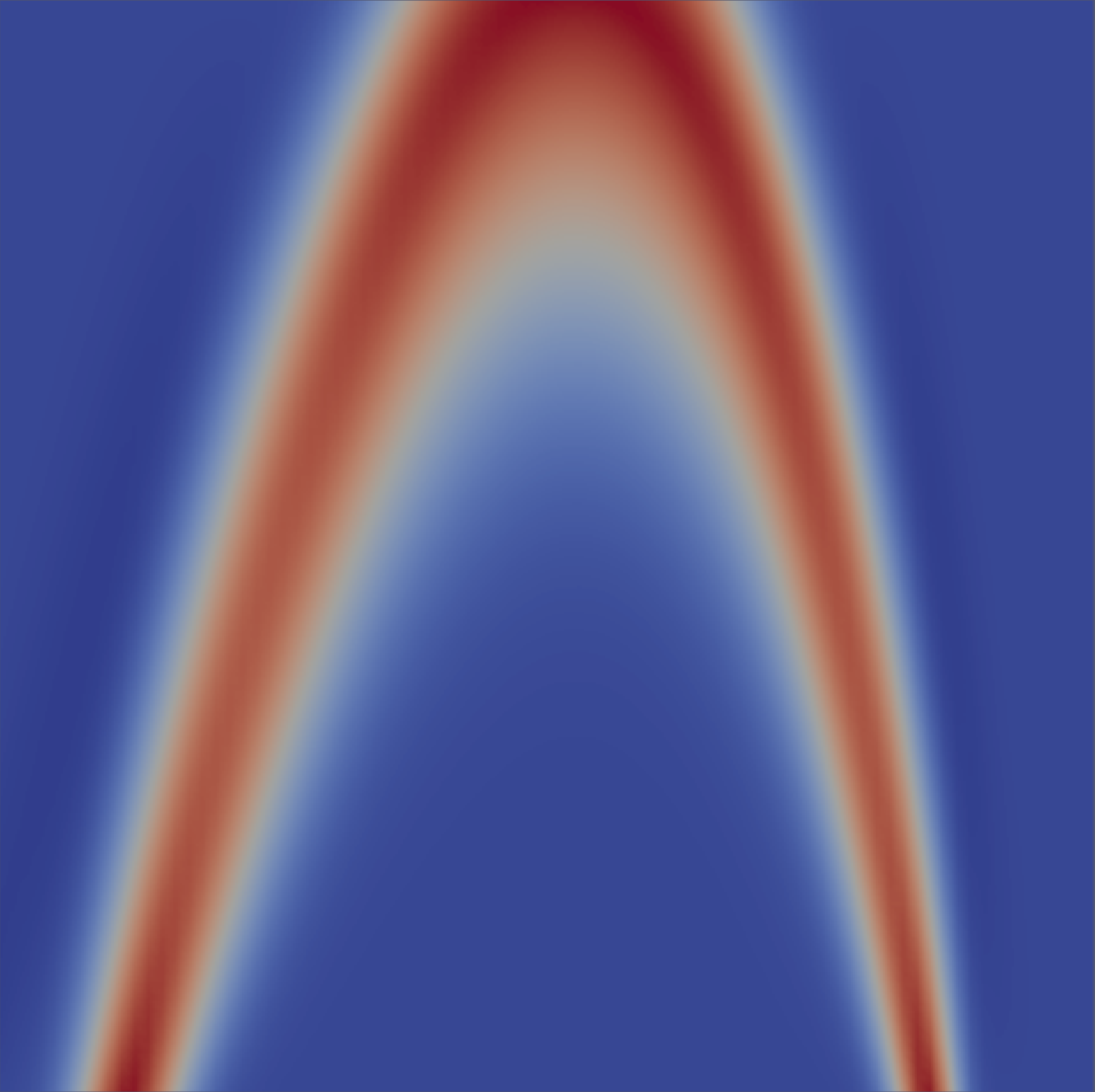}
\caption[Density matches computed from the manufactured solutions to problem
\ref{leastsquareadvection} in Table \ref{tab:cg_manufactured}.]{Density matches
computed from the manufactured solutions to problem \ref{leastsquareadvection}
in Table \ref{tab:cg_manufactured}. Left to right: example 0 through 3.}
\label{fig:cg_manufactured}
\end{figure}

\begin{figure}[h!]
\centering
\includegraphics[width=0.48\textwidth]{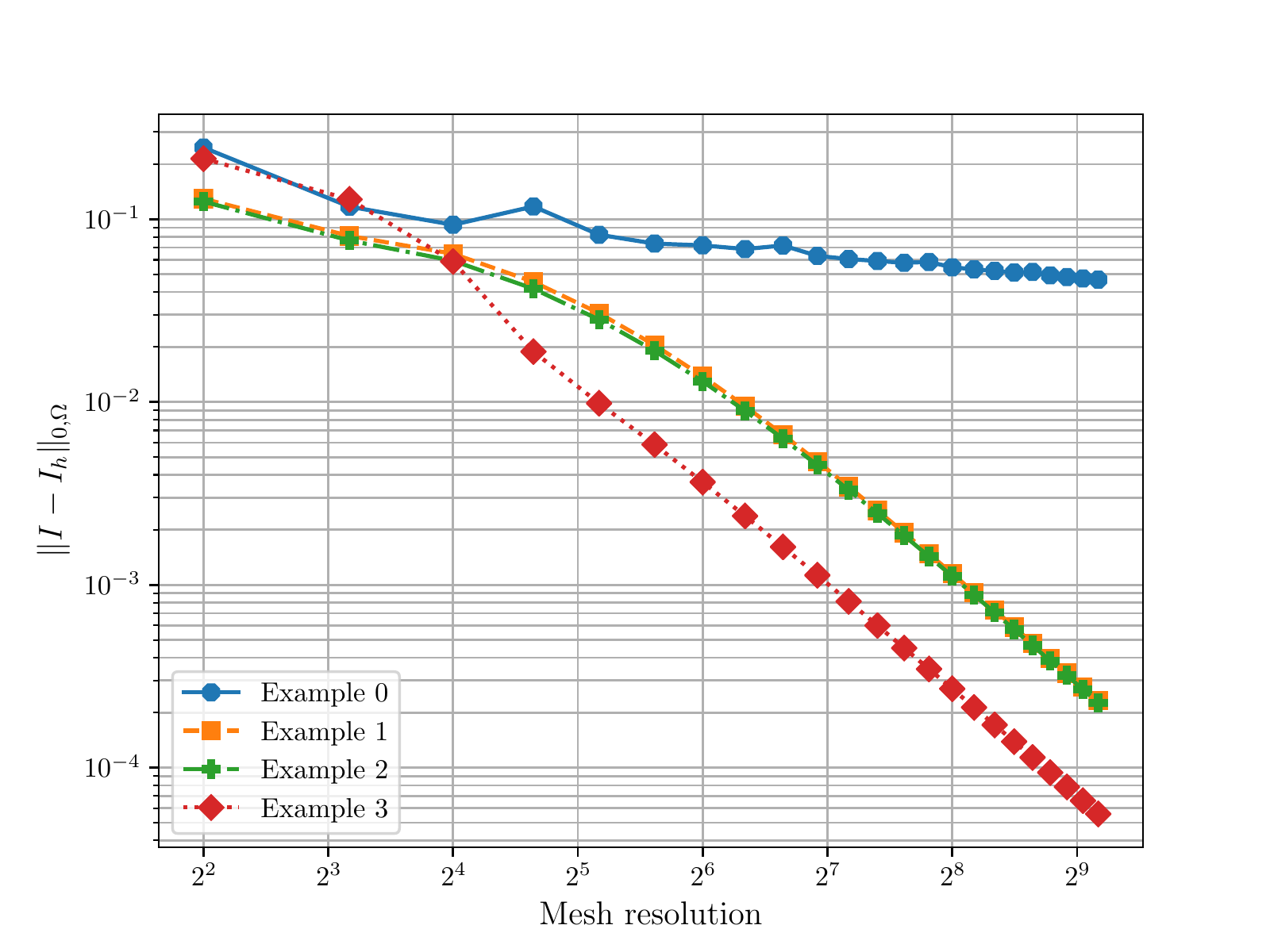}
\includegraphics[width=0.48\textwidth]{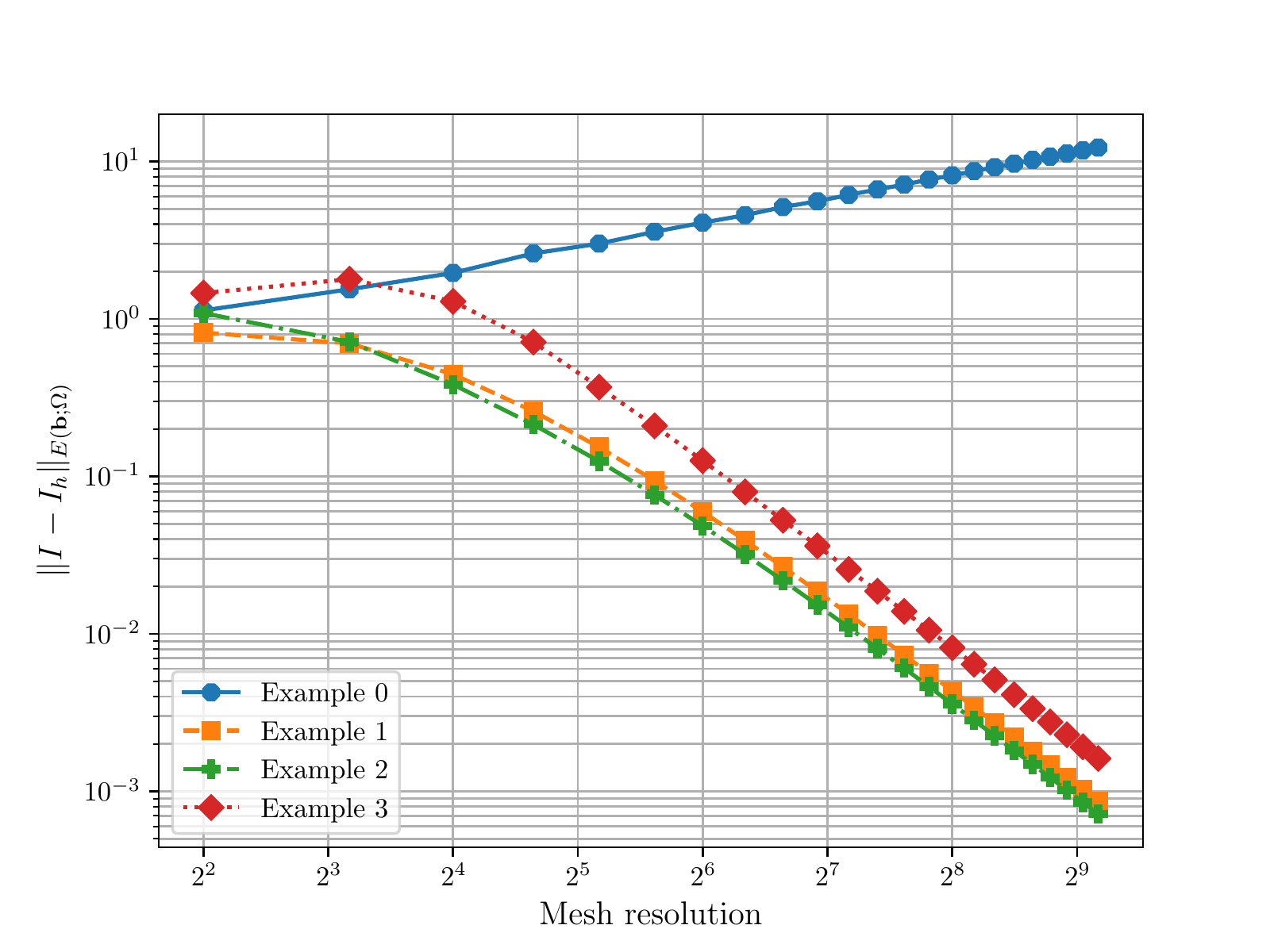}
\caption[Convergence rates for the examples in Table
\ref{tab:cg_manufactured}.]{Convergence as a function
of mesh resolution $h\inv$ for the examples in Table
\ref{tab:cg_manufactured}. Left: $L^2$ errors. Right $\Ebo$ errors.}
\label{fig:cg_inner}
\end{figure}

We manufacture solutions to problem \ref{leastsquareadvection} and show
convergence using a conforming piecewise linear finite element space. The
results presented here are for a periodic spatial domain $\Gamp$. Table
\ref{tab:cg_manufactured} shows four such manufactured solutions with exact
solution $I$ and velocity field $u$, where $\chi_M$ is the indicator function
i.e. equal to $1$ on $M$ and $0$ otherwise\footnote{These are only defined up to
a linear function which sets the value of the image to cater for periodicity in
the spatial dimension. Further, appropriate forcing terms have been added to
match the exact solutions.}. Figure \ref{fig:cg_manufactured} depicts the
computed matches for each of these manufactured solutions using the software
package Firedrake \cite{rathgeber2016firedrake,zenodo}. The intensities across
the four solutions are normalized to the same colour scale as their values are
not important here and are only presented to develop an intuition about the
$L^2$ and energy convergence rates we can observe in Fig. \ref{fig:cg_inner}.
The latter measures the semi-norm specified in \eqref{energy_norm}.  Appendix
\ref{app:innerconv1d} shows the numerical values of these errors and convergence
rates for these four manufactured solutions. We observe some interesting
phenomena. For example 0 the boundary data is discontinuous along $\Gamma$, so
even with a simple constant velocity field the method cannot recover a positive
convergence rate in the energy norm.  \begin{remark} The results for example 0
are not contradictory to our theoretical findings but rather indicating that a
linear $\mathsf{C}^0(\Omega)$ finite element space poorly approximates the
analytical solution.  \end{remark} Some convergence in the $L^2$ norm is
recovered as this norm is somehow agnostic to the discontinuity in the
derivative. Upon inspection of example 0 in Fig.  \ref{fig:cg_manufactured} we
see that the continuous approximation space fails to resolve the fine
discontinuity along the streamline of $u =-0.2$ and the smoothing effect is
clearly seen on the interior of the domain. Examples 1, 2 and 3 provide a
manufactured solution with smooth data and the finite element space is suitable
for approximating the global smoothness of the solutions. The rates are spurious
before entering an asymptotic regime.  These rates also indicate that some
higher-order convergence may be recovered with more in-depth error analysis. We
highlight that example 2 is essentially identical to example 0 albeit with
smooth data, and we clearly see that this is reflected in the convergence rate
of the method. The works \cite{bochev2001improved,houston2000posteriori} look at \emph{a
posteriori analysis} for scalar hyperbolic and transport problems and it may be
of interest to investigate these further in the context above. For completeness
appendix \ref{app:innerconv2d} shows similar results for $d=2$.  Given this
numerical evidence that the inner problems are well-posed given a velocity
satisfying assumption \ref{ass1}.
\begin{table}[h!]
\centering
\caption{Manufactured solutions to {\rm(\ref{leastsquareadvection})} where $d=1$}
{%
\begin{tabular}{c|c|c}
\hline No. & $I(x,t)$ & $u(x,t)$\\
\hline 0 & $\chi_{[0.6-0.2t<x,x<0.8-0.2t]}$ & $-0.2$\\
\hline 1 & $e^{-100(x-t(1-t)-0.5)^2}$ & $1-2t$\\
\hline 2 & $e^{\frac{-(x-0.4(1-t)-0.3)^2}{125}}$ & $-0.4$\\
\hline 3 & $e^{-25(-4.3x+5.16x^3+0.5+t)^2}$ & $ e^{-(0.5 - x)^2}$\\
\hline 
\end{tabular}
}
\label{tab:cg_manufactured}
\end{table}

\subsection{Outer Problem}\label{outer_problem}

The inner problem equips us with a \emph{forward operator} that maps the
one-parameter family of velocity fields to a solution of the inner problem.
These forward operators motivate a gradient descent method to solve the inverse
problem of finding the optimal velocity $u$ minimizing the metamorphosis
functional \eqref{abstract_metamorphosis_problem:fnl}.\\

Central to the preceding analysis was the space $\Gbo$ which depends
explicitly on the velocity $u$ via $\bbf$. If we wish to explore the space
$\ltwoV$ by means of a gradient descent method perturbing a velocity $u\mapsto u
+ \delta u$ (where $\delta u$ is some appropriate search direction) implies a
change in the Hilbert space $\Gbo$. Since a variation $\delta u$ in the velocity
is arbitrary, it is not straight-forward to characterize the derivative of image
$I\in\Gbo$ since it surely will not occupy $\Gbo$ unless the variation $\delta
u$ vanishes.  A practical solution to this problem is therefore simply to
restrict the space of images to be continuous everywhere on $\Omega$ by using
the conforming discretization \eqref{leastsquareadvection:disc}. The approach we
take here is called a \emph{discretize-then-optimize} strategy, where we replace
the function spaces with finite-dimensional surrogates and \emph{then} apply a
gradient descent method. This is in contrast with an
\emph{optimize-then-discretize} scheme, where an optimality system (sometimes
referred to as the Karush-Kuhn-Tucker system) is derived for an
infinite-dimensional optimisation problem after which the resulting equations
are discretized together. We return to this briefly later on, see also
\cite[Chapter 3]{hinze2008optimization} for details.\\

First we define a semi-discrete version of \eqref{metamorphosis_problem_relaxed}
where the velocity is infinite-dimensional and the solution image occupies the
conforming finite element space. This version of the inverse problem is as
follows:
\begin{subequations}\label{metamorphosis_problem_relaxed3} \begin{align}
\inf_{(u,I)\in \ltwoVh\times G_h}\quad & H = \half\int_0^1\| u \|_{V}^2 \dt +
\sigma^{-2} 
\ltwonorm{\bgrad I^*[u]}^2 \\
\text{subject to} \quad &
I^*[u] \triangleq \arg \inf_{\substack{I\in G_h,\\ I|_{\Gamma_0}=I_0,\\
I|_{\Gamma_1}=I_1}} \| \bgrad I \|^2_\ltwoomega.\label{metamorphis_problem_relaxed:inner3}
\end{align}
\end{subequations}
We write this as a \emph{reduced} problem for the control variable $u$:
\begin{subequations}\label{metamorphosis_problem_relaxed4} \begin{align}
\inf_{u\in V}\quad & H_R = \half\int_0^1\| u \|_{V}^2 \dt +
\sigma^{-2} 
\ltwonorm{\bgrad I^*[u]}^2 \label{metamorphosis_problem_relaxed4:fnl}\\
 \text{subject to} \quad &
I^*[u] \text{ solves problem } \ref{leastsquareadvection}.\label{metamorphis_problem_relaxed:inner4}
\end{align}
\end{subequations}
$I^*[u]$ is the \emph{optimal} source term given $u$ in sense of the inner
problem. First, we note that
\eqref{metamorphosis_problem_relaxed3} and
\eqref{metamorphosis_problem_relaxed4} lead to the same stationary points, so
the reformulation is idempotent. Indeed, setting $\delta H=0$ yields:
\begin{align*}
& \int_0^1 \langle u, \delta u \rangle_{V} \dt = 
- \sigma^{-2}\langle \bgrad I, \delta u \cdot \nabla I \rangle_{0,\Omega},\quad & \forall
\delta u \in \ltwoV,\\
& \langle\bbf \cdot \nabla I, \bbf \cdot \nabla \delta I\rangle_{0,\Omega} = 0,\quad & \forall
\delta I \in G_h,
\end{align*}
with the image at time $t=0,1$ fixed. Since $G_h$ is a
globally continuous piecewise finite-dimensional space, the spatial derivative
is well-defined in a weak sense. Carrying out the same steps for $H_R$ we get:
\begin{align}\label{eq:HR}
& \int_0^1\langle u, \delta u \rangle_{V}\dt  = 
- \sigma^{-2}\langle \bgrad I^*[u], \bgrad \frac{\delta I^*}{\delta u}[u]\delta
u+ \delta u\cdot\nabla I^*[u] \rangle_{0,\Omega},\quad & \forall
\delta u \in \ltwoV,
\end{align}
where $\frac{\delta I^*}{\delta u}[u] \delta u $ is the sensitivity of $I^*[u]$
at $u$ with respect to $u$. We notice that the reduced nature of the
reformulation \eqref{metamorphosis_problem_relaxed4} introduces a new term in
the equation for $u$. However, since $\frac{\delta I^*}{\delta u}[u] \delta u
\in G_h$ (since the variational derivative is taken in sense of Fr\'echet),
then:
\[
\langle\bgrad I^*[u],\bgrad\frac{\delta I^*}{\delta u}[u]\rangle_{0,\Omega} = 0,
\]
since $I^*[u]$ solves problem \ref{leastsquareadvection}, which is a convex
problem, so the equivalence follows.\\

Equipped with a reduced problem we now discretise the space of the velocity
field.  The norm in \eqref{Vnorm} is equivalent to an $H^3$ norm in space in
line with assumption \ref{ass1}. Ideally, this property must be preserved under
discretization, \emph{independently} of the mesh refinement parameter. An $H^3$
conforming finite element space is globally $\mathsf{C}^2$ and an implementation
is, to the best of the authors' knowledge, not available. In this paper we focus
our mathematical analysis on the inner problem and in this section concentrate
our work on developing and implementing a gradient descent scheme on the
velocity field. Deferring convergence analysis pertaining to the velocity field
to future work we therefore discretize the velocity field by continuous
piecewise affine functions. In practice we use the iterated nature of the
operator $L$ in \eqref{Vnorm} and solve a modified version of
\eqref{metamorphosis_problem_relaxed4} as follows. Using the space-time bilinear
form:
\[
a(w,v) = \langle w, v \rangle_{0,\Omega} + \alpha^2\langle \nabla w,
\nabla v \rangle_{0,\Omega},
\]
this version given by \eqref{metamorphosis_problem_relaxed5} below:
\begin{subequations}\label{metamorphosis_problem_relaxed5} \begin{align}
\inf_{v_h^0\in \pone}\quad & H_R = \half \| v_h^0 \|_{0,\Omega}^2 +
\sigma^{-2} 
\ltwonorm{\bgrad I_h^*[u_h]}^2 \label{metamorphosis_problem_relaxed4:fnl}\\
 \text{subject to} \quad &
G_h \ni I_h^*[u_h] \text{ solves problem }
\ref{leastsquareadvection:disc},\label{metamorphis_problem_relaxed:inner4}\\
&
a(v^1_h, w) = \langle v_h^0, w\rangle_{0,\Omega}, \quad \forall
w \in \pone,\\
& a(v_h^2, w) = \langle v_h^1, w
\rangle_{0,\Omega}, \quad \forall
w \in \pone,\\
& a(u_h, w) = \langle v_h^2, w
\rangle_{0,\Omega}, \quad \forall
w \in \pone,
\end{align}
\end{subequations}
where $\pone$ is a $\mathsf{C}^0$ conforming finite element space consisting of
piecewise affine functions defined over the mesh $\Omega$. Although this
precludes $u_h$ from having a \emph{mesh-independent} Lipschitz constant,
\eqref{metamorphosis_problem_relaxed5} forms the basis for a useful numerical
method. The next section discusses how this can be realized numerically using
the Firedrake software package
\cite{rathgeber2016firedrake} as well as how knowledge of the
existence of the gradient from the previous section motivates an adjoint-based
automatic differentiation algorithm.

\subsection{Numerical Examples}

\emph{Automatic differentiation} is a software abstraction allowing for the
numerical evaluation of the Jacobian of expressions that are given as functions
in a programming language, see \cite{naumann2012art}. These essentially
articulate chain and product rules for computer programmes. The implementation
of such an abstraction permits a user to evaluate gradients of for instance
complex composite functions for use in an optimisation loop. For instance,
formally speaking, the derivative of the functional in
\eqref{metamorphosis_problem_relaxed4} with respect to  the velocity requires
some sense to be made about the expression $\frac{\delta z^*[u]}{\delta u}$. In
the PDE-constrained optimisation community this expression is called the
\emph{adjoint} and it is precisely the machinery of automatic differentiation
that (upon suitable discretization) allows for its numerical evaluation. This
machinery is indeed necessary as no explicit gradient can be evaluated from
$z^*[u]$. Firedrake supports an automatic differentiation engine (see
\cite{mitusch2019dolfin}) so with only a few lines of code expresses the
optimisation problem \eqref{metamorphosis_problem_relaxed4} which can then be
solved by e.g.  a Broyden-Fletcher-Goldfarb-Shanno (BFGS) algorithm
\cite{nocedal2006numerical} (we used the Rapid Optimization Library (ROL), see
\texttt{trilinos.github.io}). Here we provide but a few numerical examples to
support our work in section \ref{inner:lsq}.\\

\begin{figure}
  \centering
  \subfigure{\includegraphics[width=.24\linewidth]{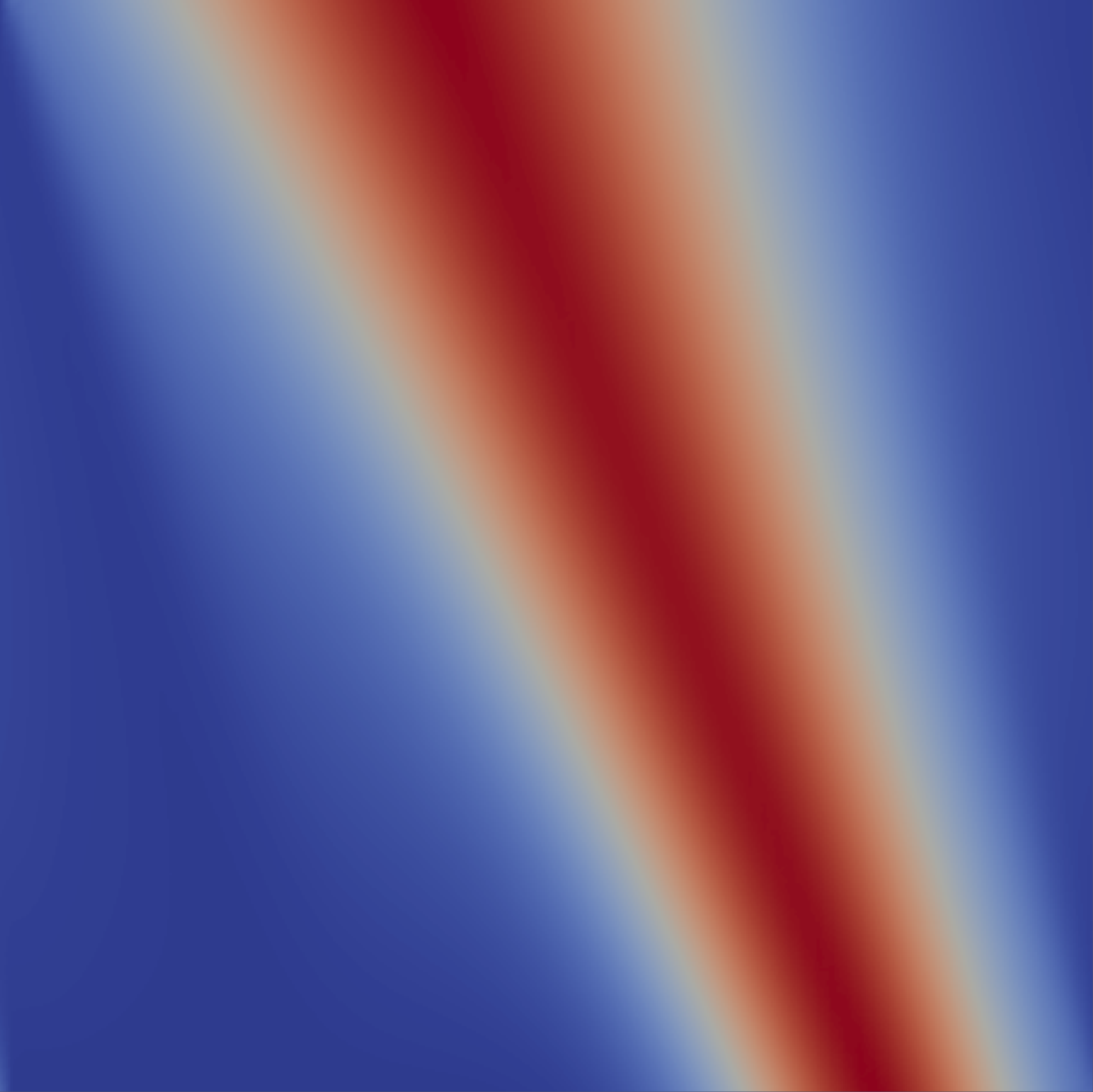}}
  \subfigure{\includegraphics[width=.24\linewidth]{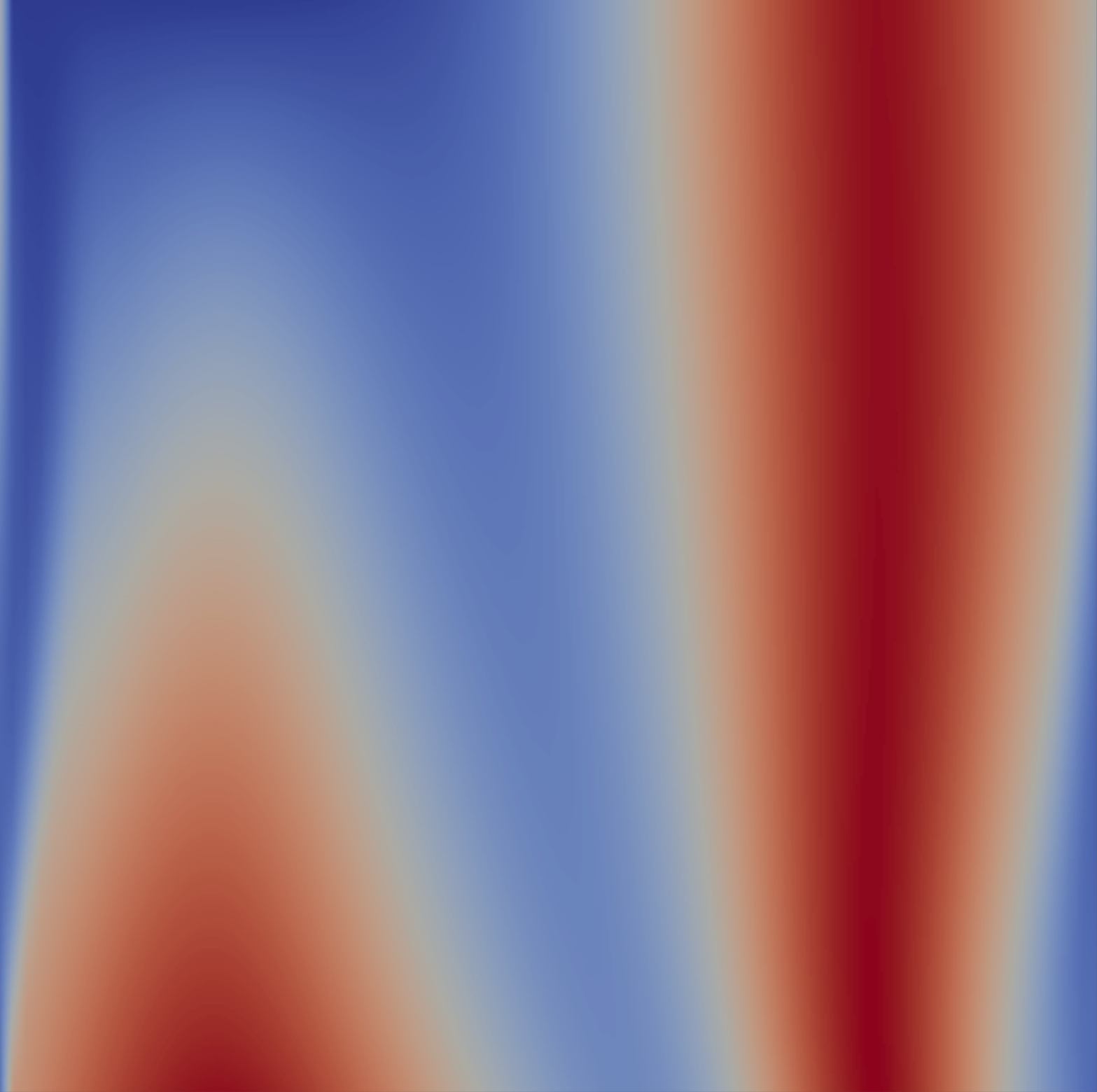}}
  \subfigure{\includegraphics[width=.24\linewidth]{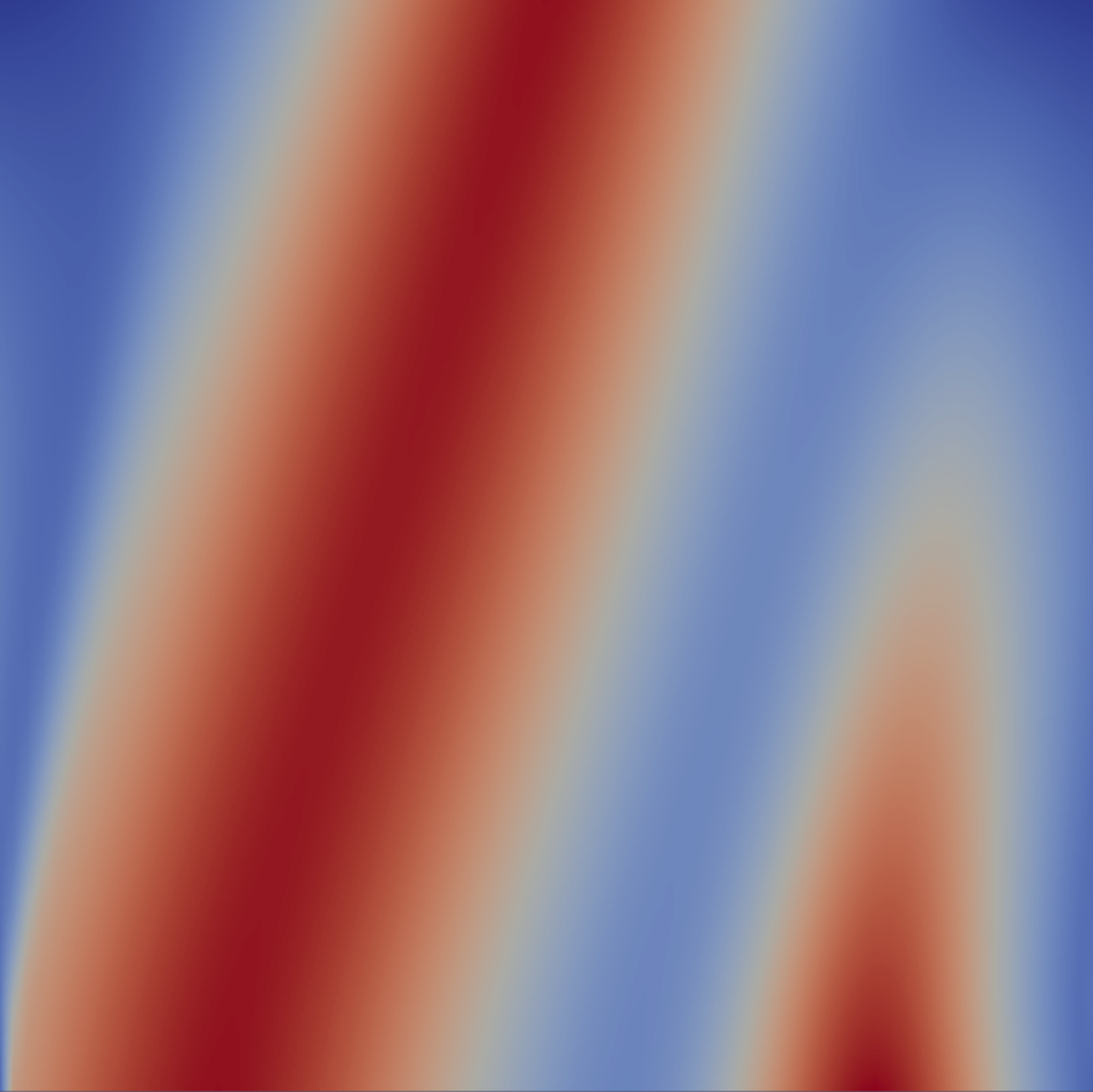}}
  \subfigure{\includegraphics[width=.24\linewidth]{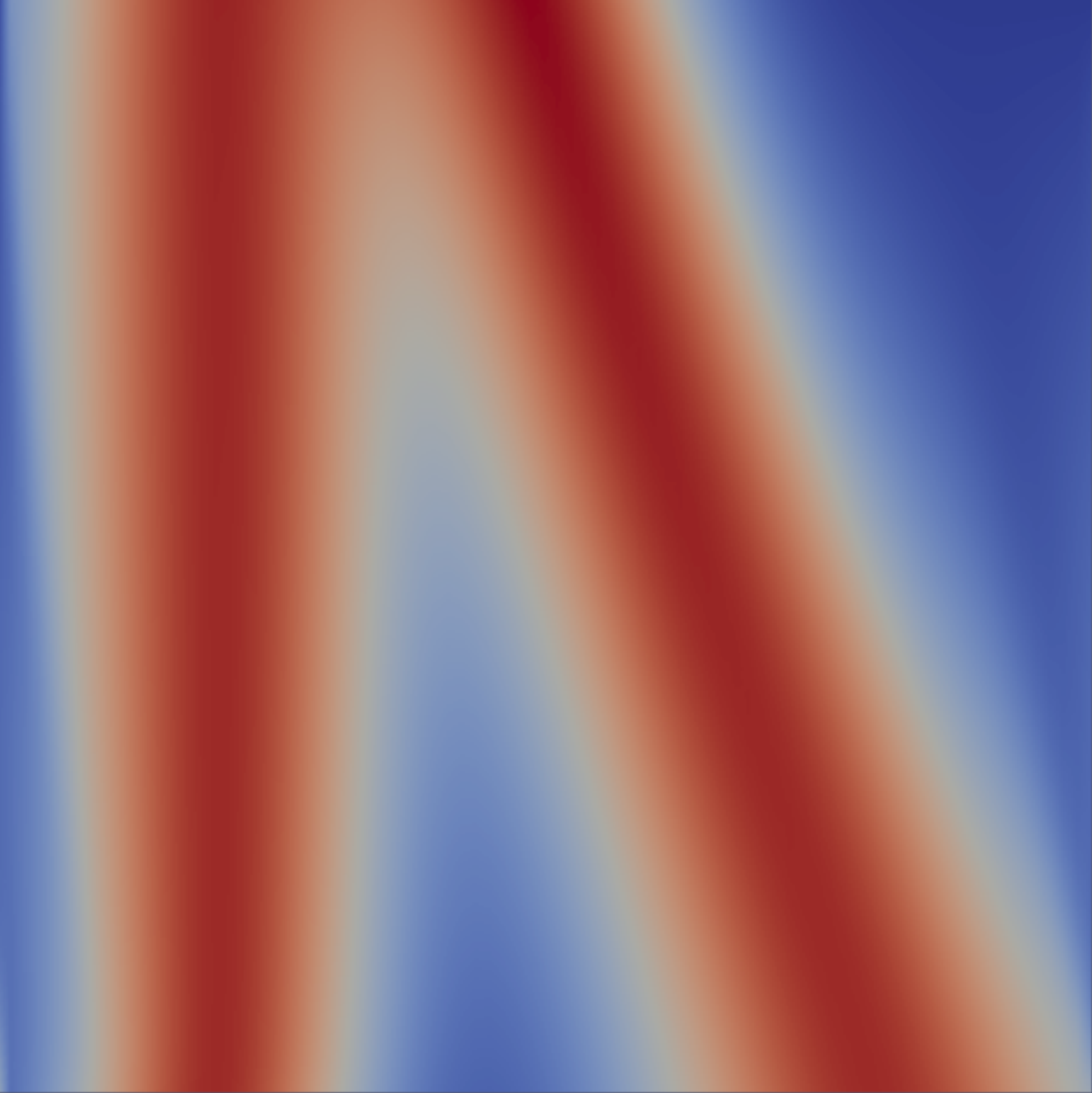}}
  \subfigure{\includegraphics[width=.24\linewidth]{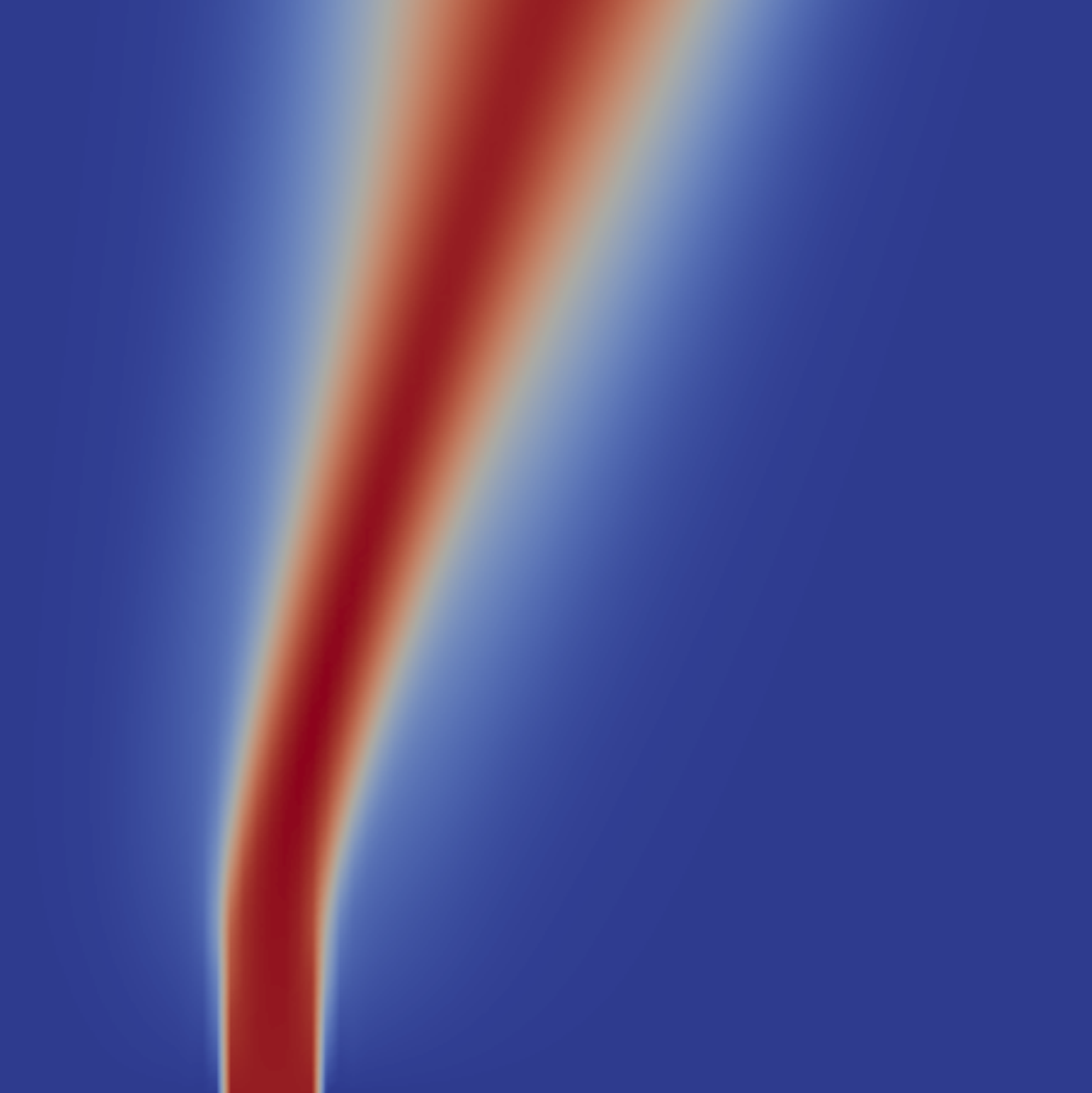}}
  \subfigure{\includegraphics[width=.24\linewidth]{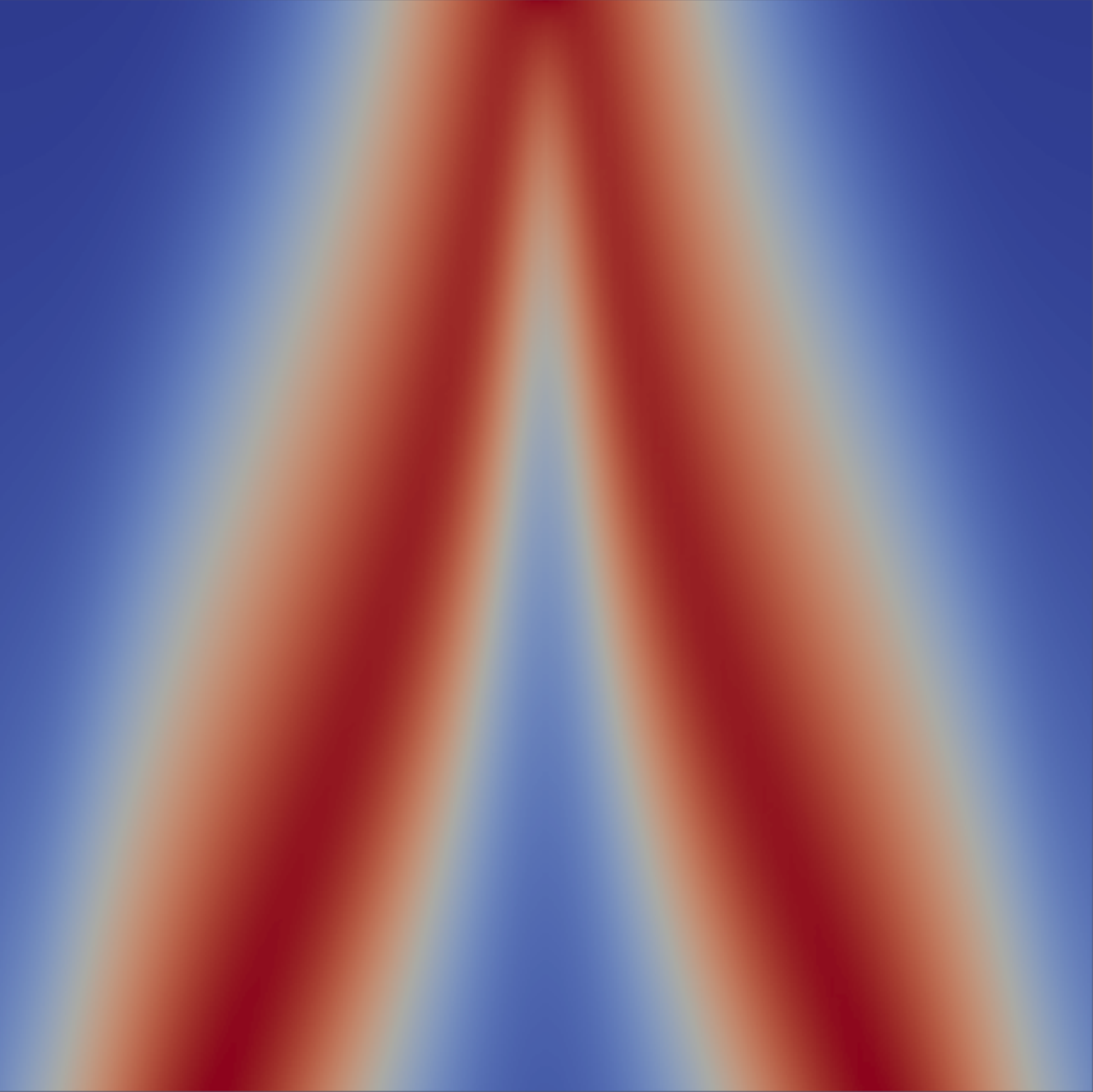}}
  \subfigure{\includegraphics[width=.24\linewidth]{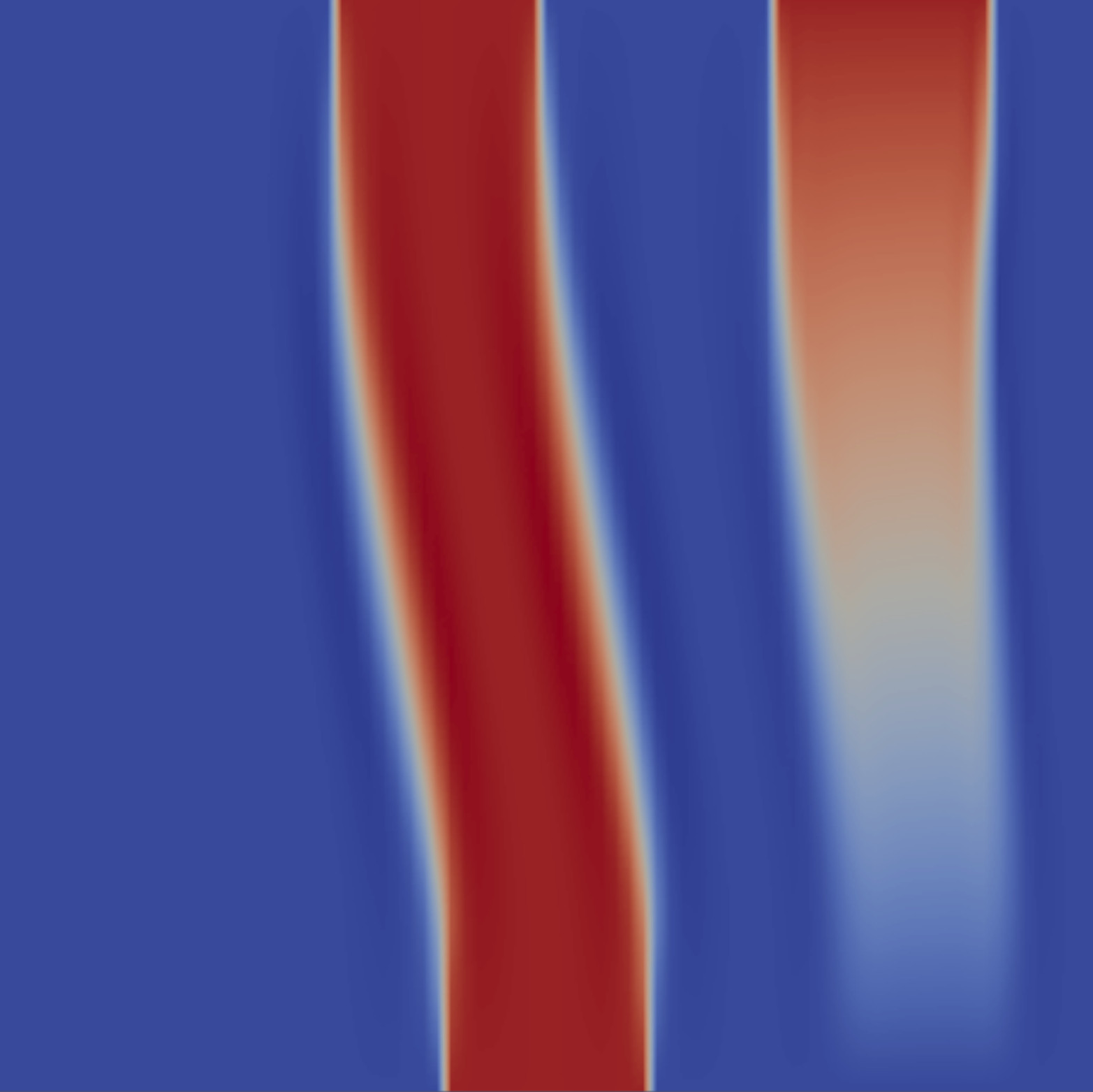}}
  \subfigure{\includegraphics[width=.24\linewidth]{imgs/results1d_c=1/test3.pdf}}
\caption[Results for the outer problem for various template and target pairs in
one dimension using $\sigma^{-2}=1$.]{Results for the outer problem for
various template and target pairs in one dimension. The vertical dimension is
time and the horizontal dimension is space, so the \emph{north} and \emph{south}
faces of the domain represent the target and template, respectively. Here we
have taken $\sigma^{-2}=1$.}
\label{fig:adjoint1d_c=1}
\end{figure}

\begin{figure}
  \centering
  \subfigure{\includegraphics[width=.24\linewidth]{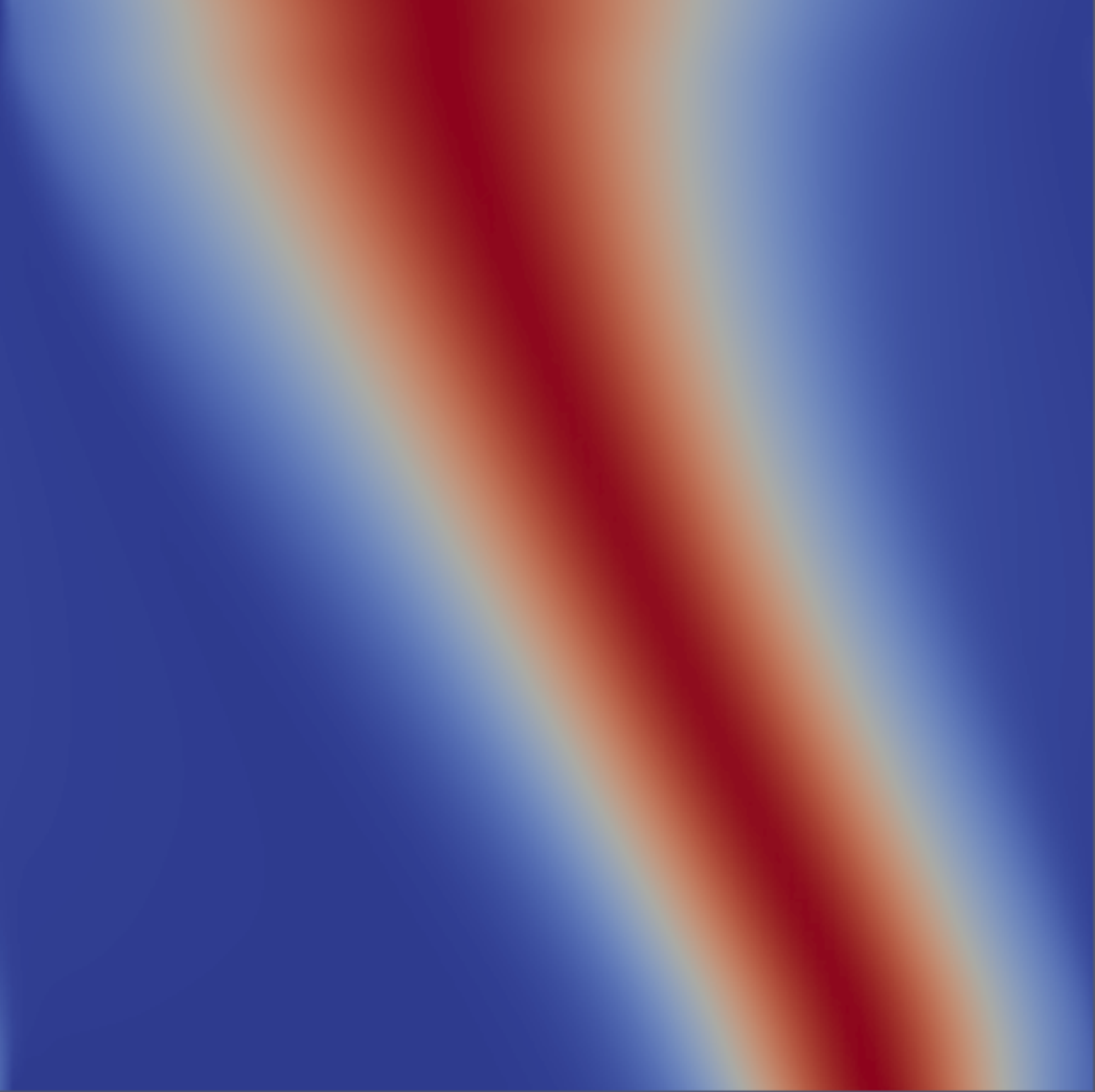}}
  \subfigure{\includegraphics[width=.24\linewidth]{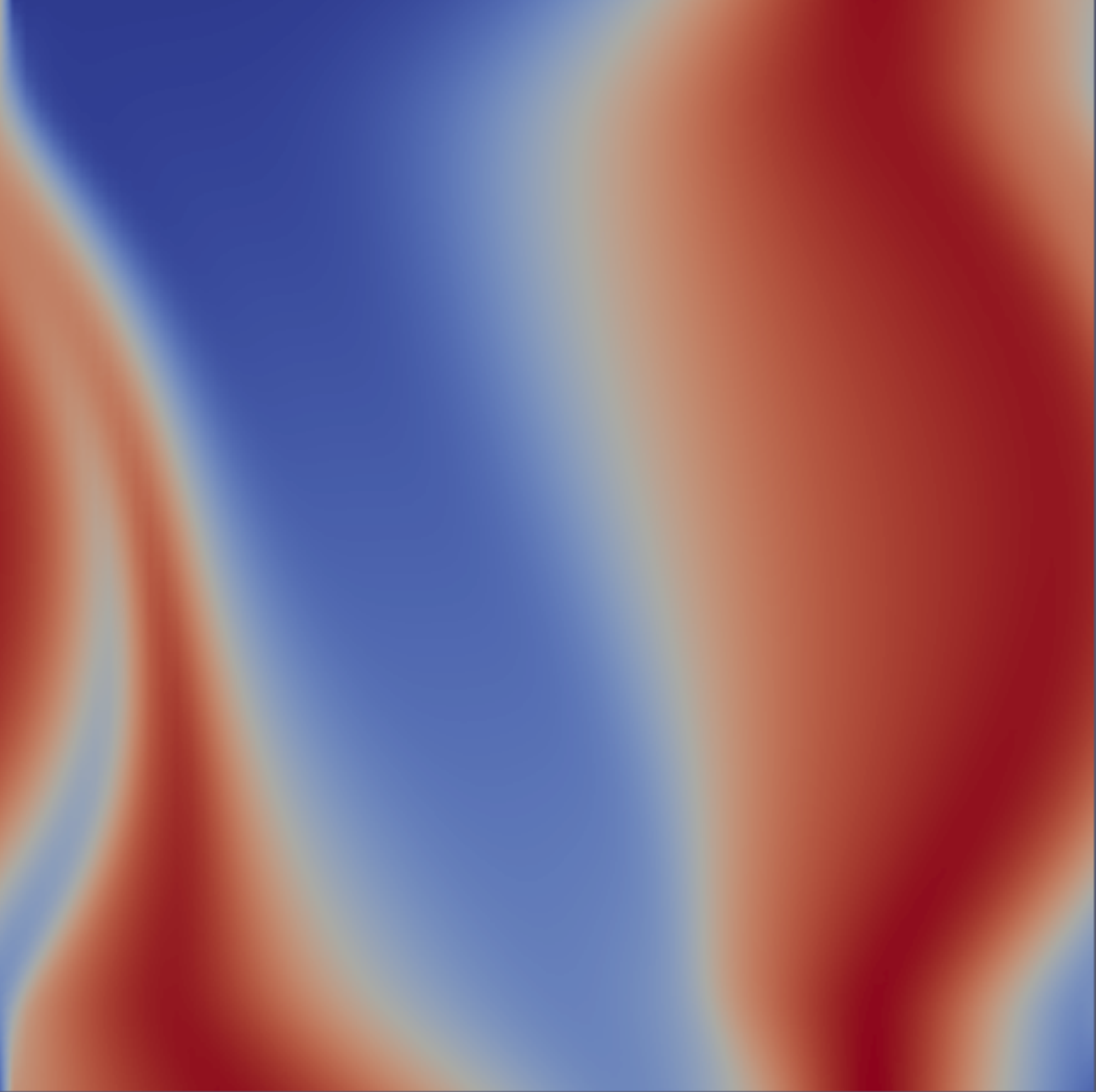}}
  \subfigure{\includegraphics[width=.24\linewidth]{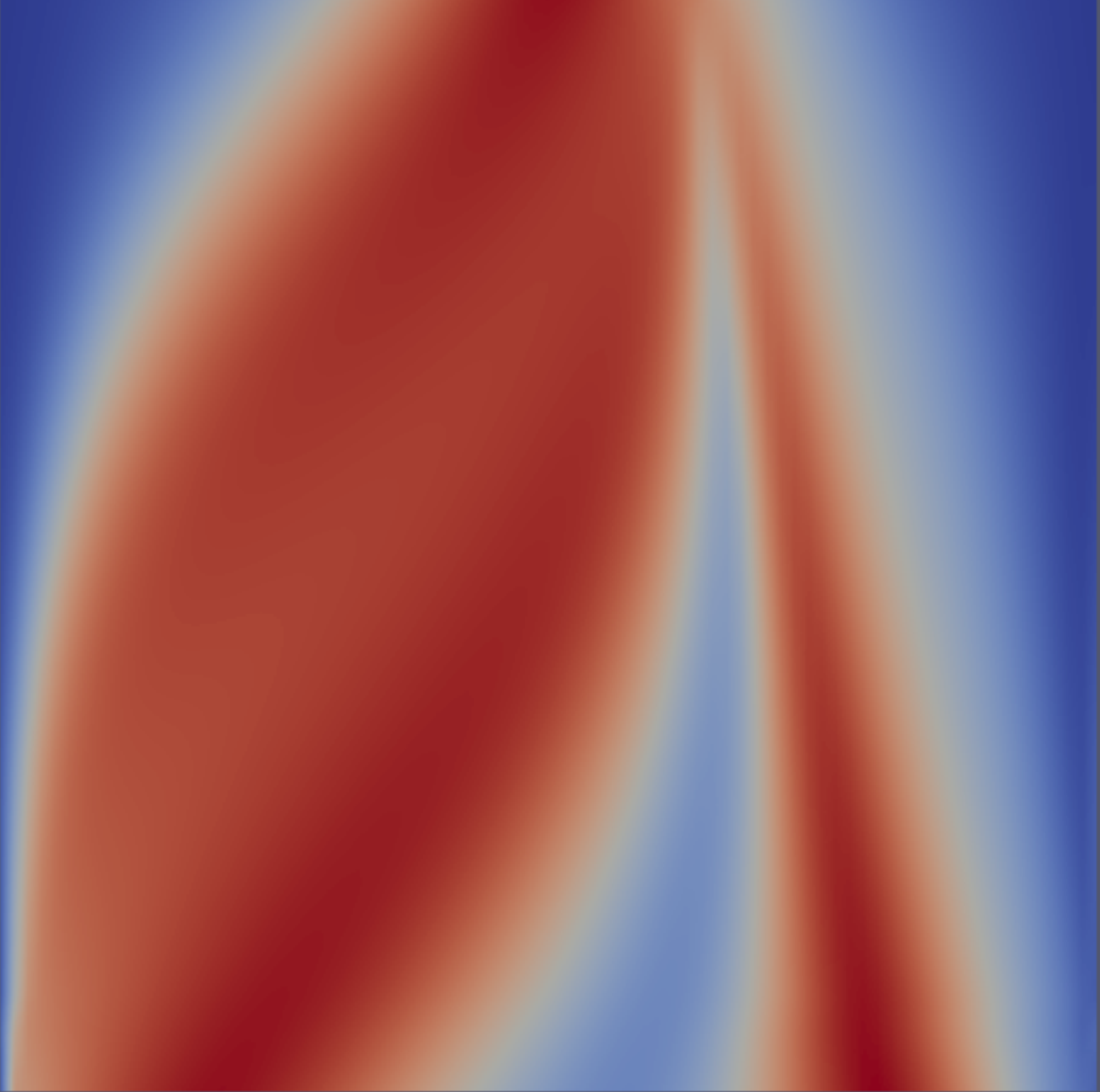}}
  \subfigure{\includegraphics[width=.24\linewidth]{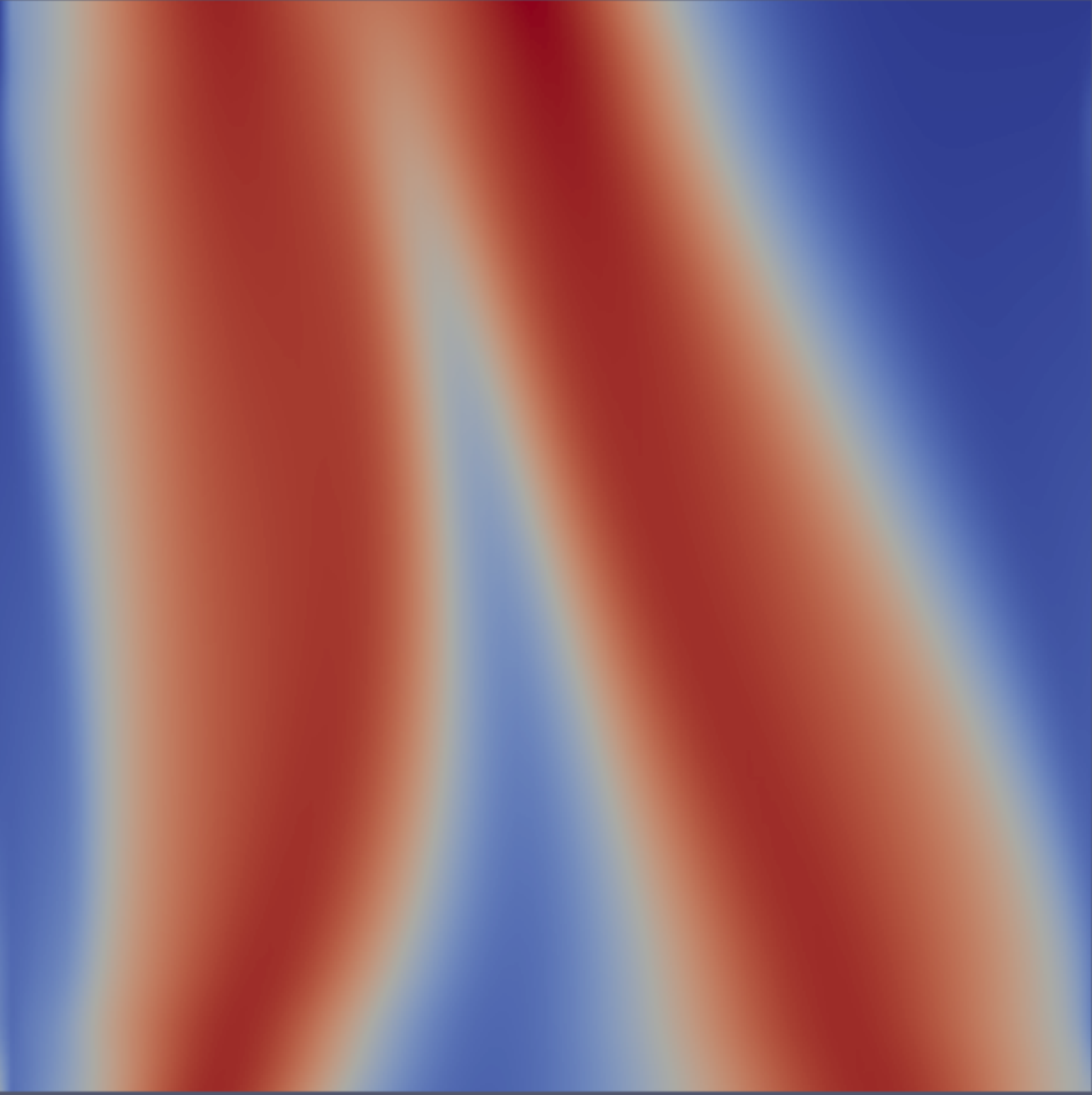}}
  \subfigure{\includegraphics[width=.24\linewidth]{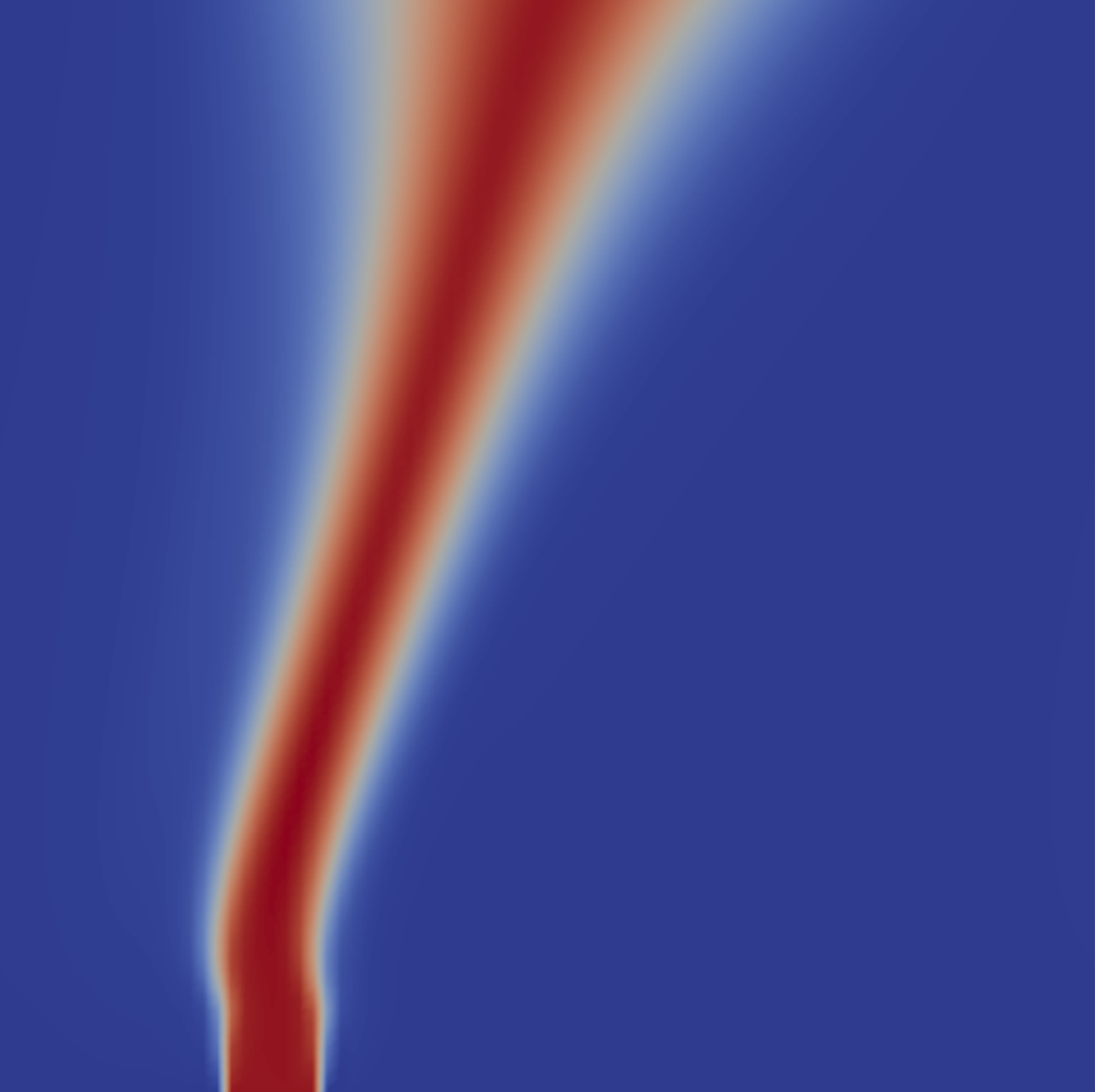}}
  \subfigure{\includegraphics[width=.24\linewidth]{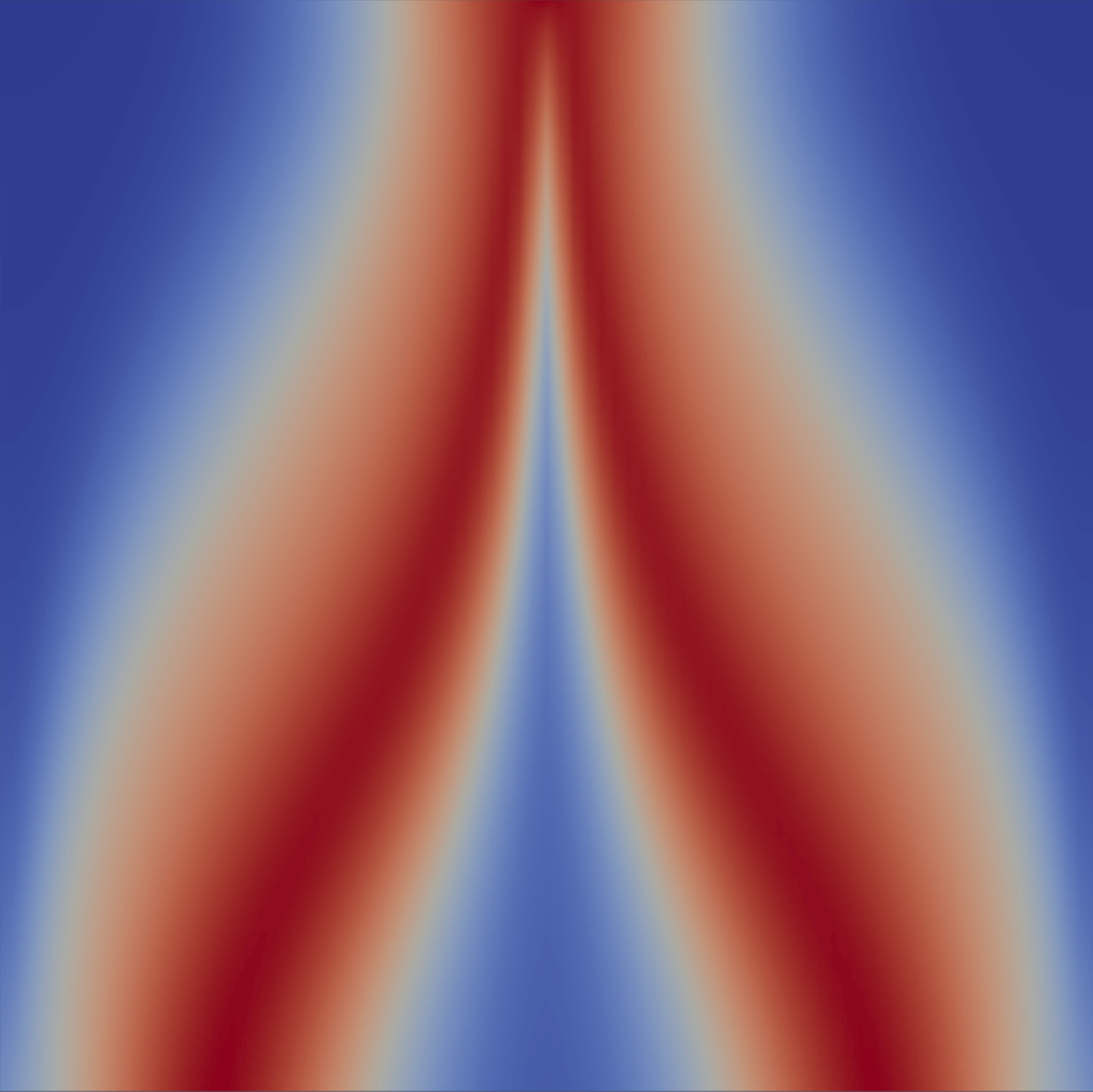}}
  \subfigure{\includegraphics[width=.24\linewidth]{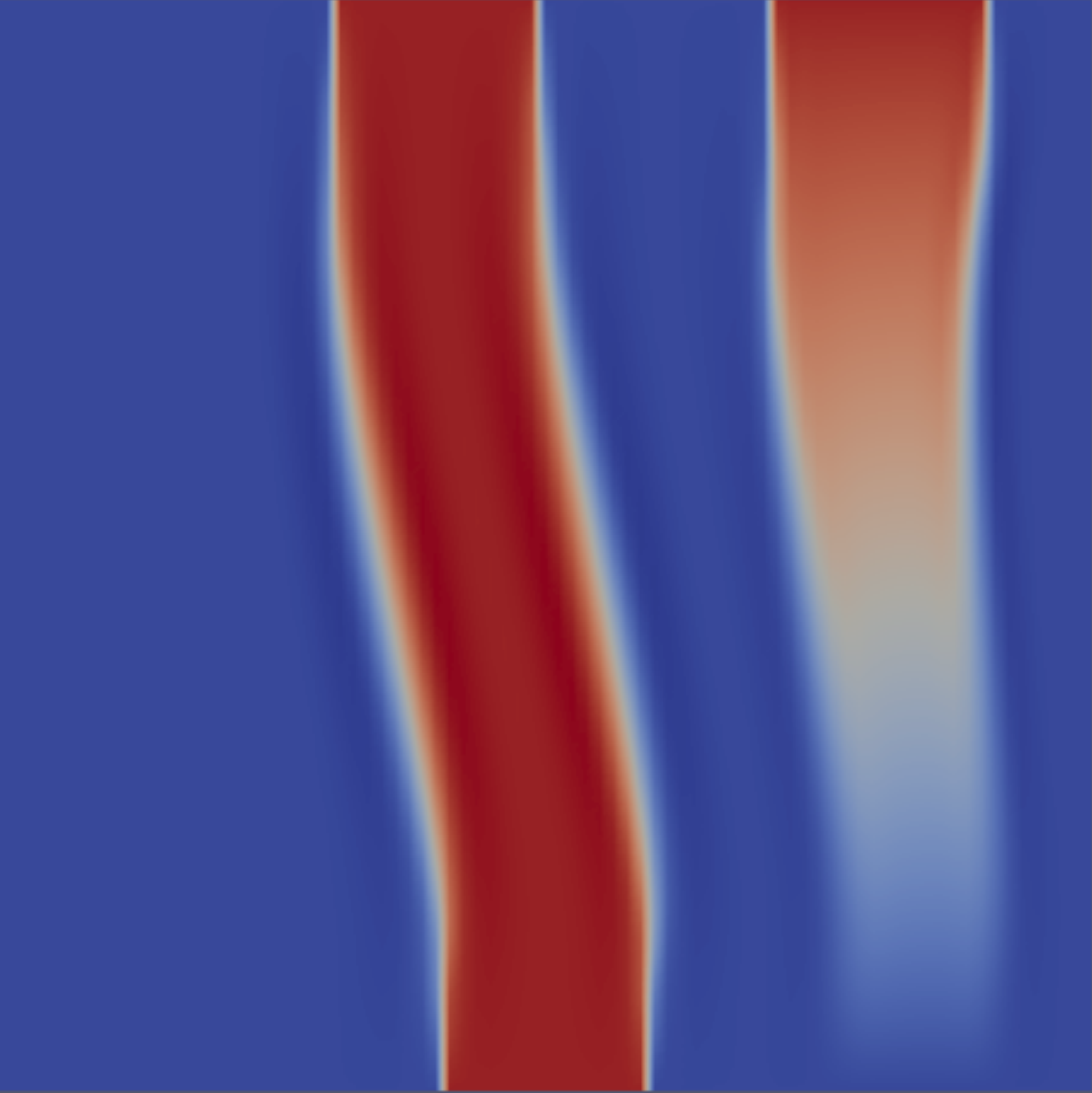}}
  \subfigure{\includegraphics[width=.24\linewidth]{imgs/results1d/test3.pdf}}
\caption[Results for the outer problem for various template and target pairs in
one dimension using $\sigma^{-2}=10^5$.]{Results for the outer problem for various
template and target pairs in one dimension with the same template and target
densities as in Fig. \ref{fig:adjoint1d}. Here we have taken
$\sigma^{-2}=10^5$ thereby forcing the algorithm to seek more diffeomorphic geodesics.}
\label{fig:adjoint1d}
\end{figure}

Figures \ref{fig:adjoint1d_c=1}--\ref{fig:adjoint1d} show the image (or
density) $I$ results of a gradient descent applied to the problem in
\eqref{metamorphosis_problem_relaxed4} using Firedrake's automatic
differentiation abstraction and the conforming discretization from section
\ref{leastsquareadvection}. In these figures, $h=0.0125$. Both the image and the
velocity field is discretized using a globally continuous piecewise linear
finite element space on a spatially periodic domain. These figures show that our
method is able to compute meaningful matches between different densities. In
Fig. \ref{fig:adjoint1d_c=1}.A and E we see that transport solutions are
recovered, even in the presence of a template with a sharp spatial gradient. For
the other cases we have chosen densities which can be matched without pure
transport. In some instances, i.e. F and D, the gradient descent finds intuitive
matches between the data, while B, C and G struggle to recover the same visual
effect.  For these three configurations we find that advection is recovered for
some parts of the image but we do not see a \emph{merging} effect as before. For
instance, in Fig. \ref{fig:adjoint1d_c=1}.C it is more economical to seek an
approximately constant velocity field to match the left-most mode of the
template while causing the right-most mode to vanish.\\

We also comment on the role of $\sigma$ which we recall acts as a penalisation
parameter for advection term. In Fig. \ref{fig:adjoint1d_c=1} the two terms in
\eqref{metamorphosis_problem_relaxed4} are weighted equally.  Figure
\ref{fig:adjoint1d} uses the same data as in Fig. \ref{fig:adjoint1d_c=1},
albeit with $\sigma^{-2}=10^5$. We observe similar results as before with the
exception of B and C. Adjusting the value of the penalty parameter allows the
gradient method to recover, with moderate success, some of the density merging
behaviour observed previously. With a more complete theoretical understanding of
the space-time gradient descent method, and more expertise in numerical
optimisation, we can begin to carry out more careful parameter studies for
$\sigma$. This is deferred to future work.

\begin{figure}[h!]
\begin{minipage}{\textwidth}
\centering
\includegraphics[width=.19\textwidth]{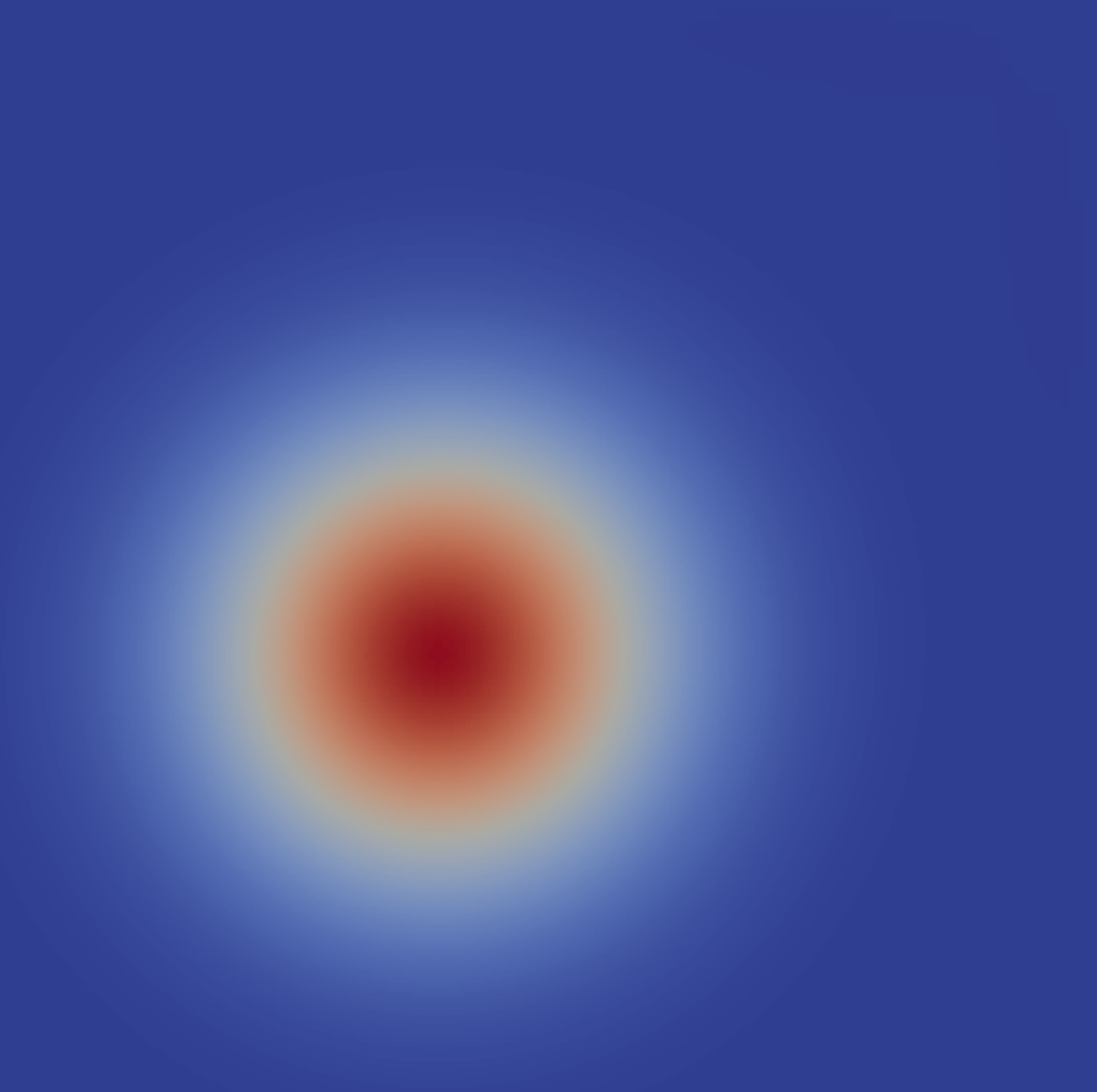}
\includegraphics[width=.19\textwidth]{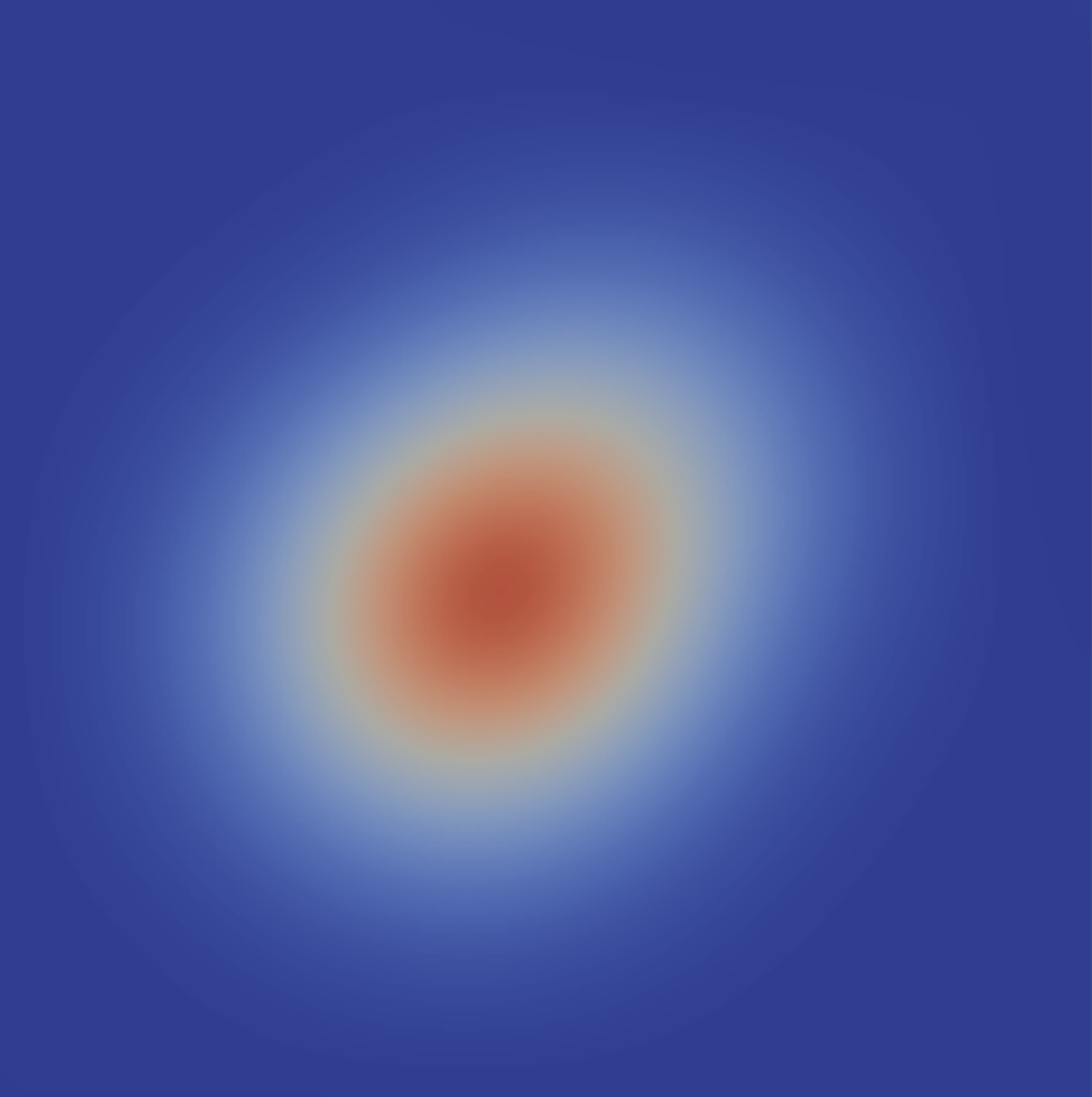}
\includegraphics[width=.19\textwidth]{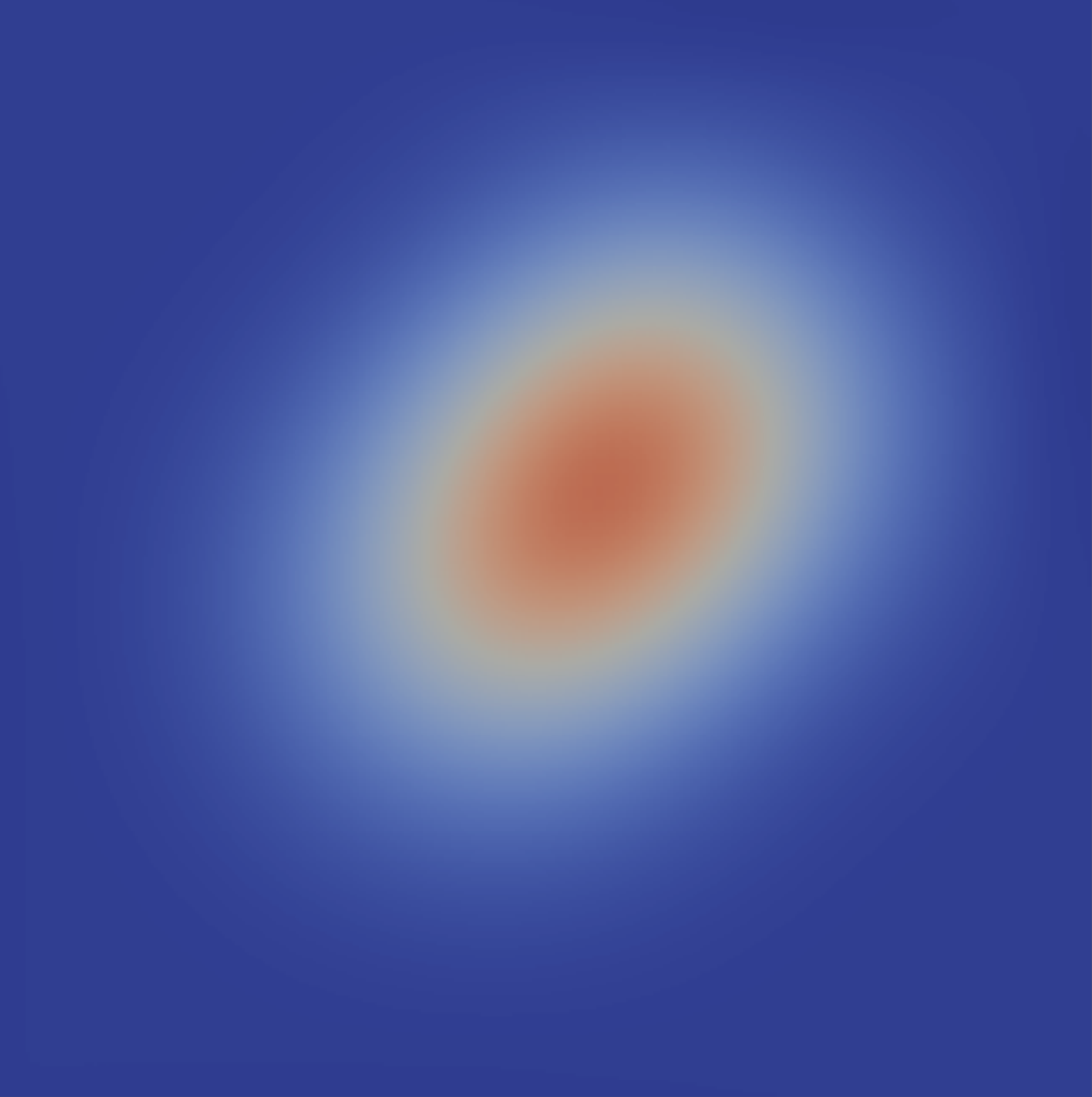}
\includegraphics[width=.19\textwidth]{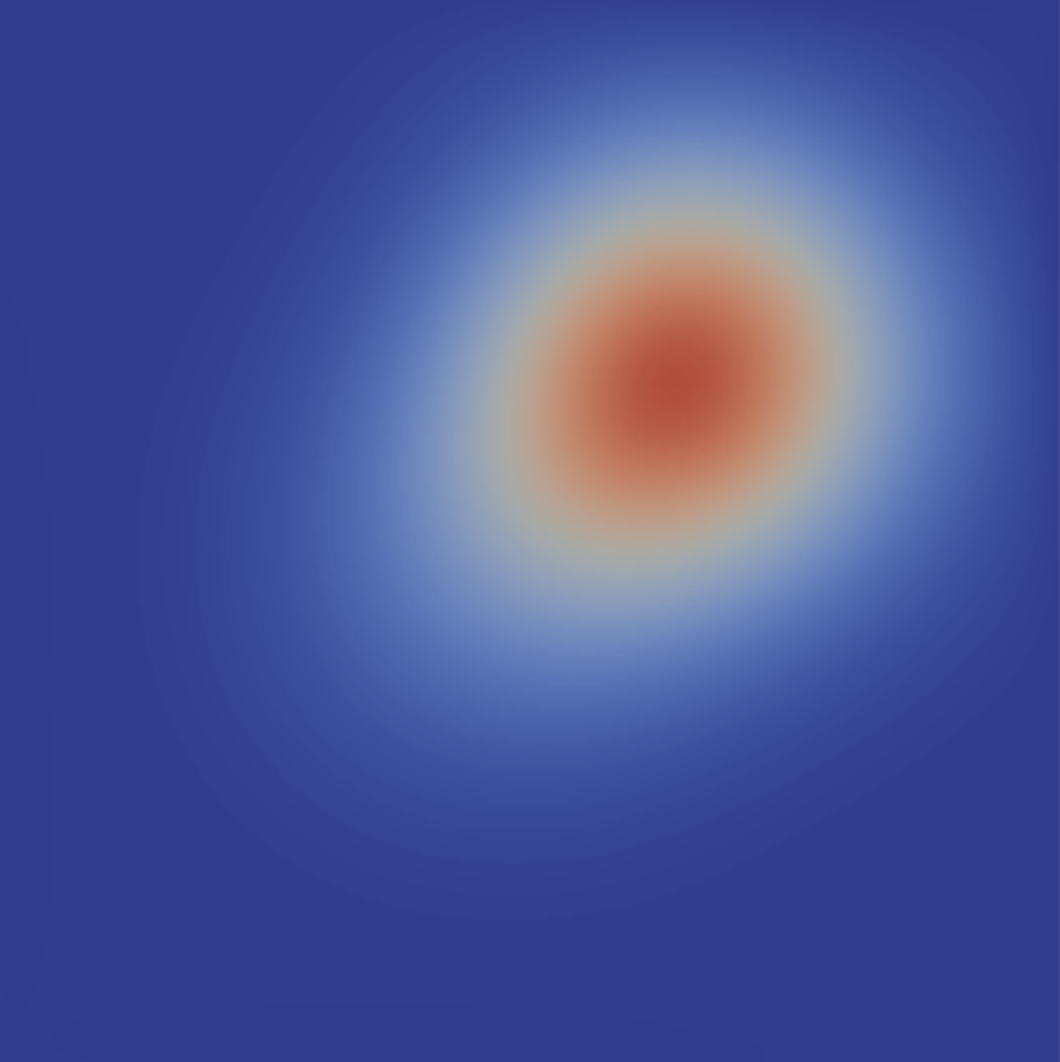}
\includegraphics[width=.19\textwidth]{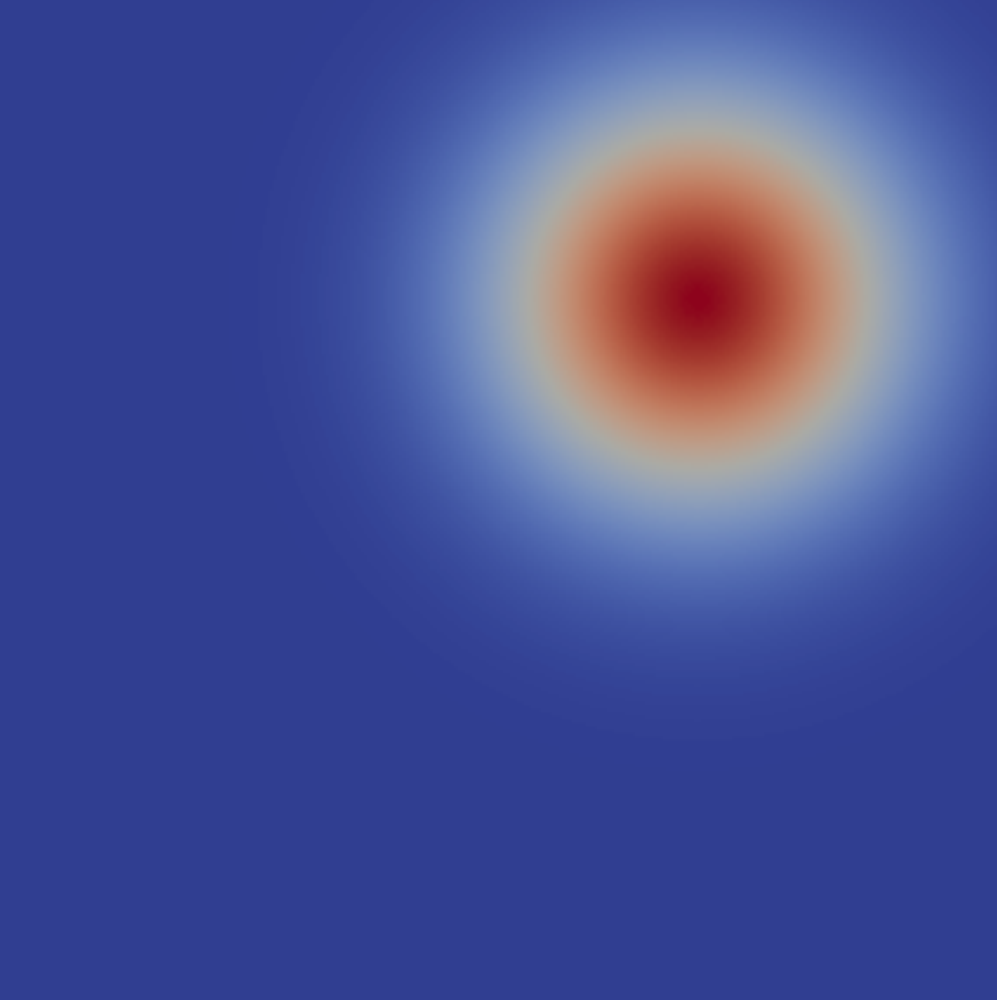}\\[0.07cm]
\includegraphics[width=.19\textwidth,angle=90,origin=c]{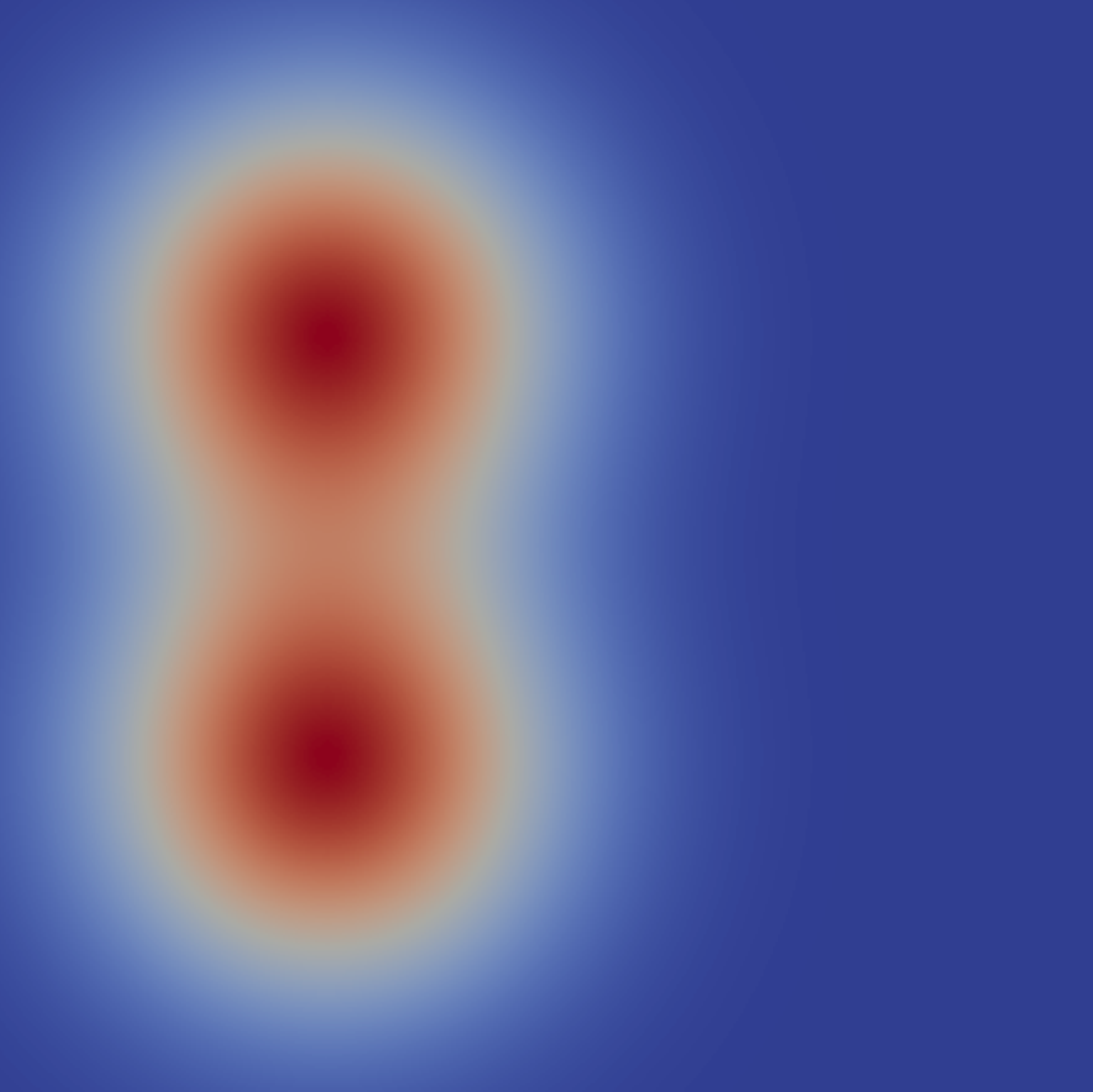}
\scalebox{-1}[1]{\includegraphics[width=.19\textwidth,angle=90,origin=c]{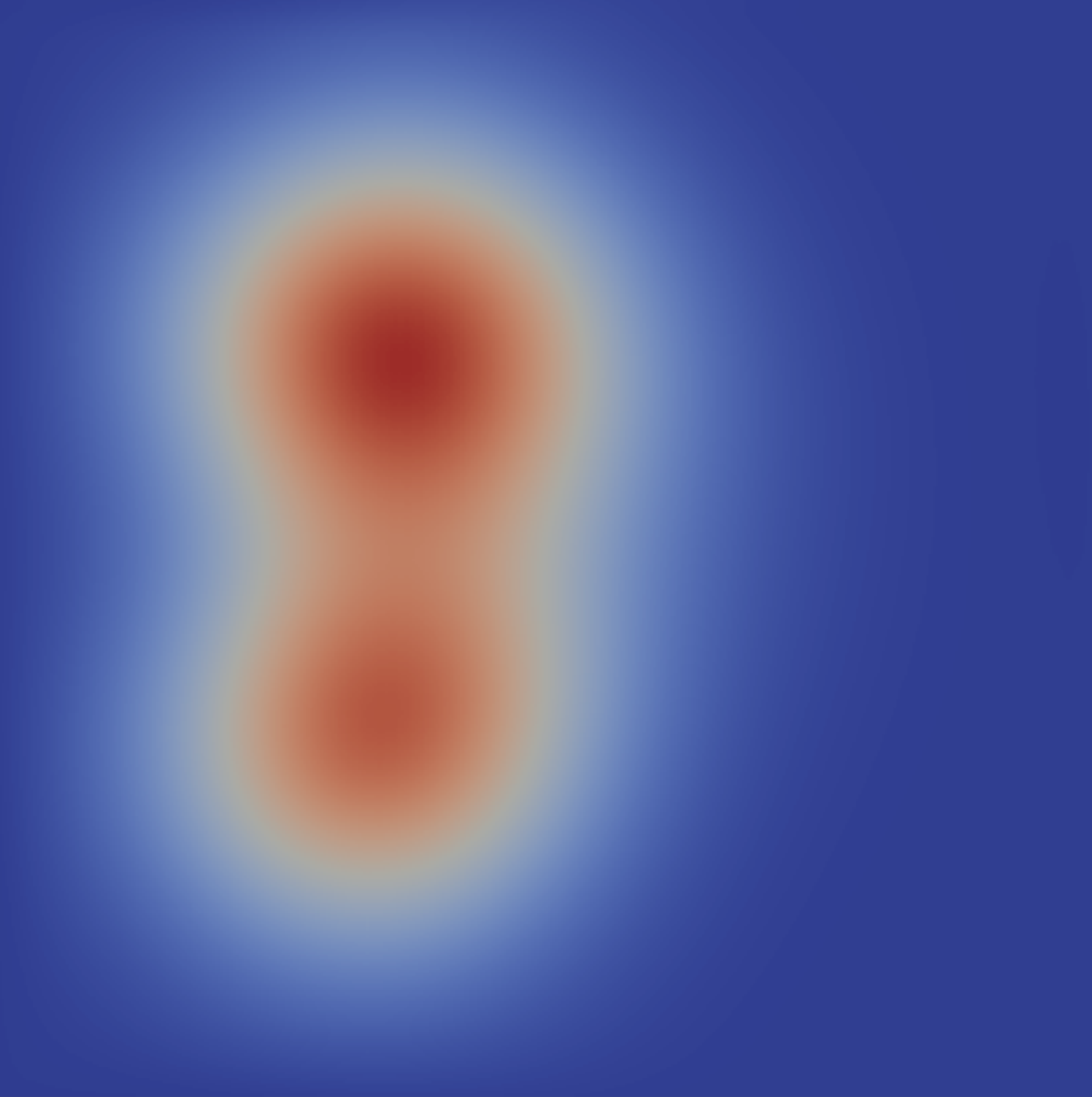}}
\scalebox{-1}[1]{\includegraphics[width=.19\textwidth,angle=90,origin=c]{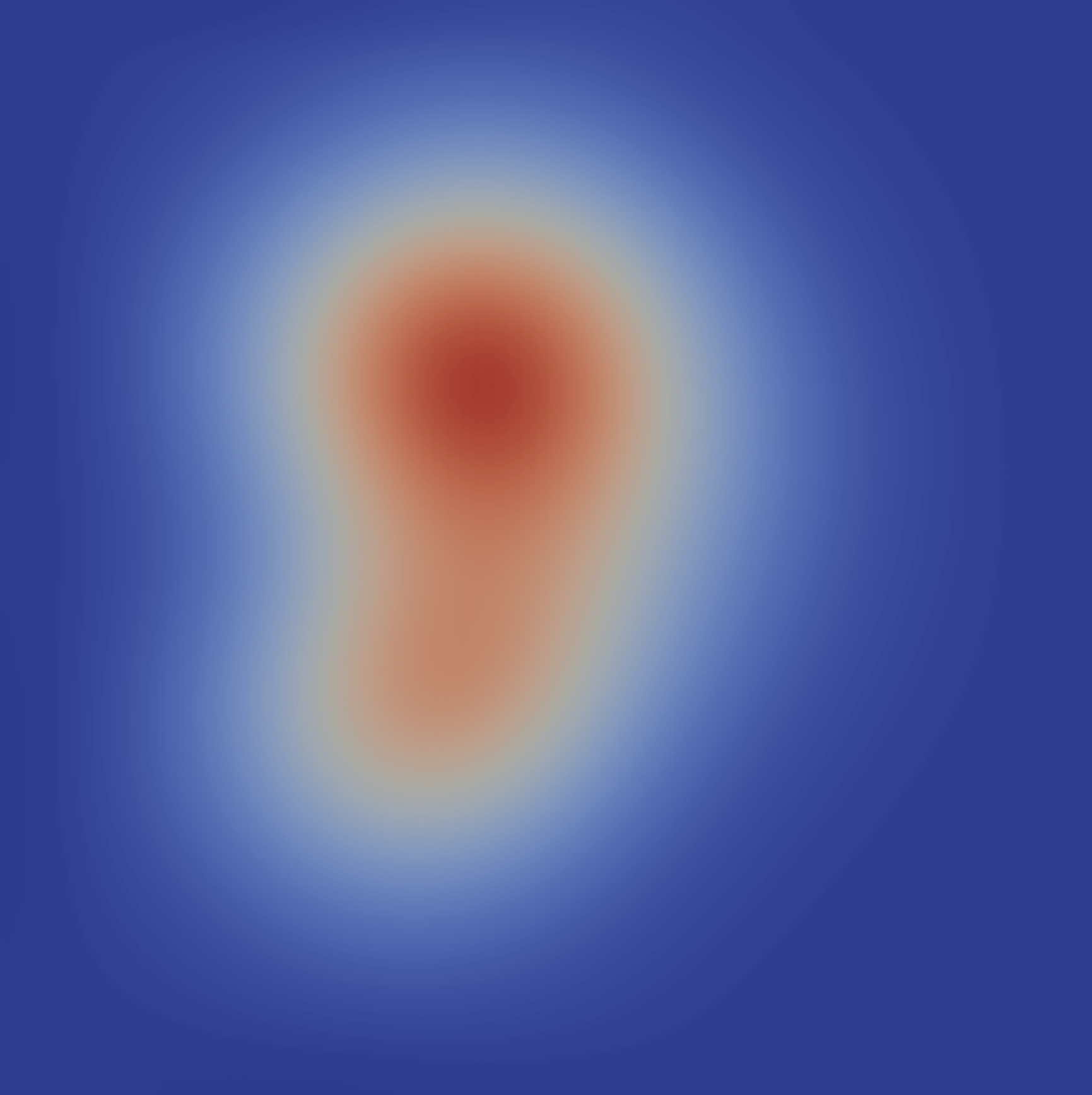}}
\scalebox{-1}[1]{\includegraphics[width=.19\textwidth,angle=90,origin=c]{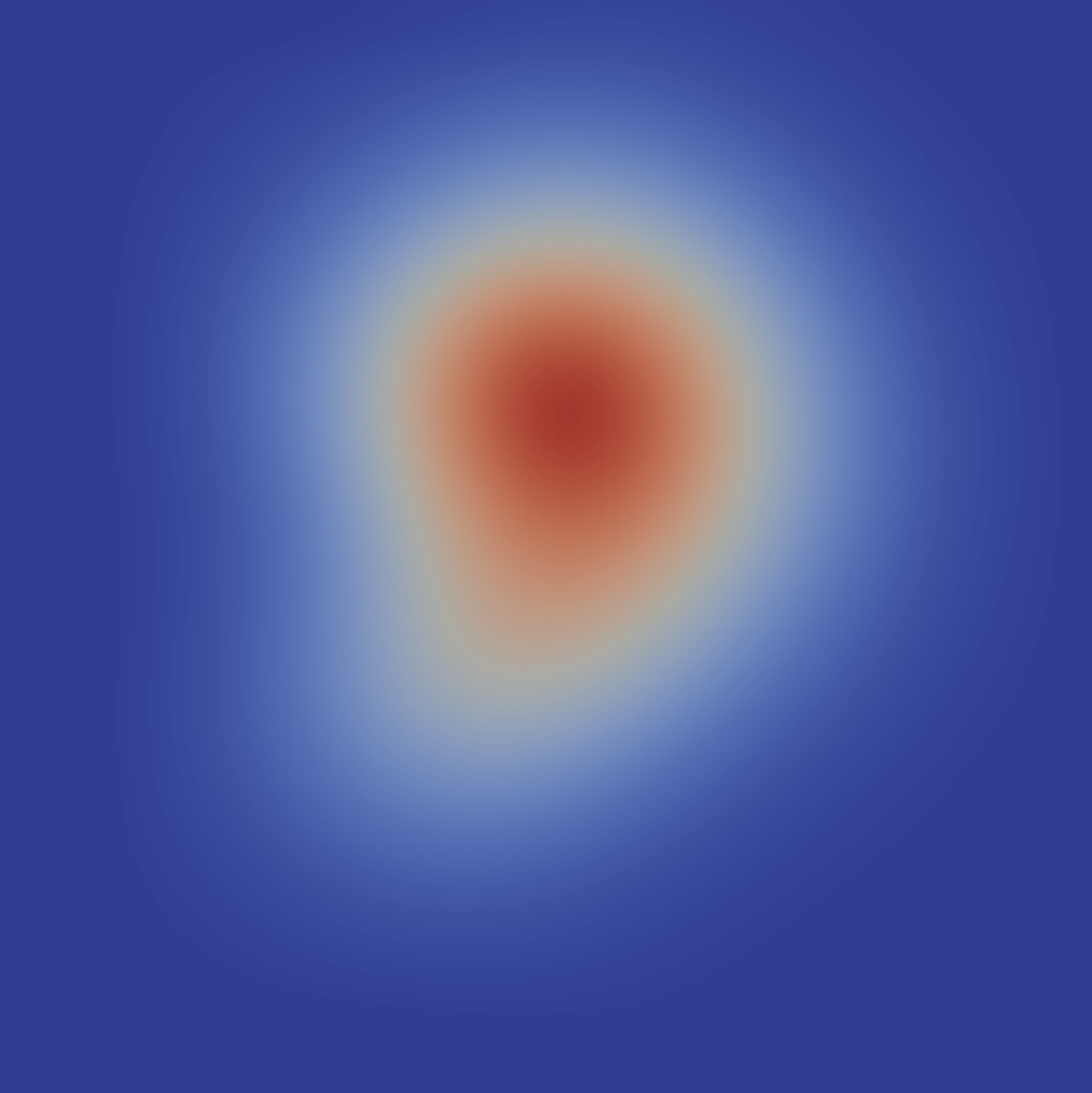}}
\includegraphics[width=.19\textwidth]{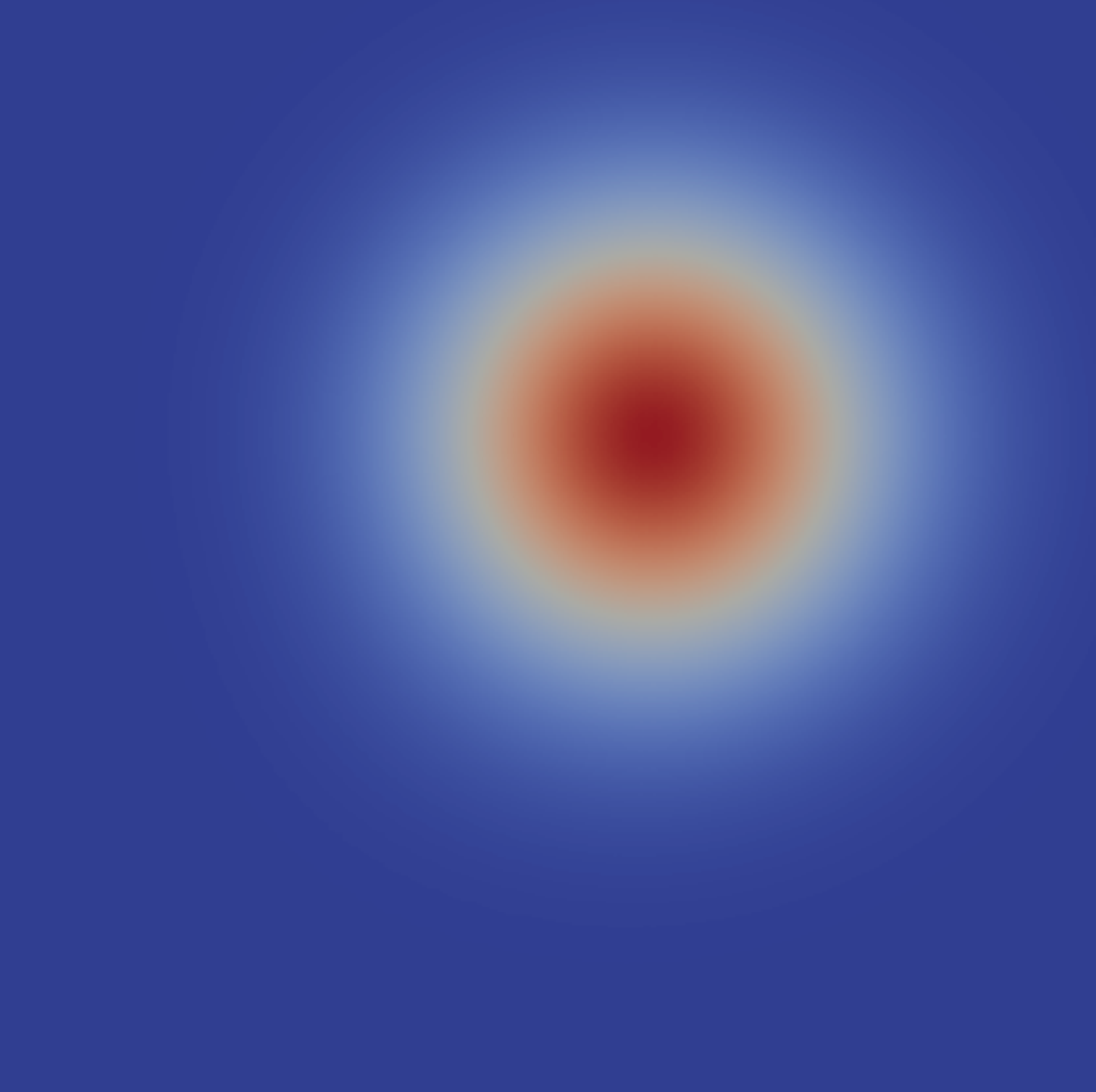}\\[0.07cm]
\includegraphics[width=.19\textwidth]{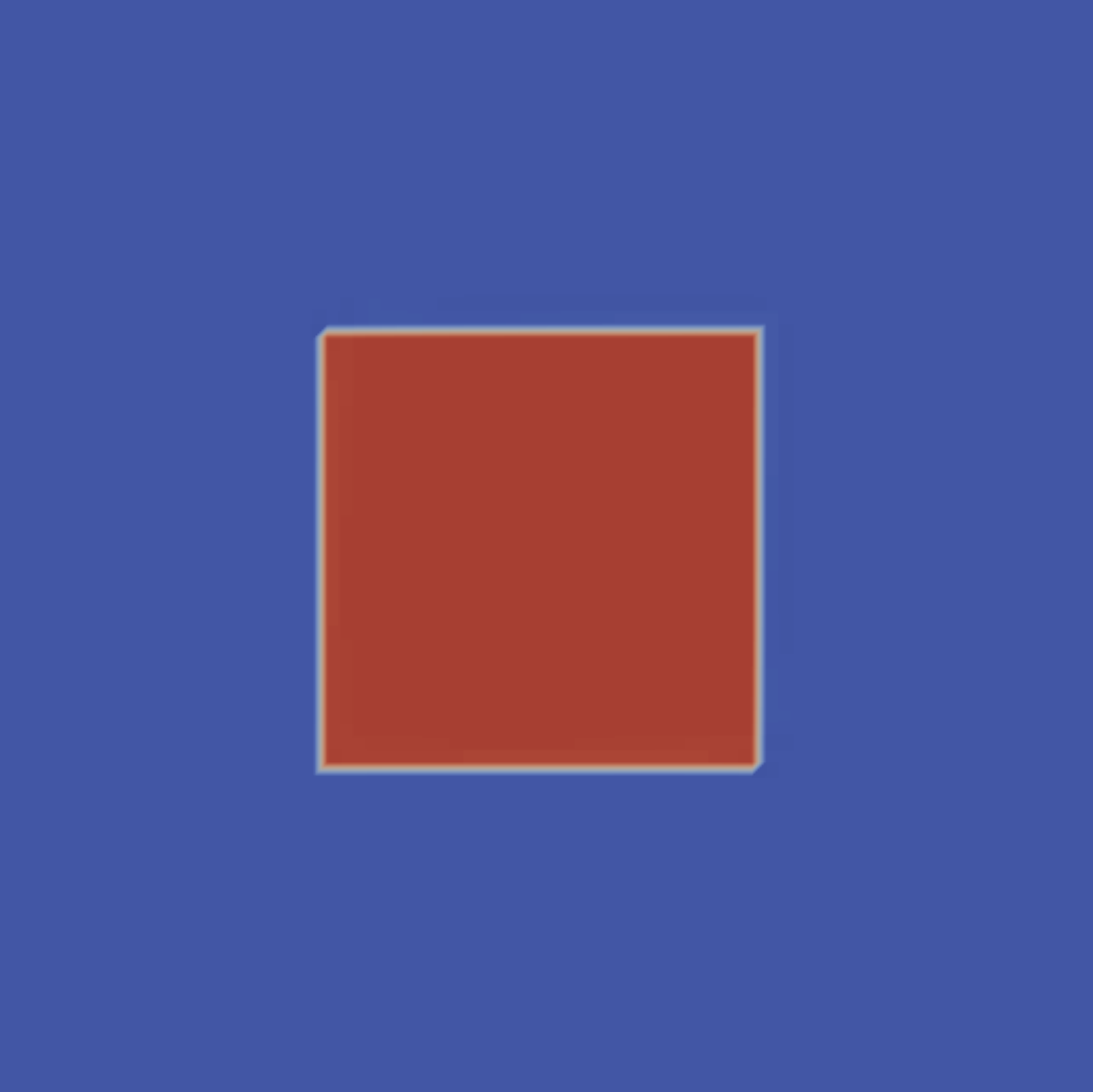}
\includegraphics[width=.19\textwidth]{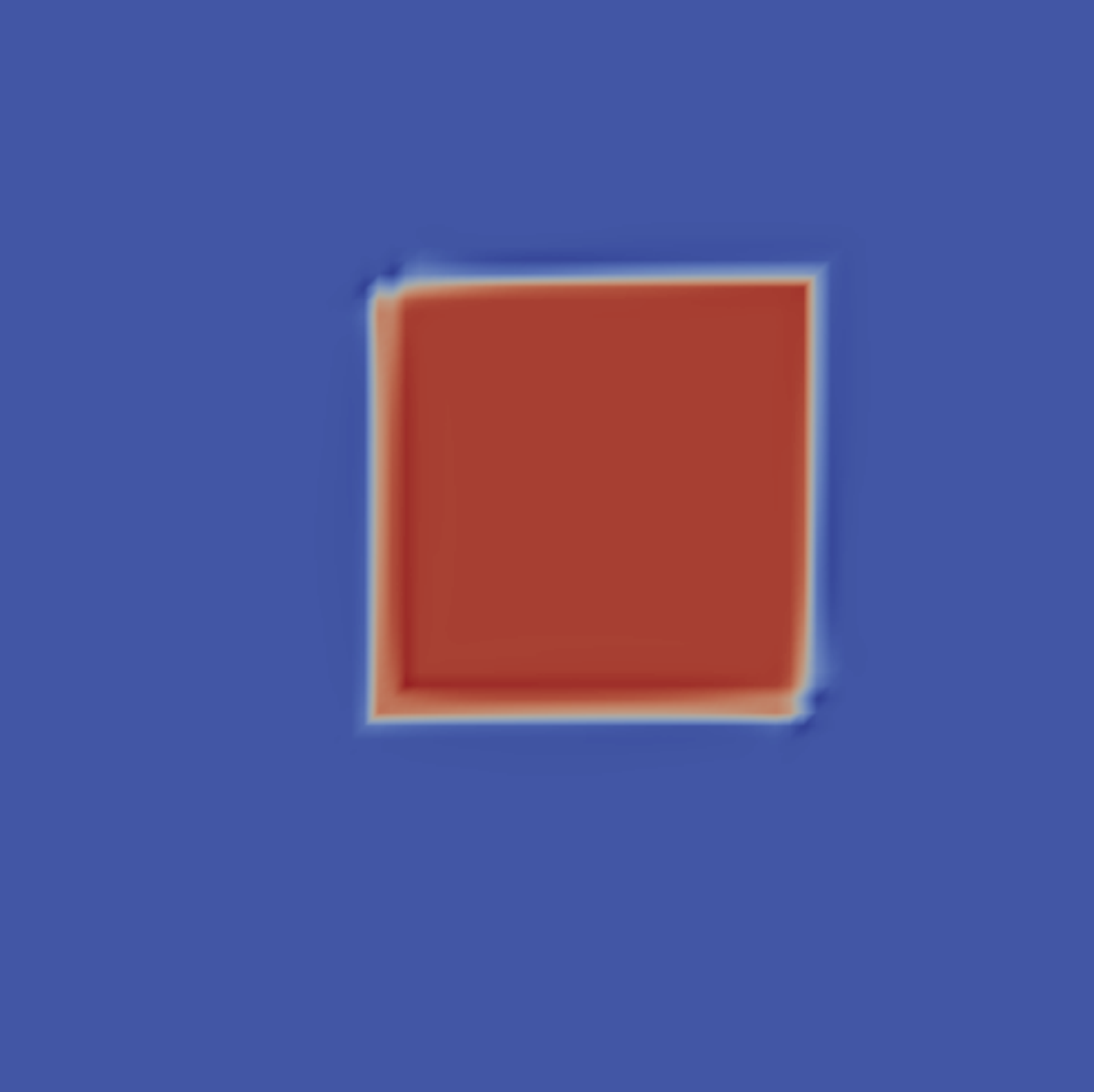}
\includegraphics[width=.19\textwidth]{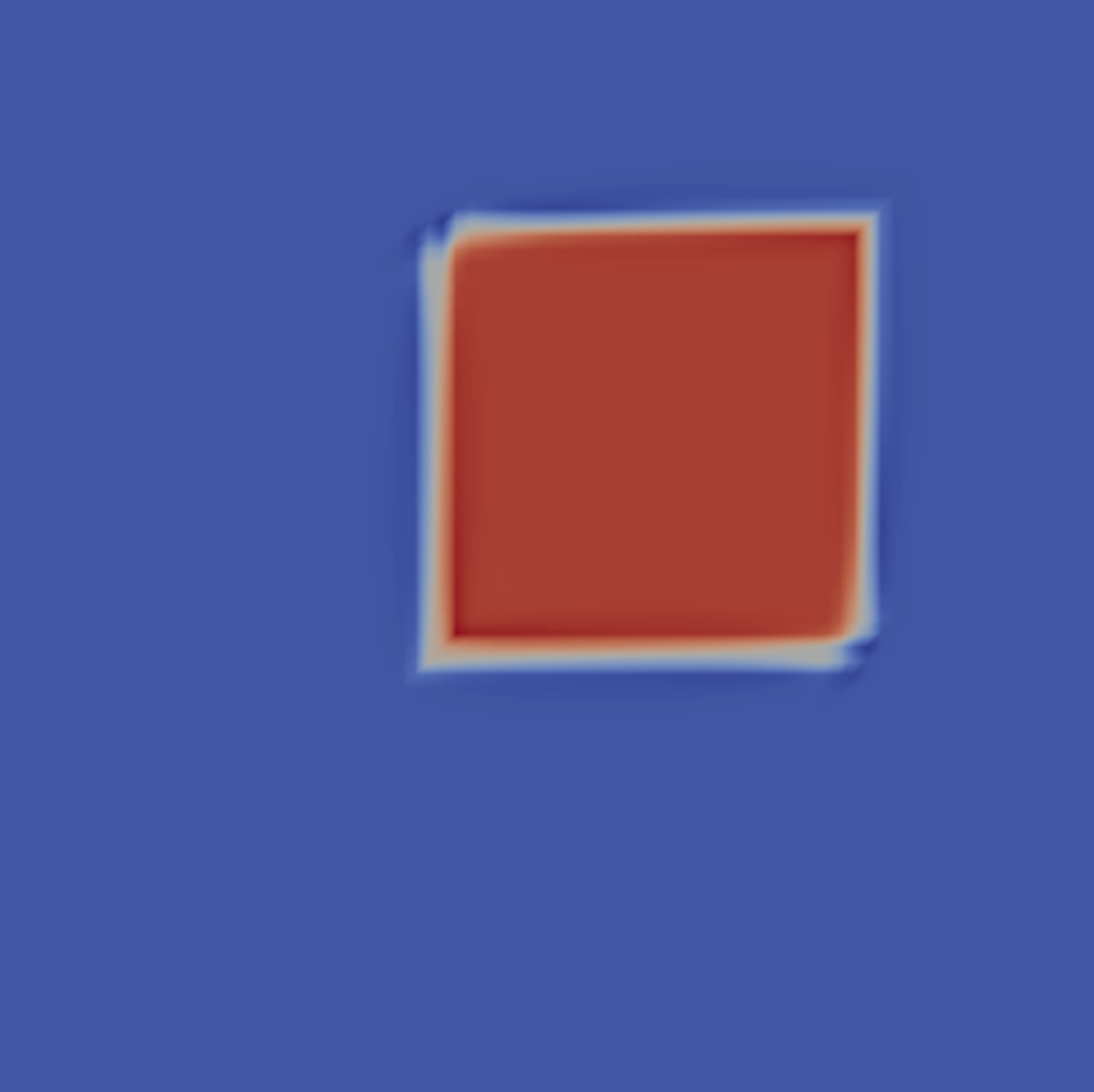}
\includegraphics[width=.19\textwidth]{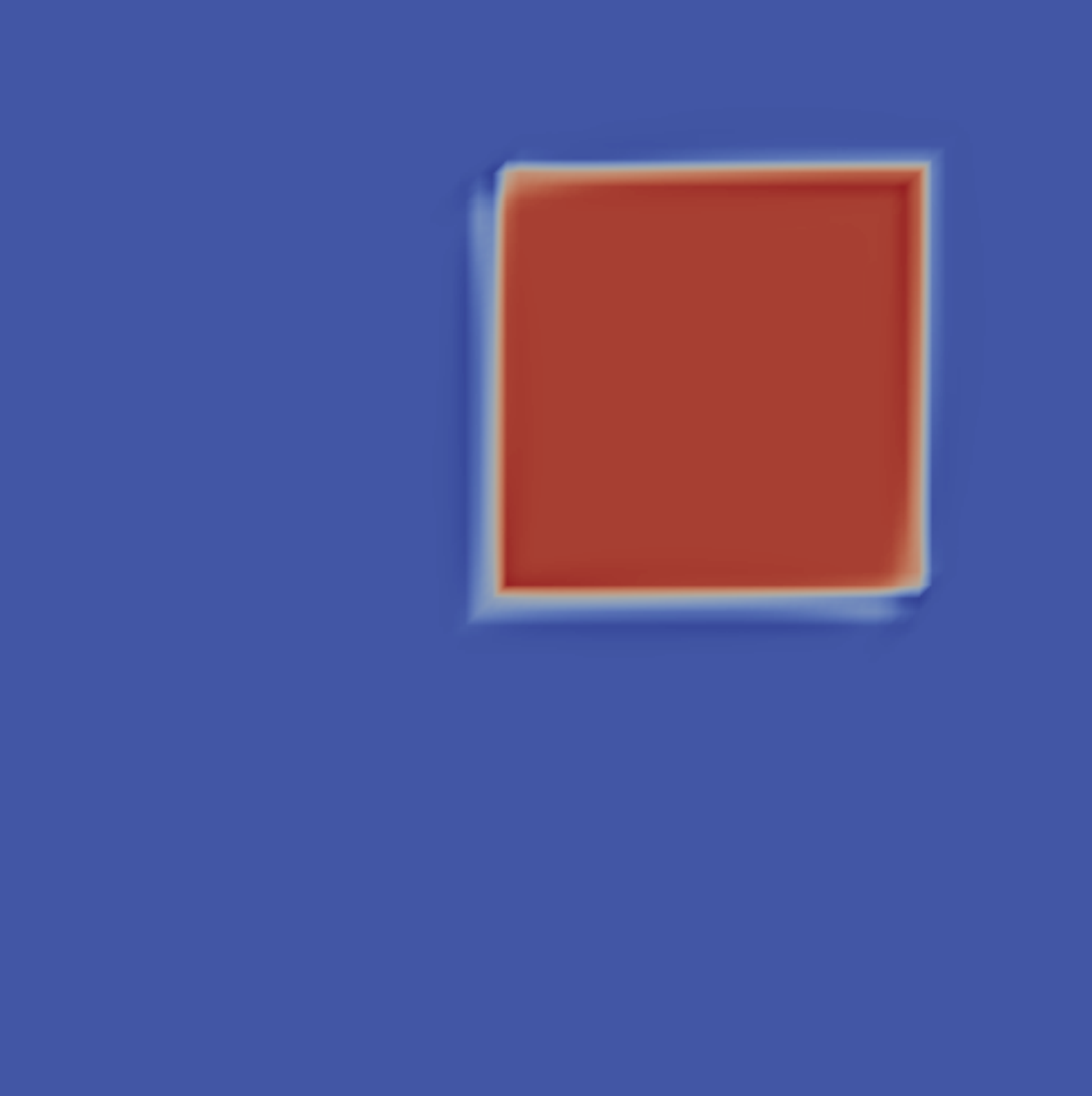}
\includegraphics[width=.19\textwidth]{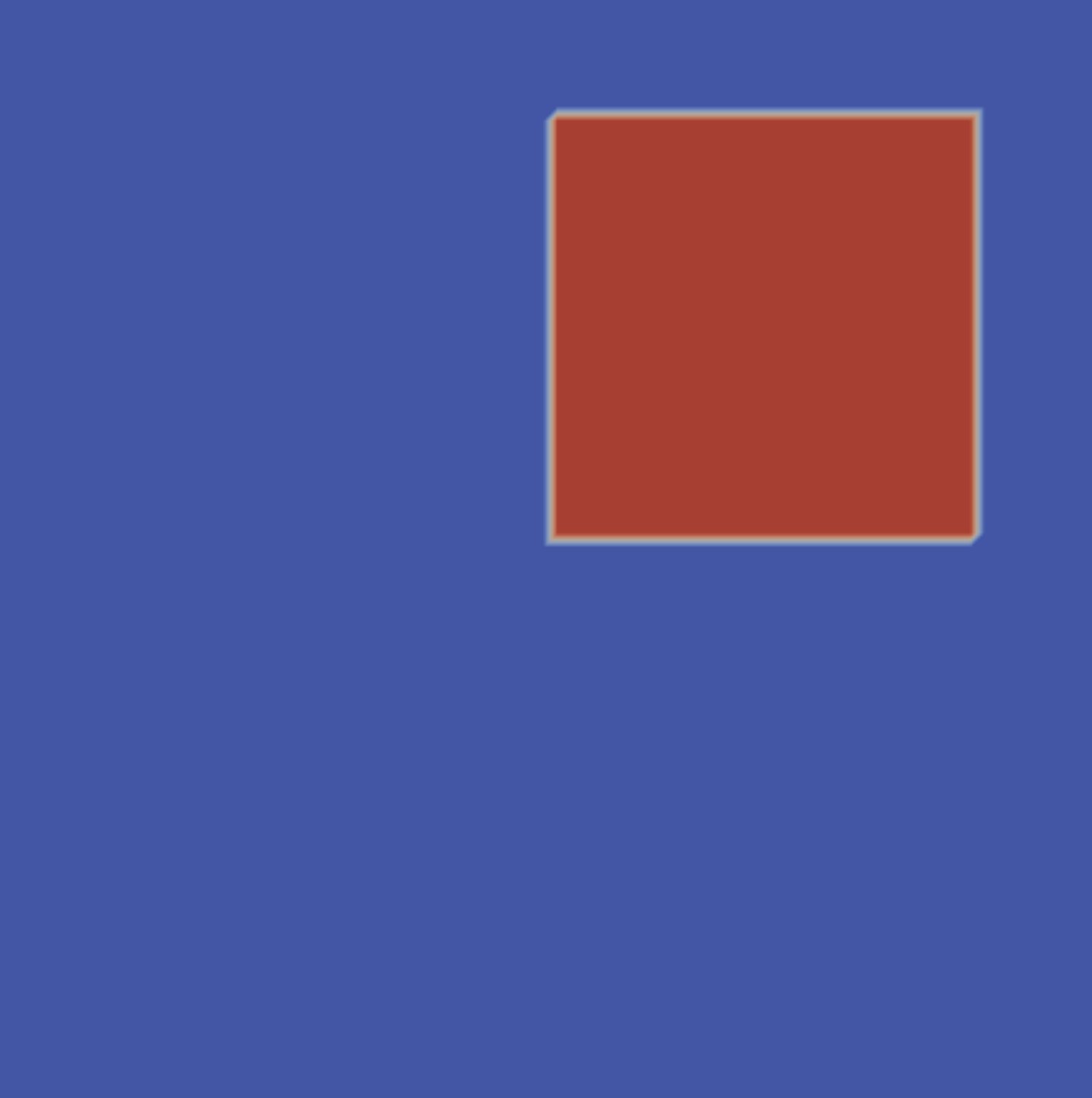}\\[0.07cm]
\includegraphics[width=.19\textwidth]{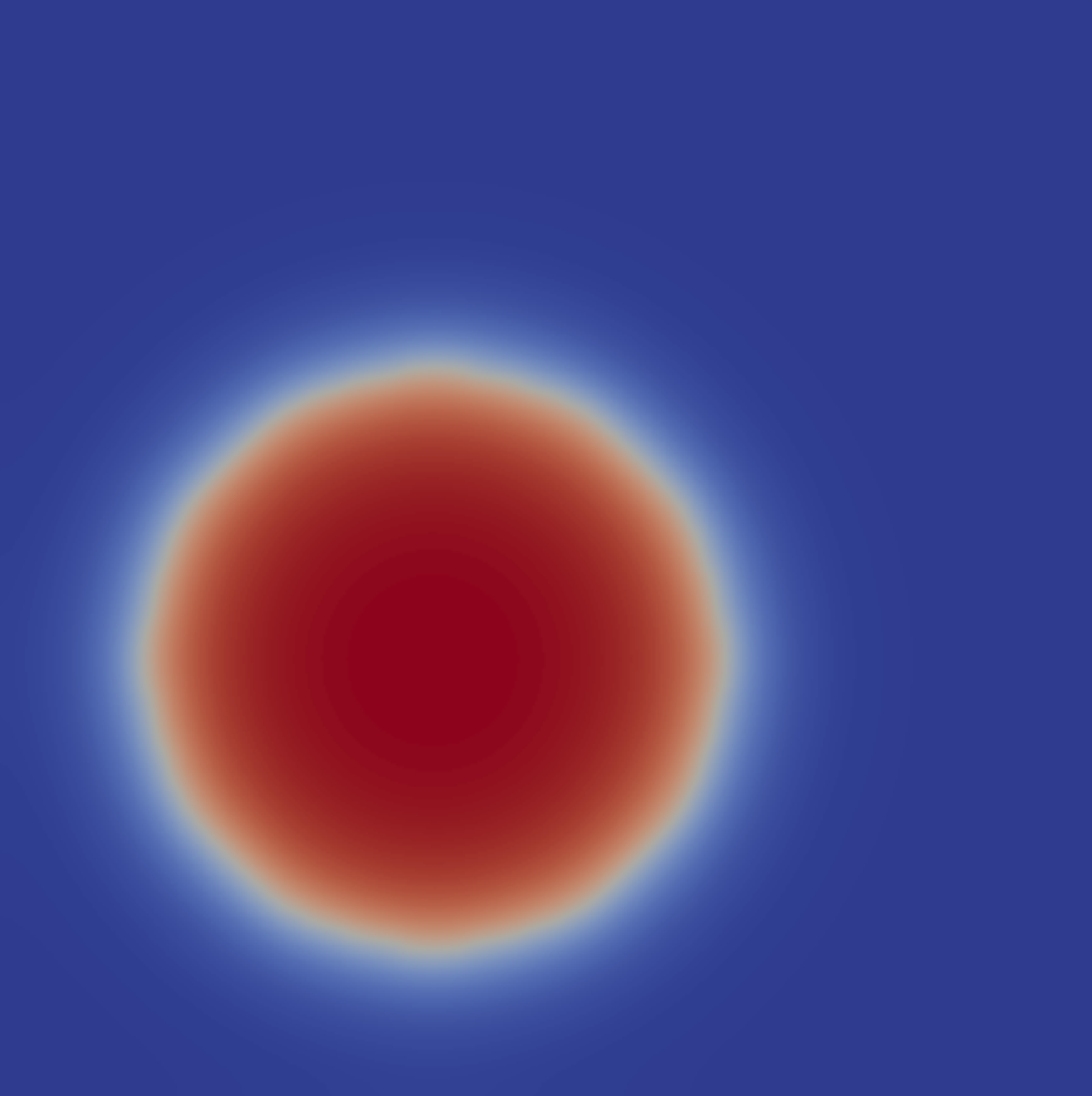}
\includegraphics[width=.19\textwidth]{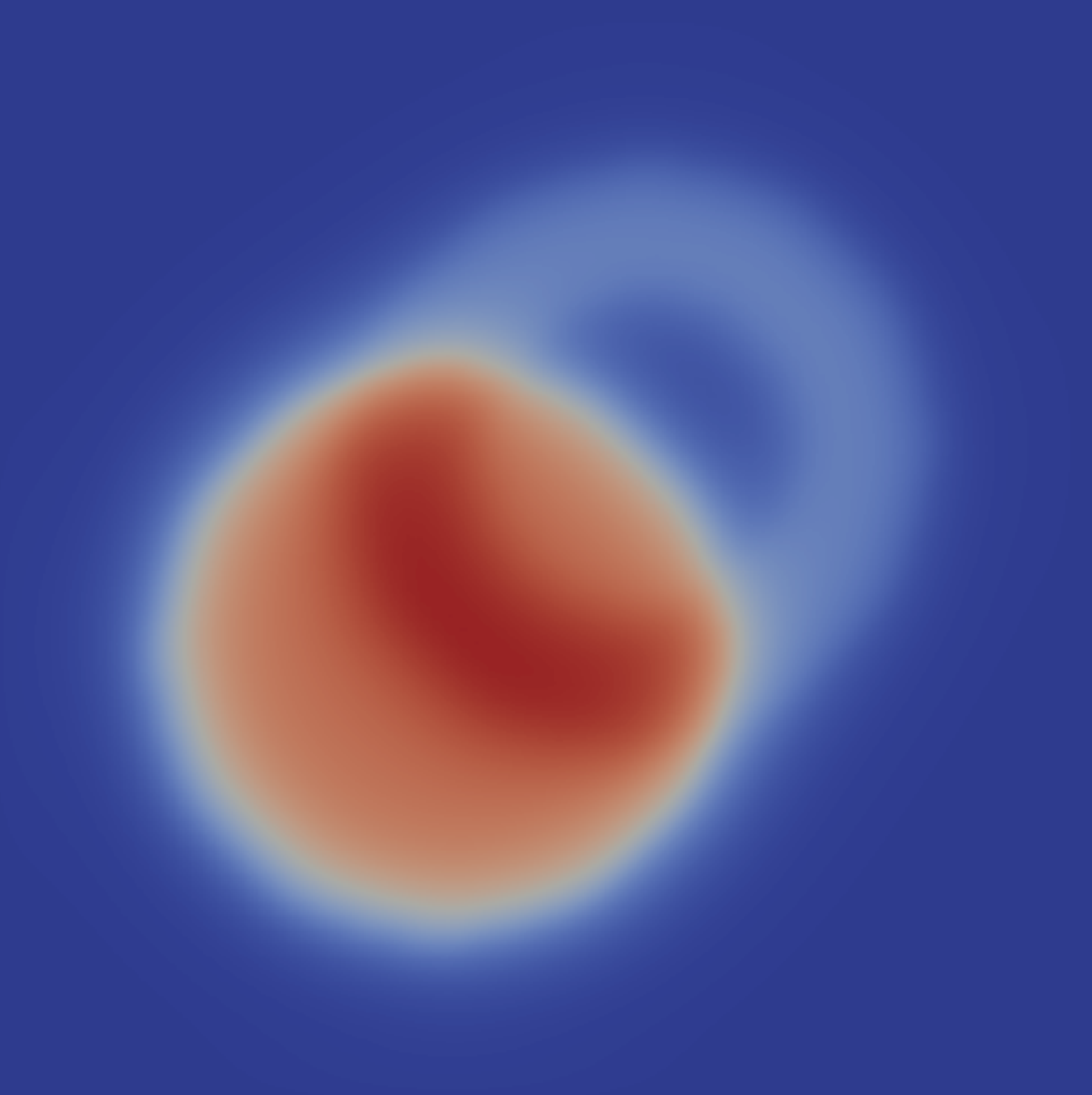}
\includegraphics[width=.19\textwidth]{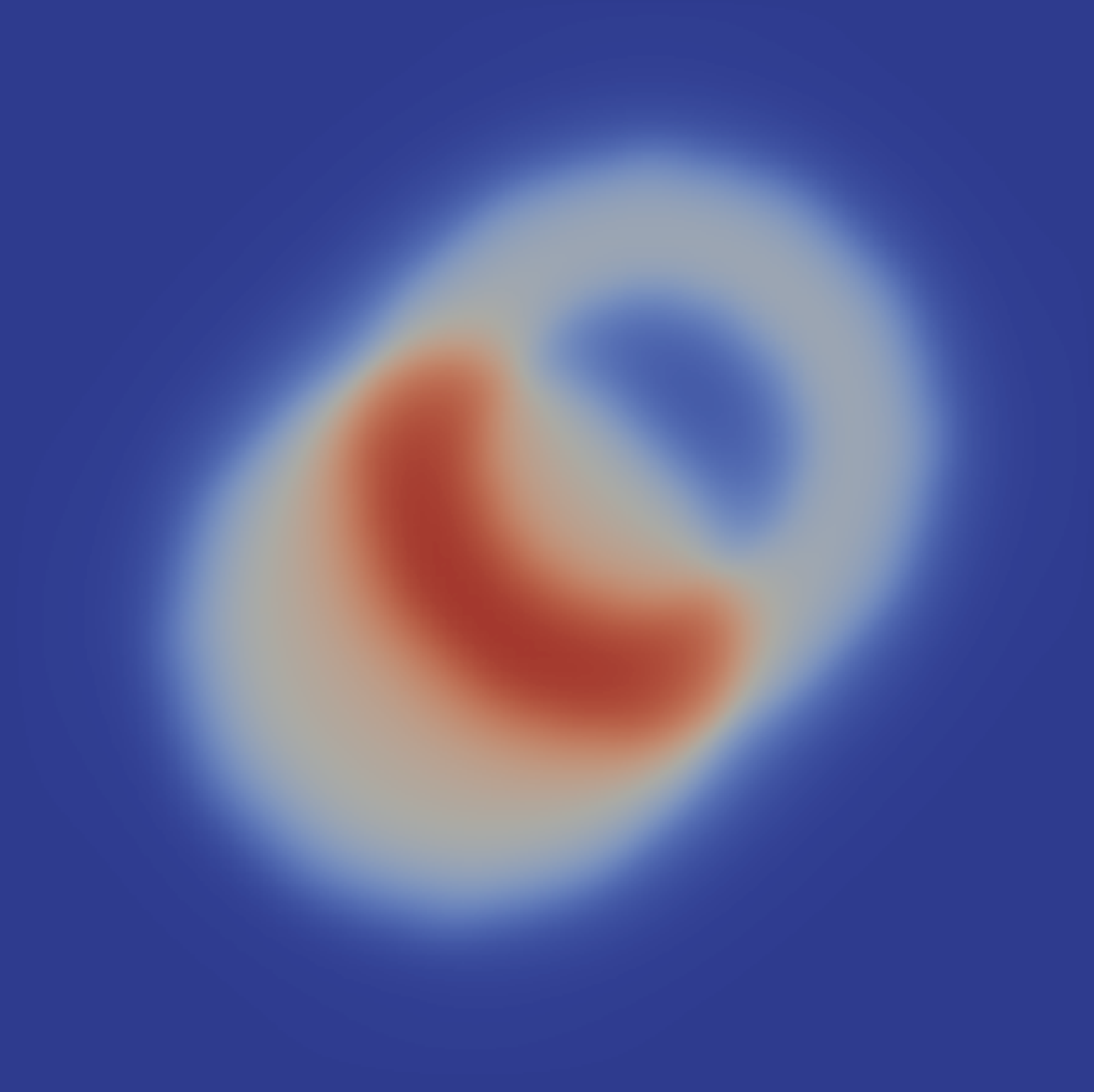}
\includegraphics[width=.19\textwidth]{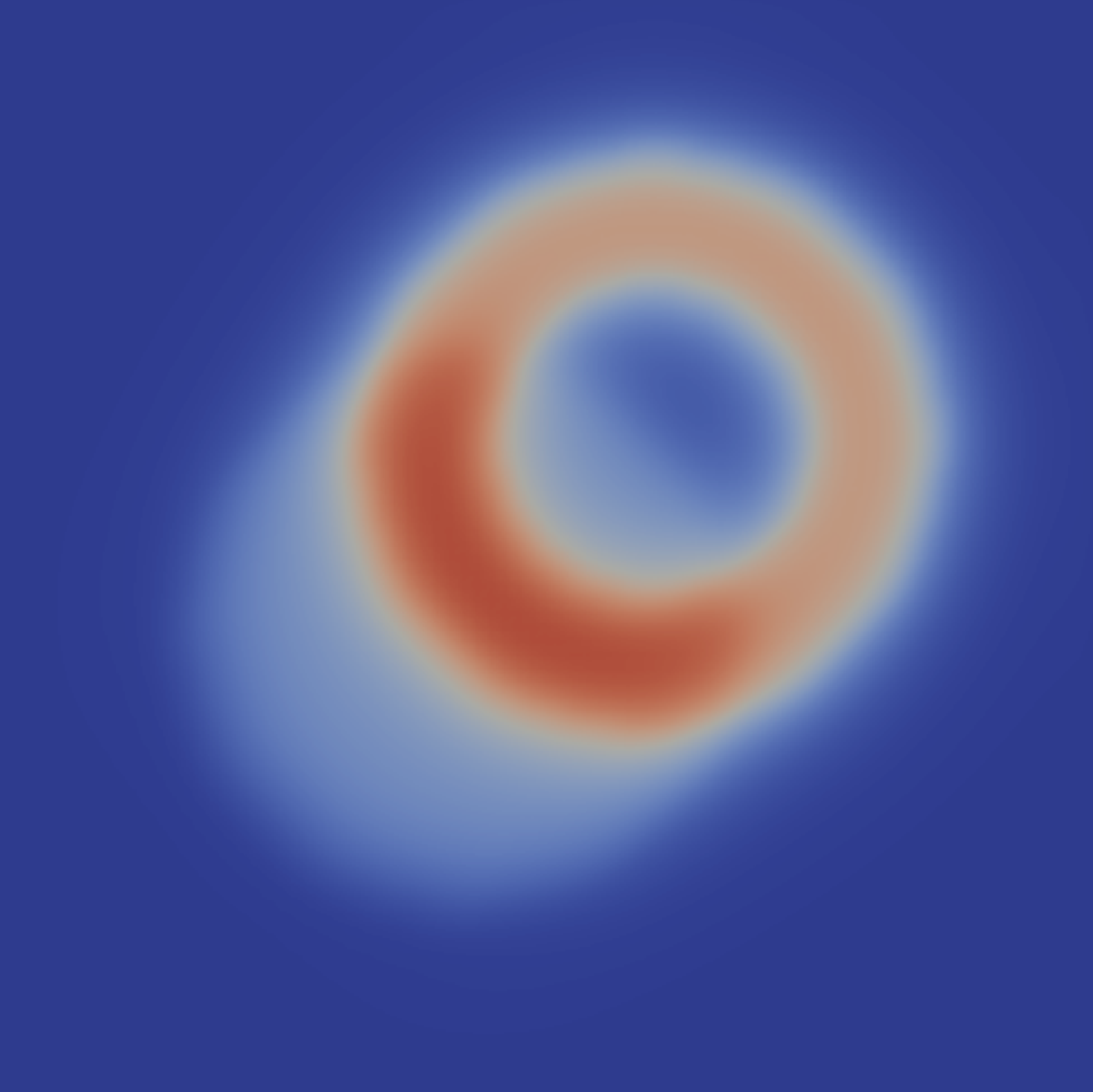}
\includegraphics[width=.19\textwidth]{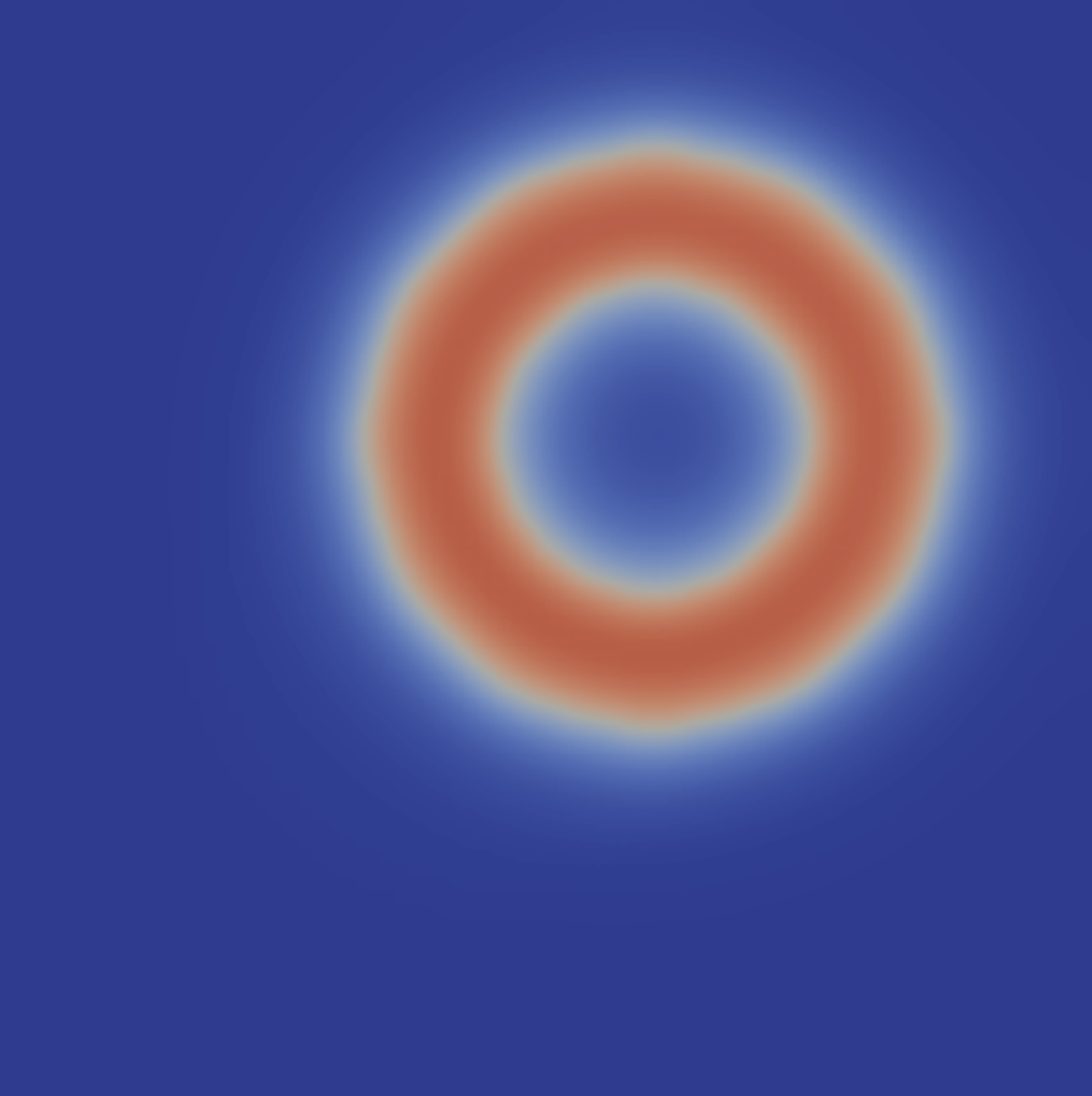}\\[0.07cm]
\includegraphics[width=.19\textwidth]{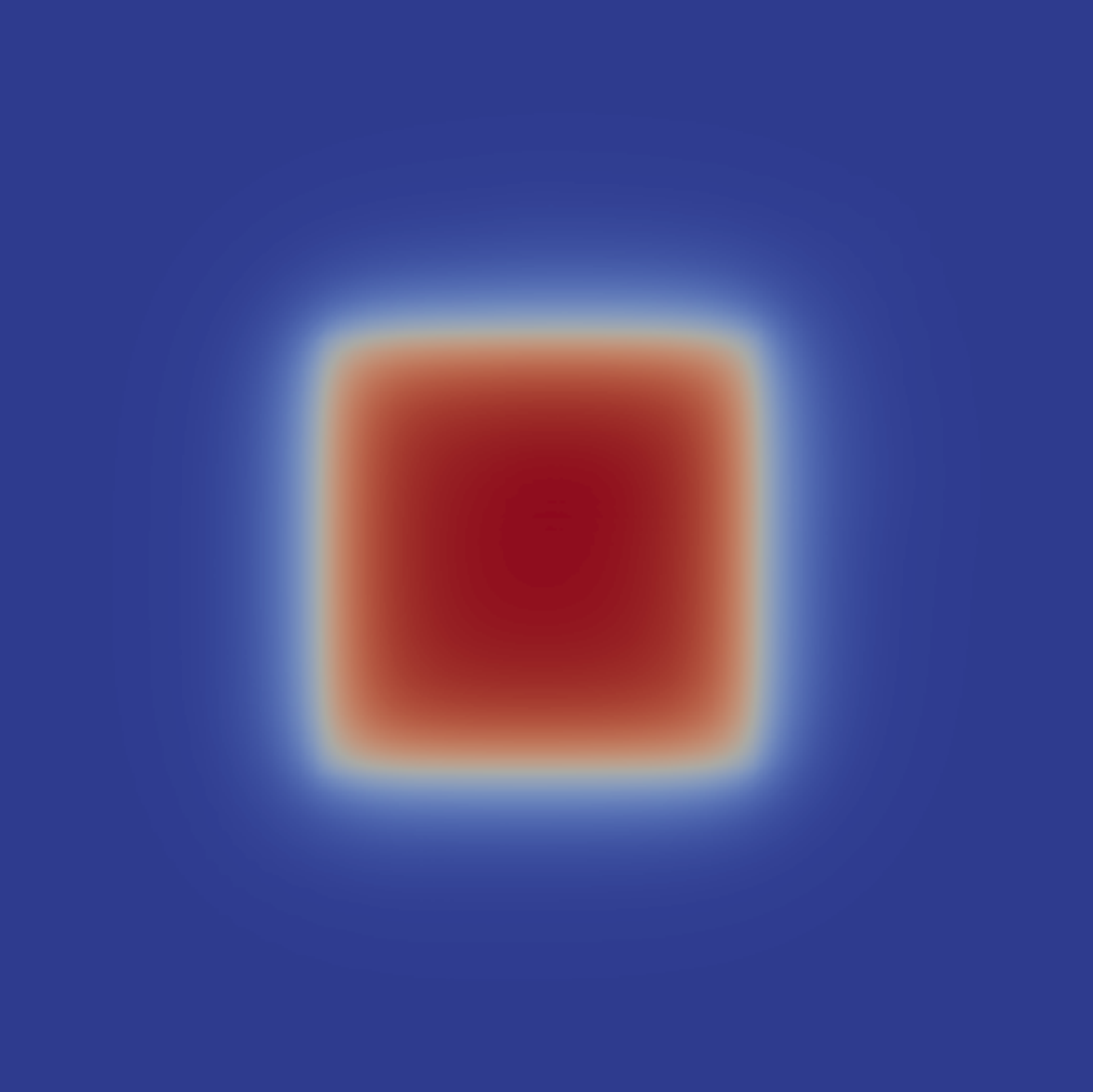}
\includegraphics[width=.19\textwidth]{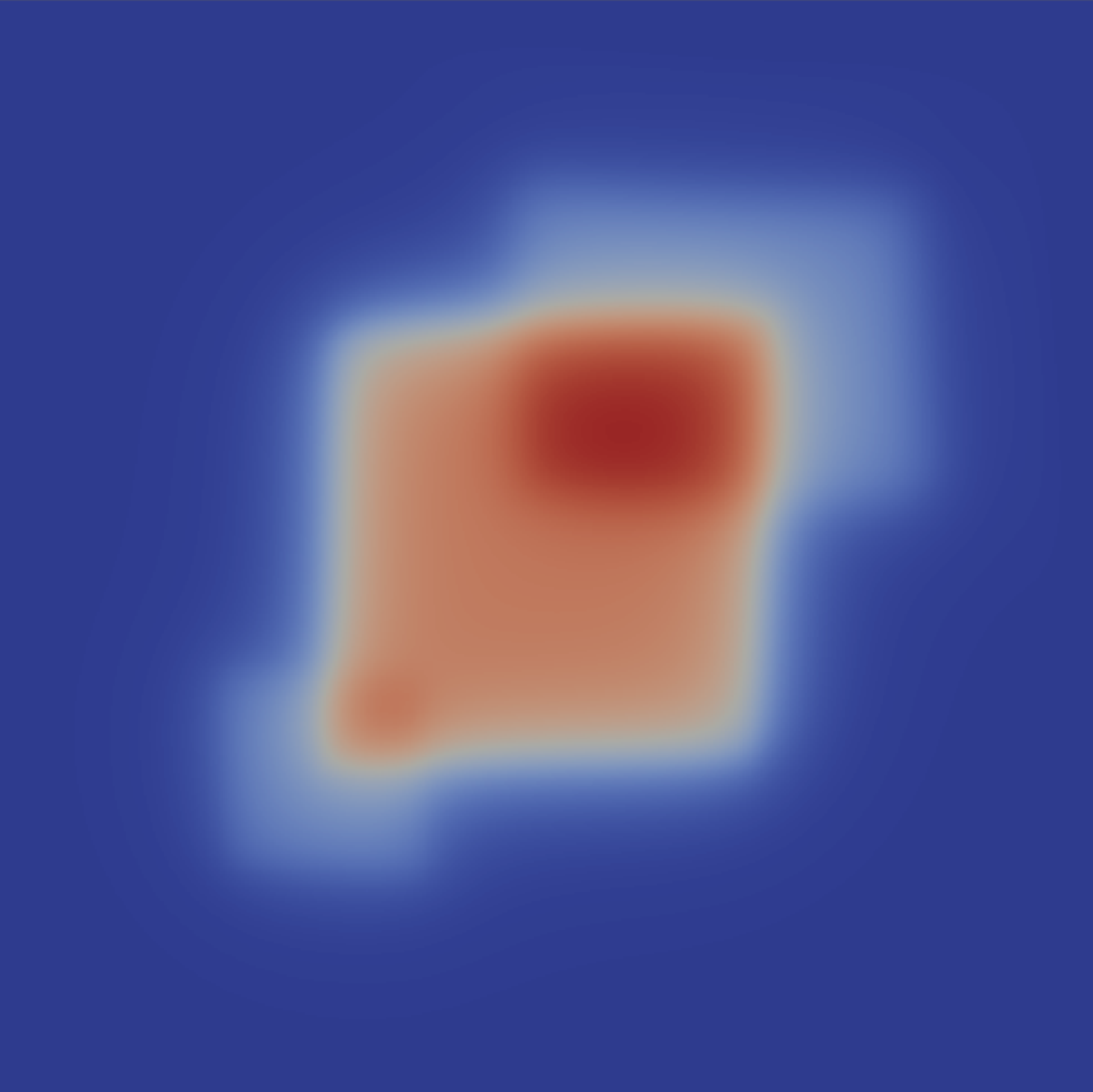}
\includegraphics[width=.19\textwidth]{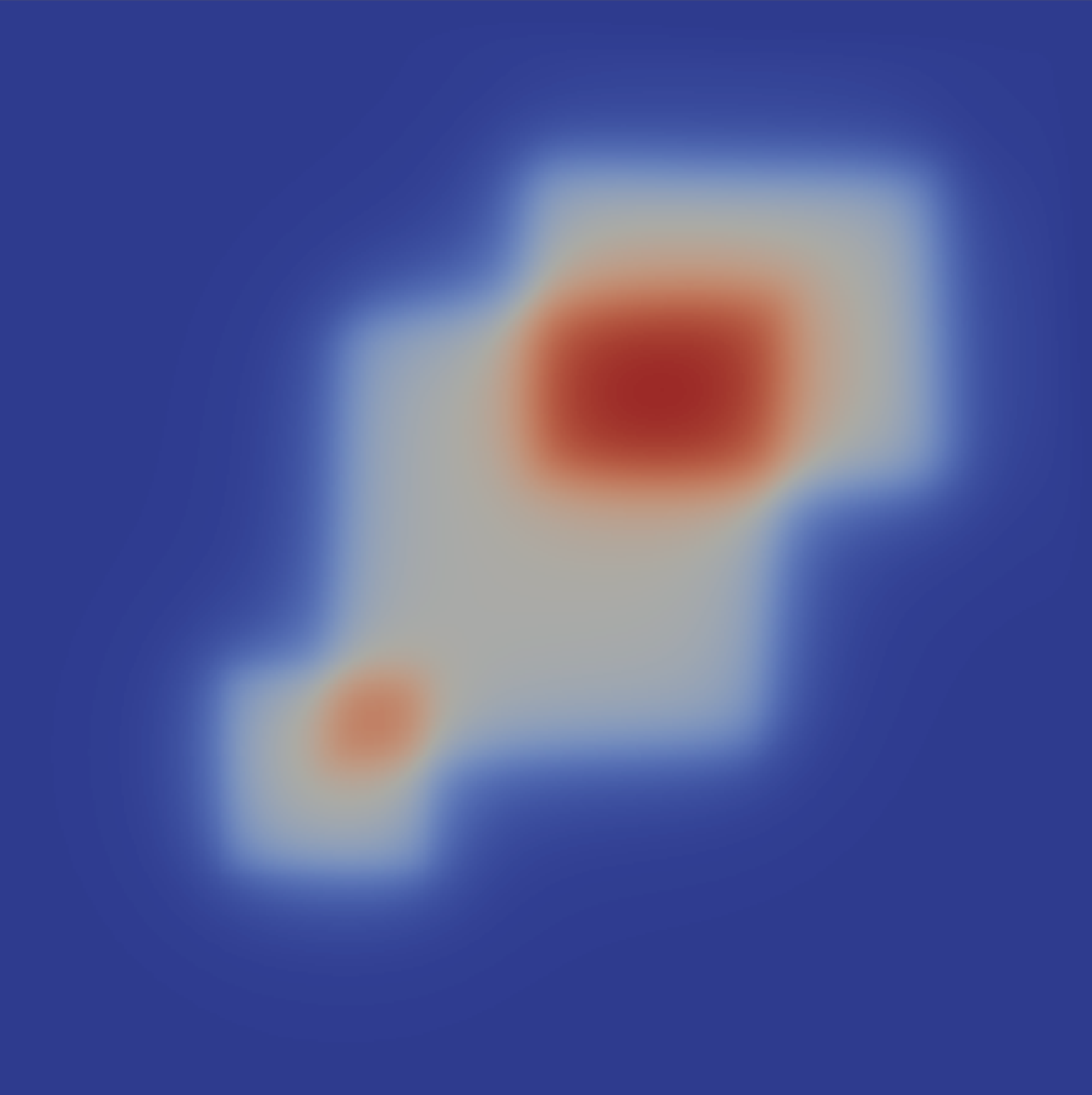}
\includegraphics[width=.19\textwidth]{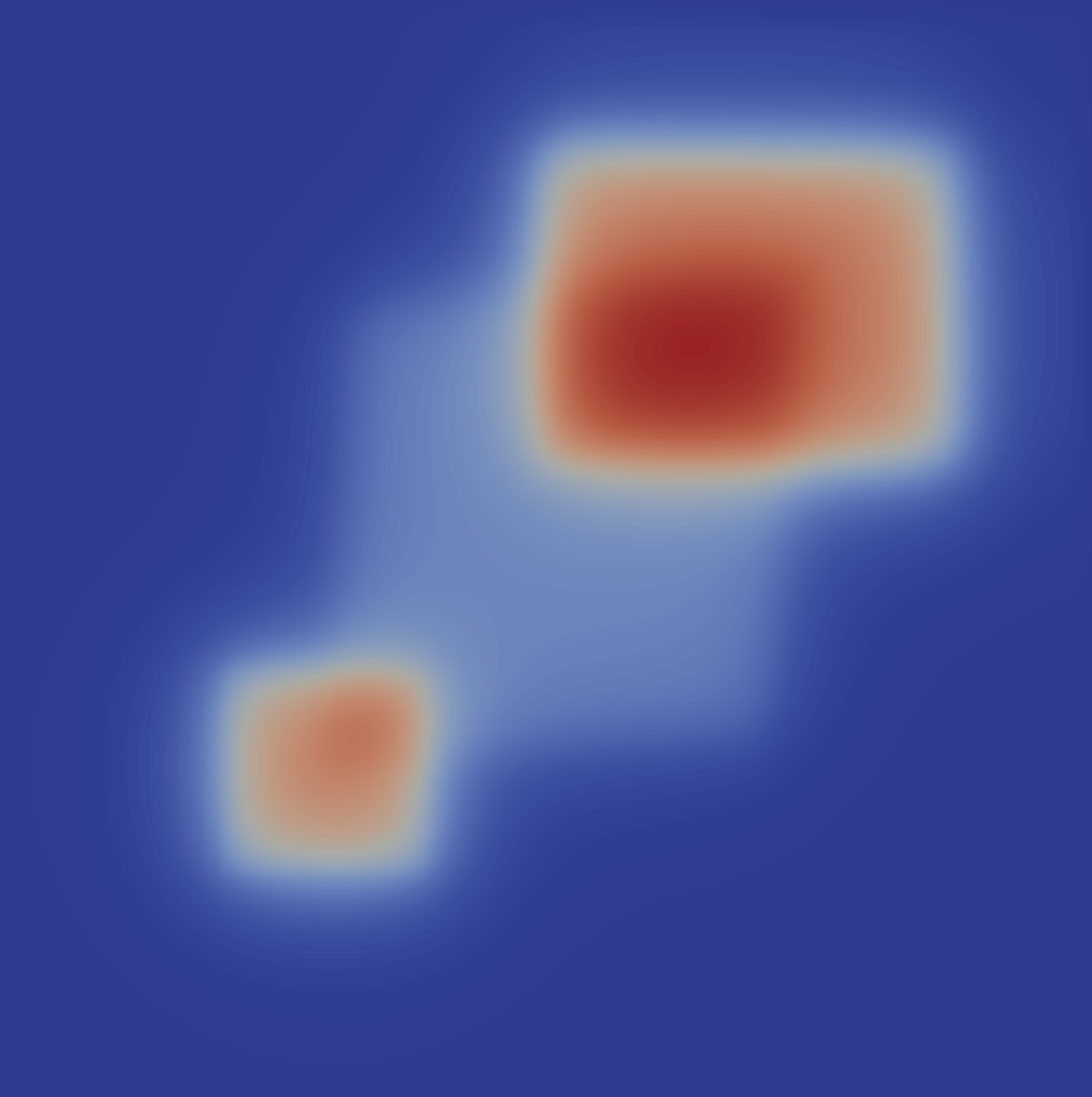}
\includegraphics[width=.19\textwidth]{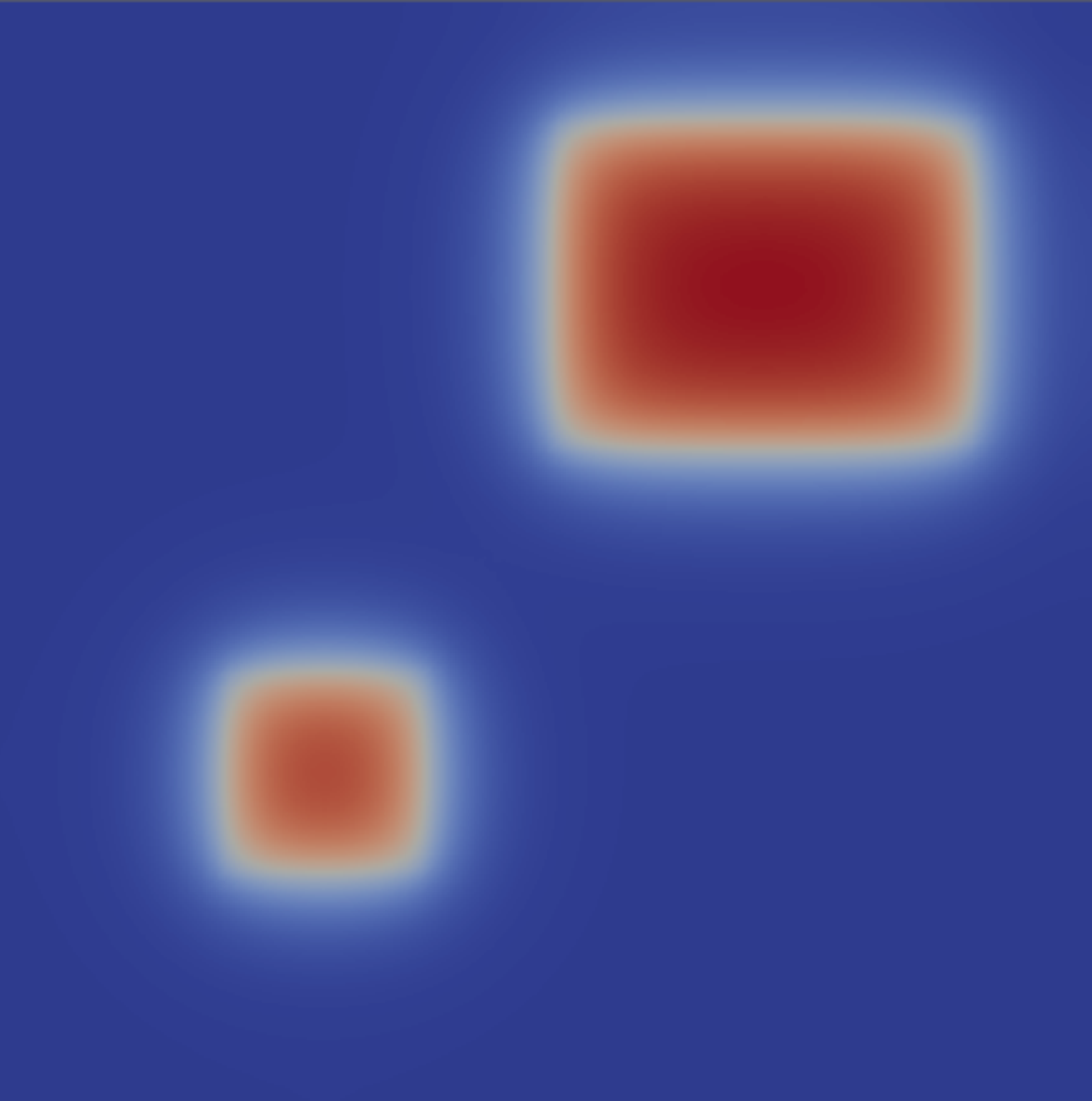}
\end{minipage}
\caption[Results for the outer problem for various template and target pairs in
two dimensions using $\sigma^{-2}=1$.]{Results for the outer problem for
various template and target pairs in two dimensions for $\sigma^{-2}=1$. Each
row shows the evolution at time $t=0,0.25,0.5,0.75,1$.}
\label{fig:adjoint2d_c=1}
\end{figure}

\begin{figure}[h!]
\begin{minipage}{\textwidth}
\centering
\includegraphics[width=.19\textwidth]{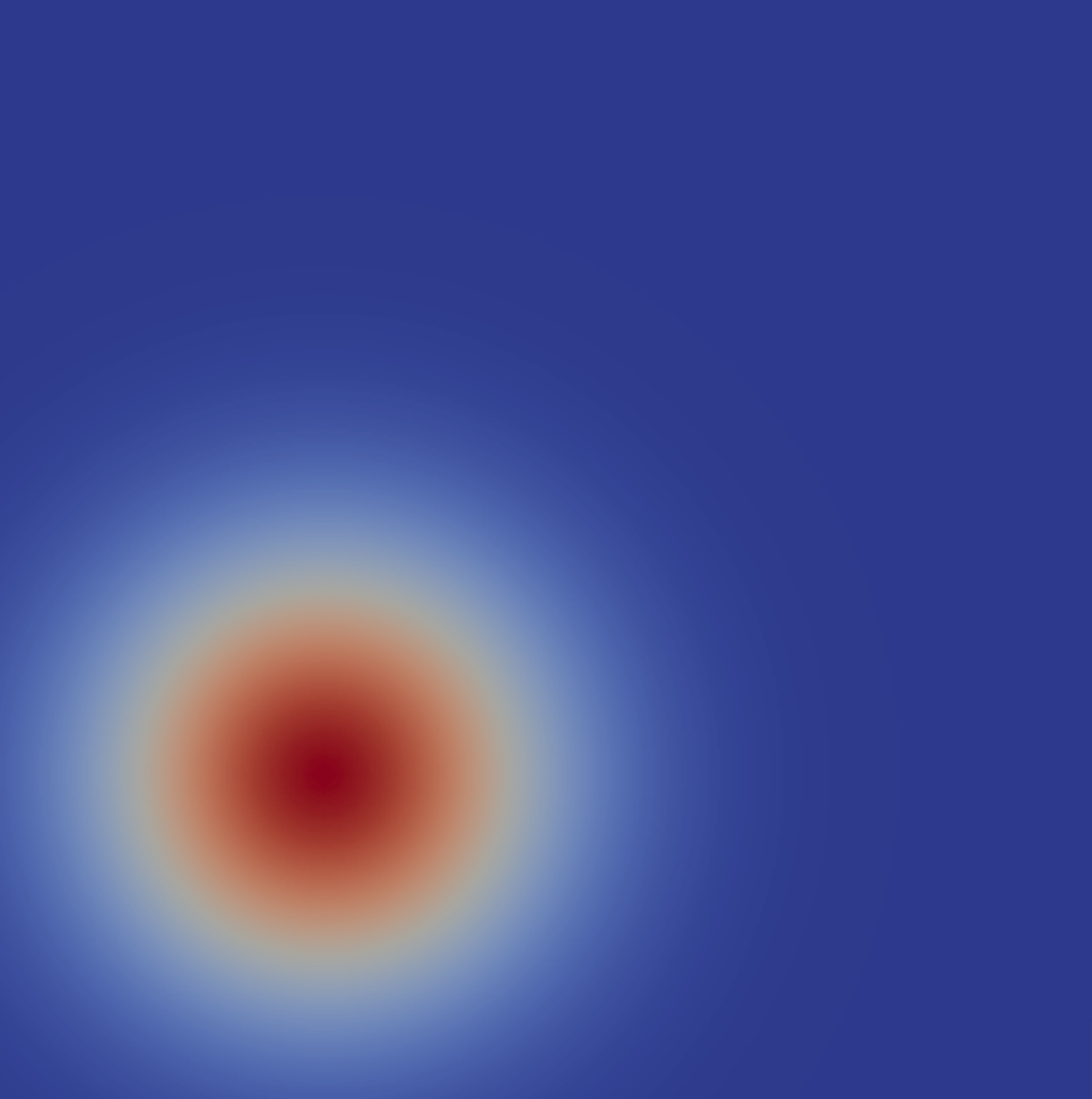}
\includegraphics[width=.19\textwidth]{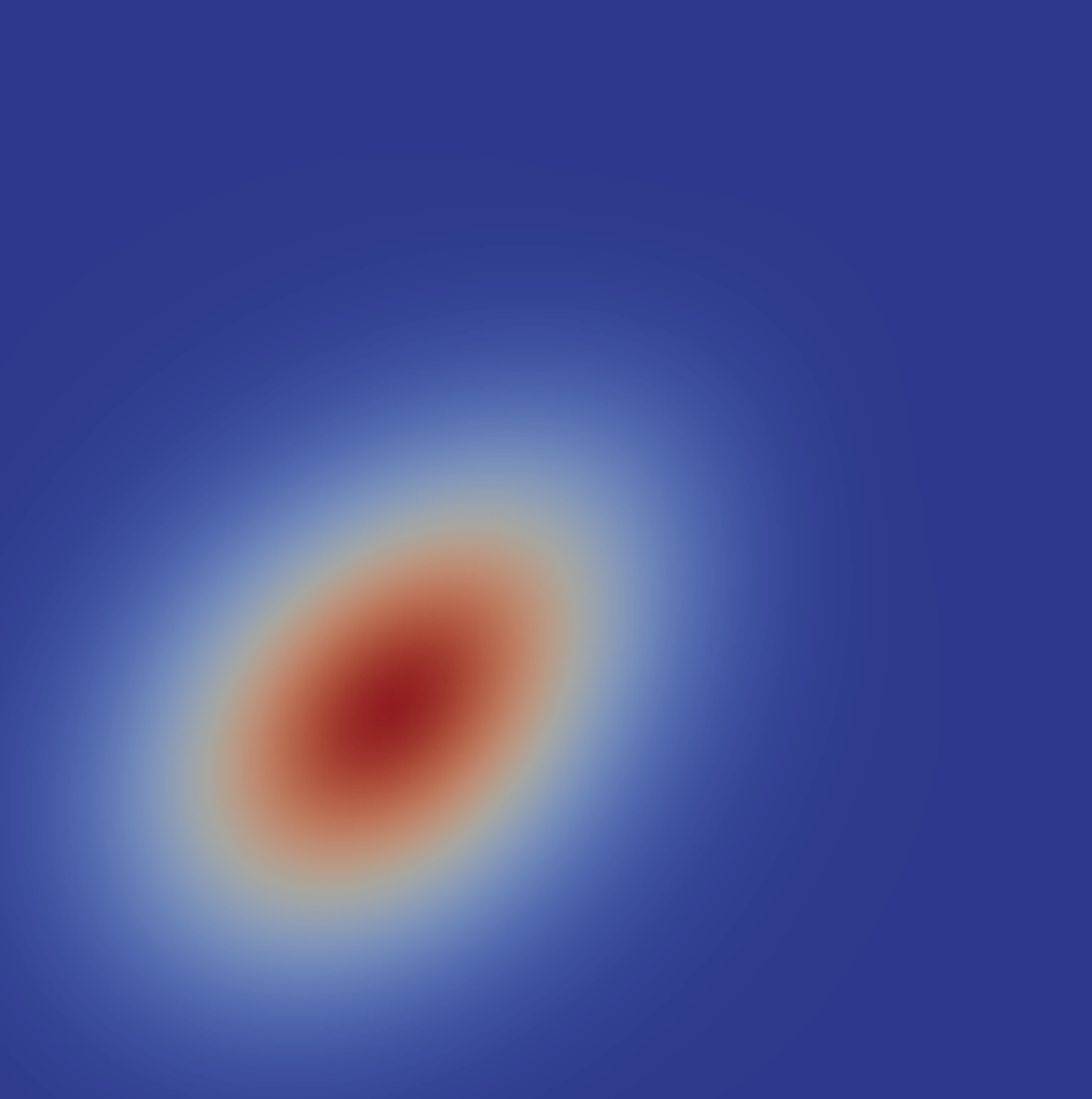}
\includegraphics[width=.19\textwidth]{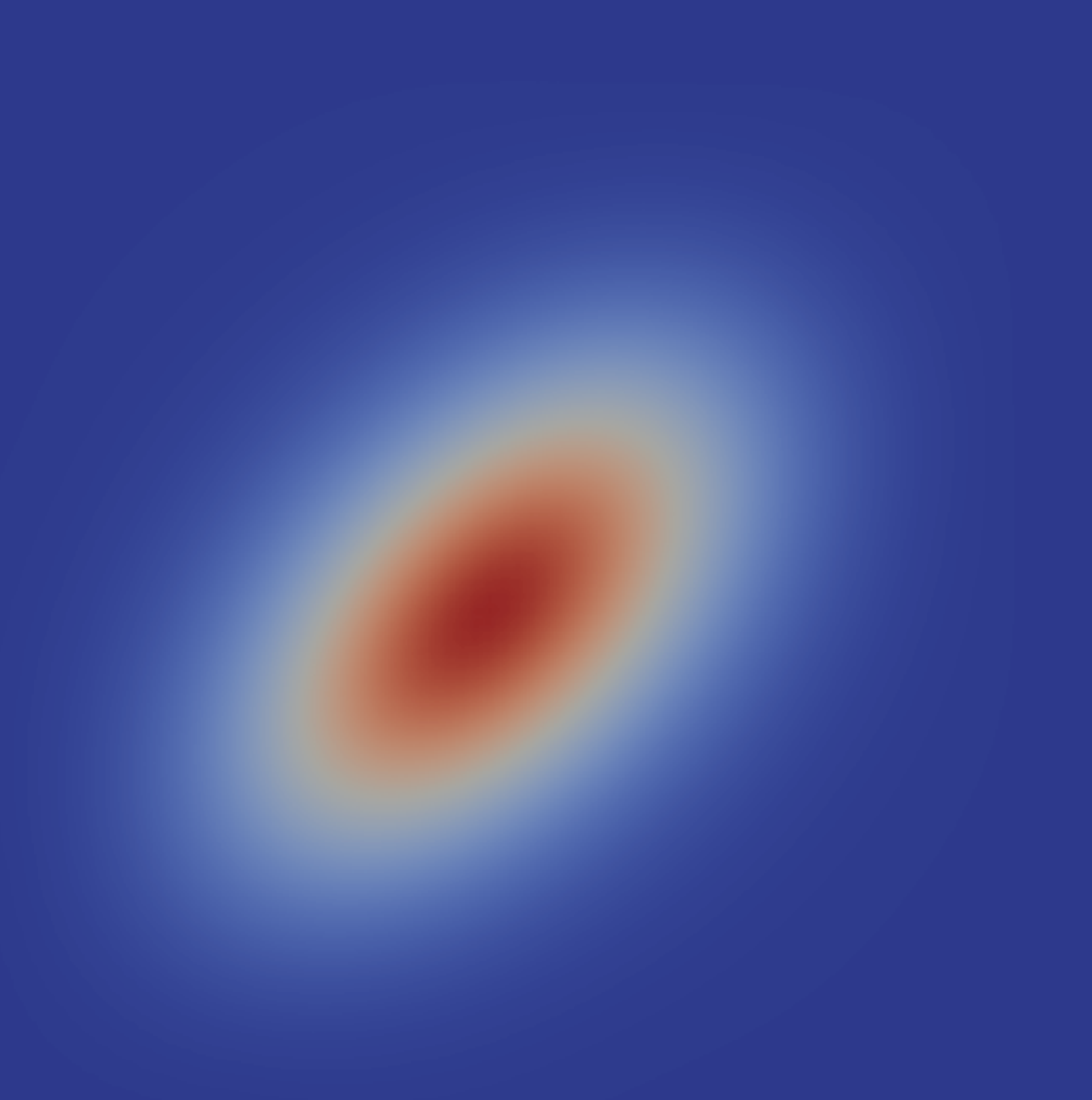}
\includegraphics[width=.19\textwidth]{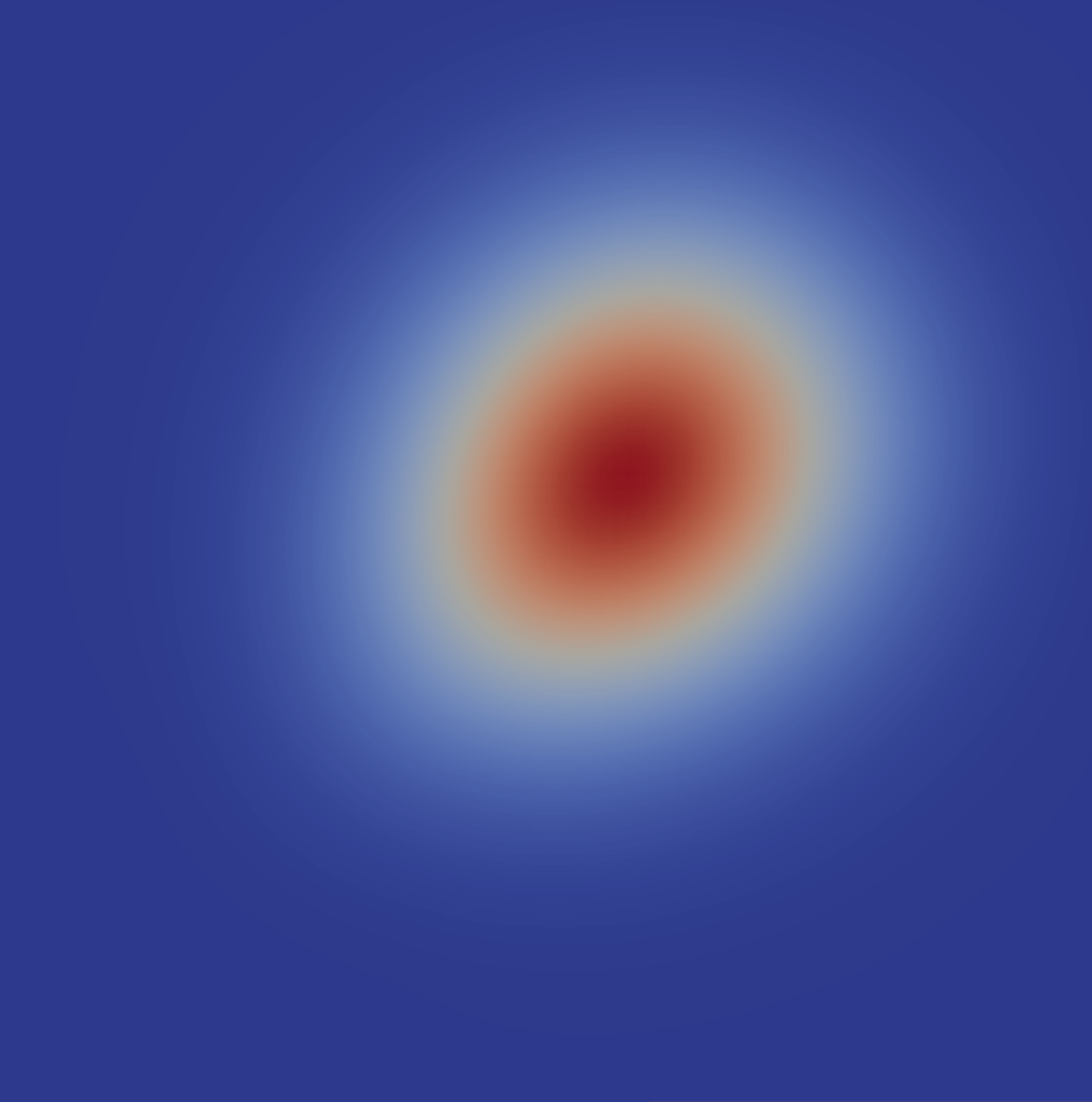}
\includegraphics[width=.19\textwidth]{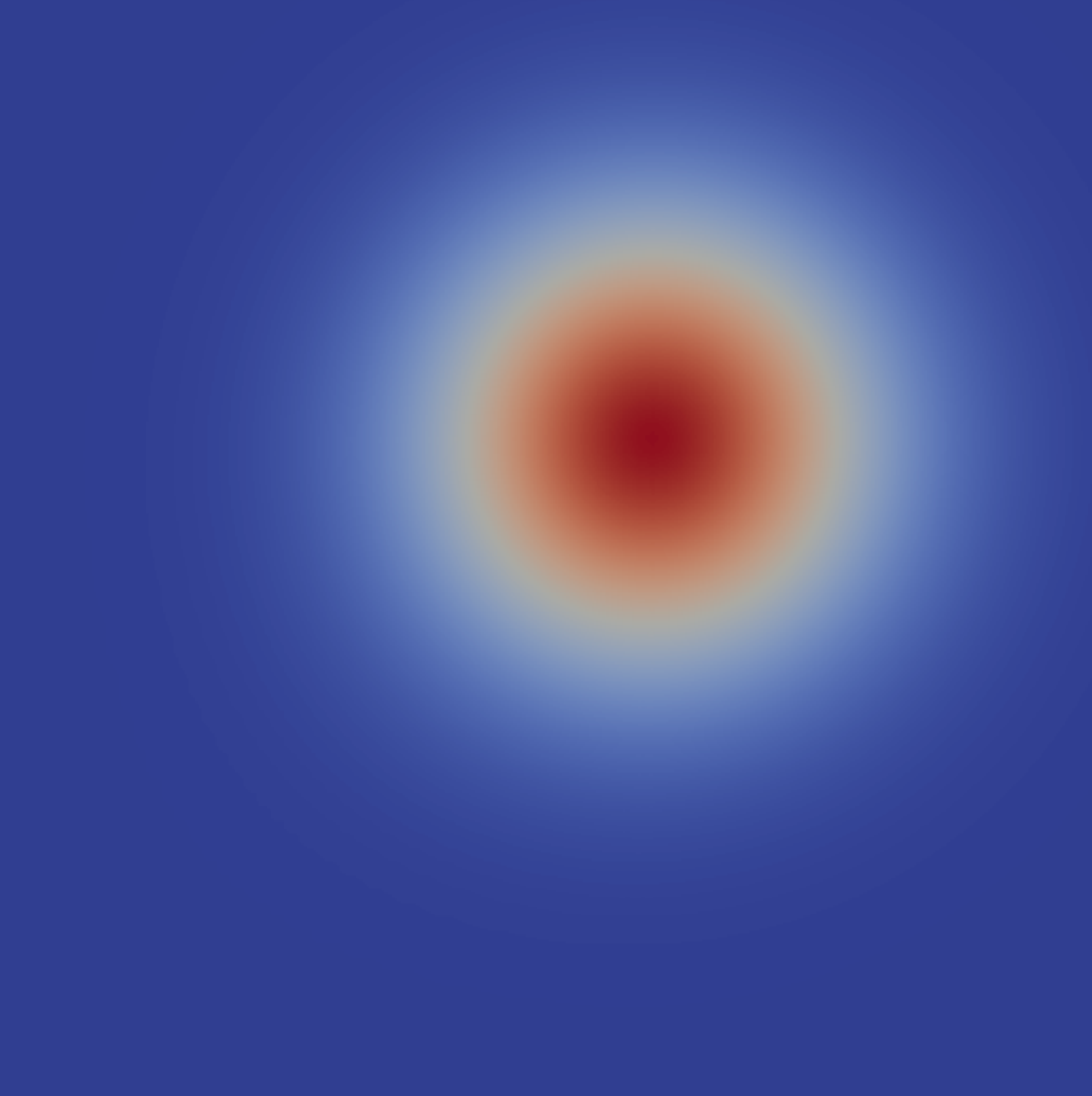}\\[0.07cm]
\includegraphics[width=.19\textwidth]{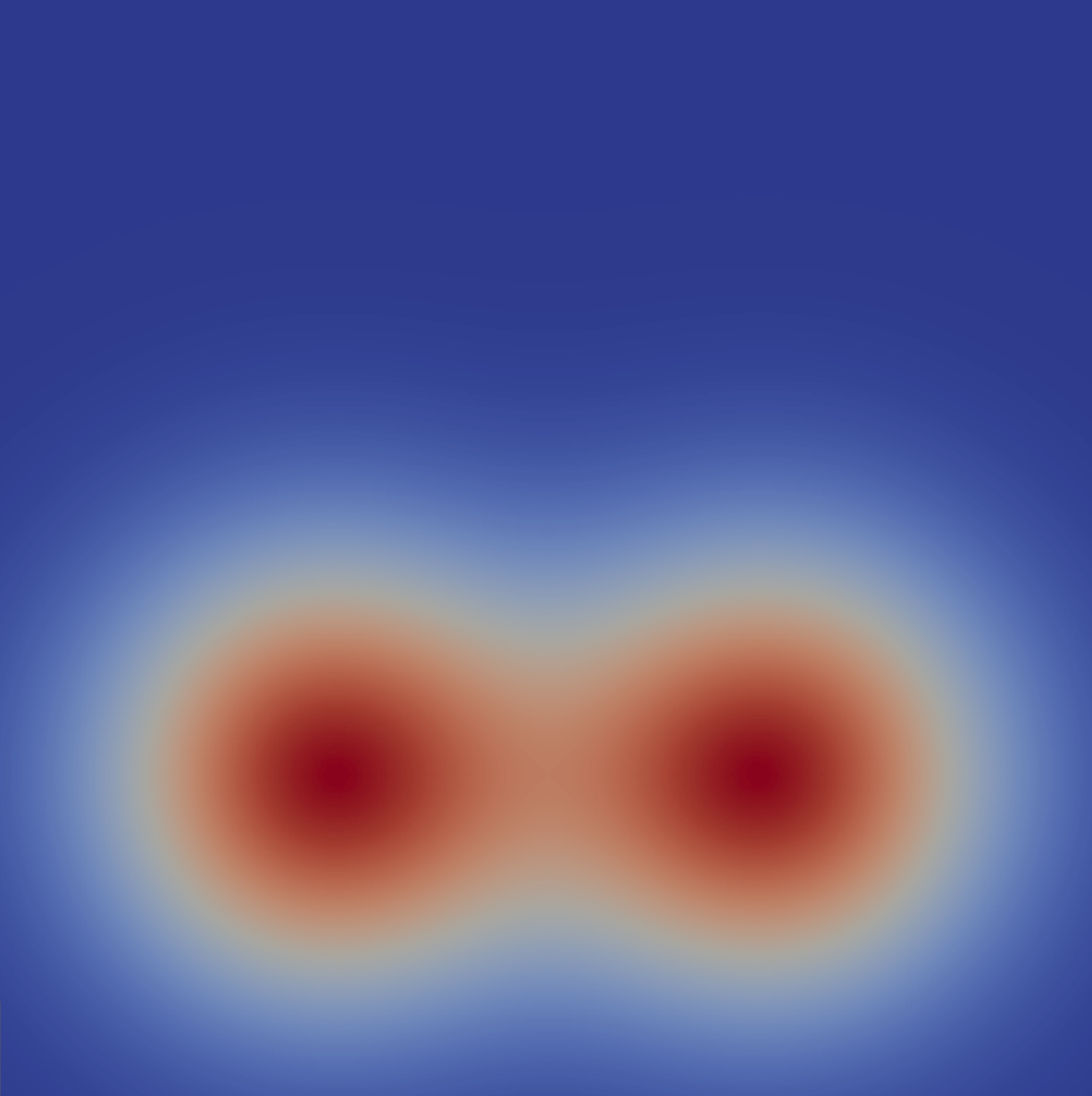}
\includegraphics[width=.19\textwidth]{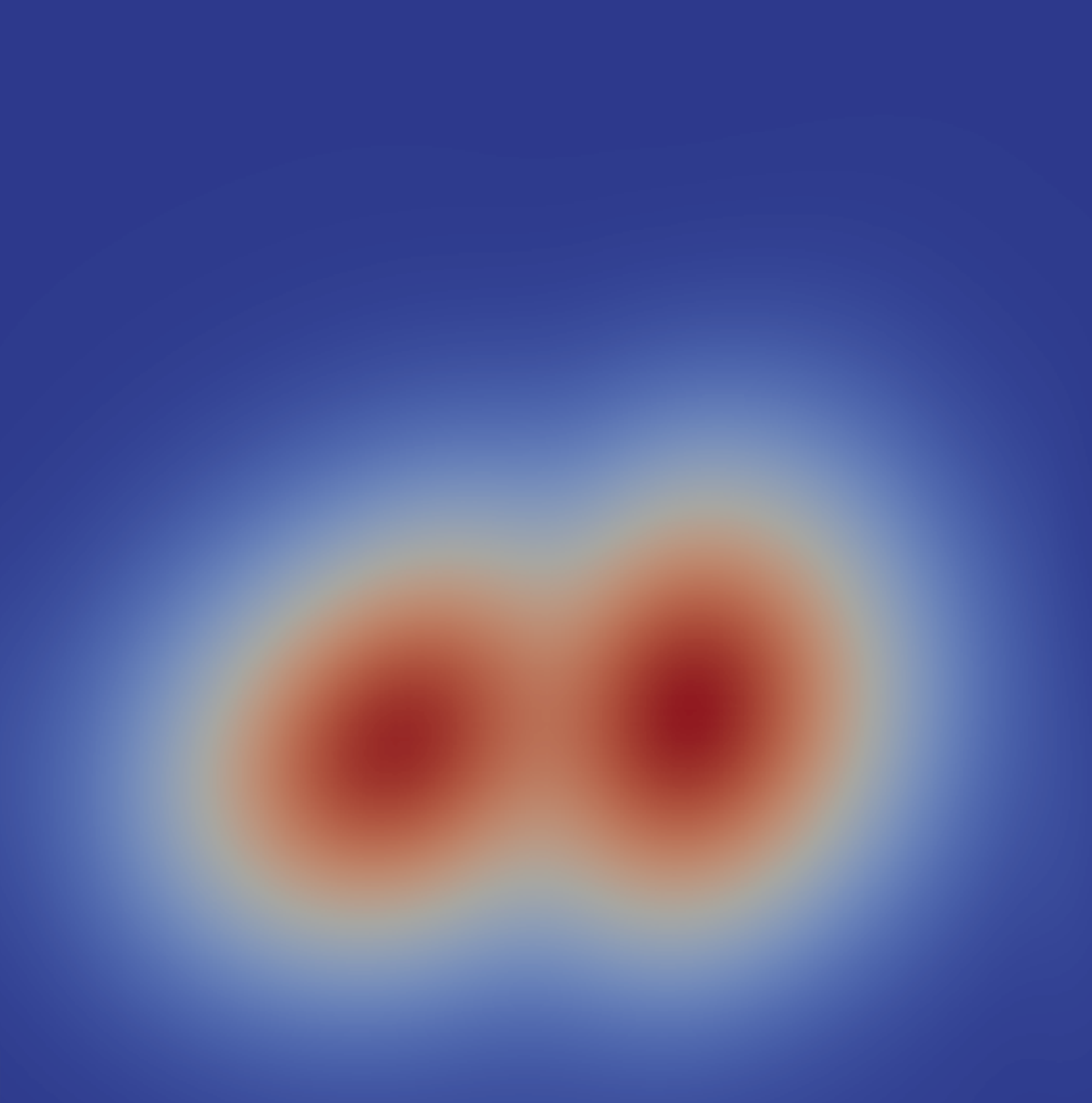}
\includegraphics[width=.19\textwidth]{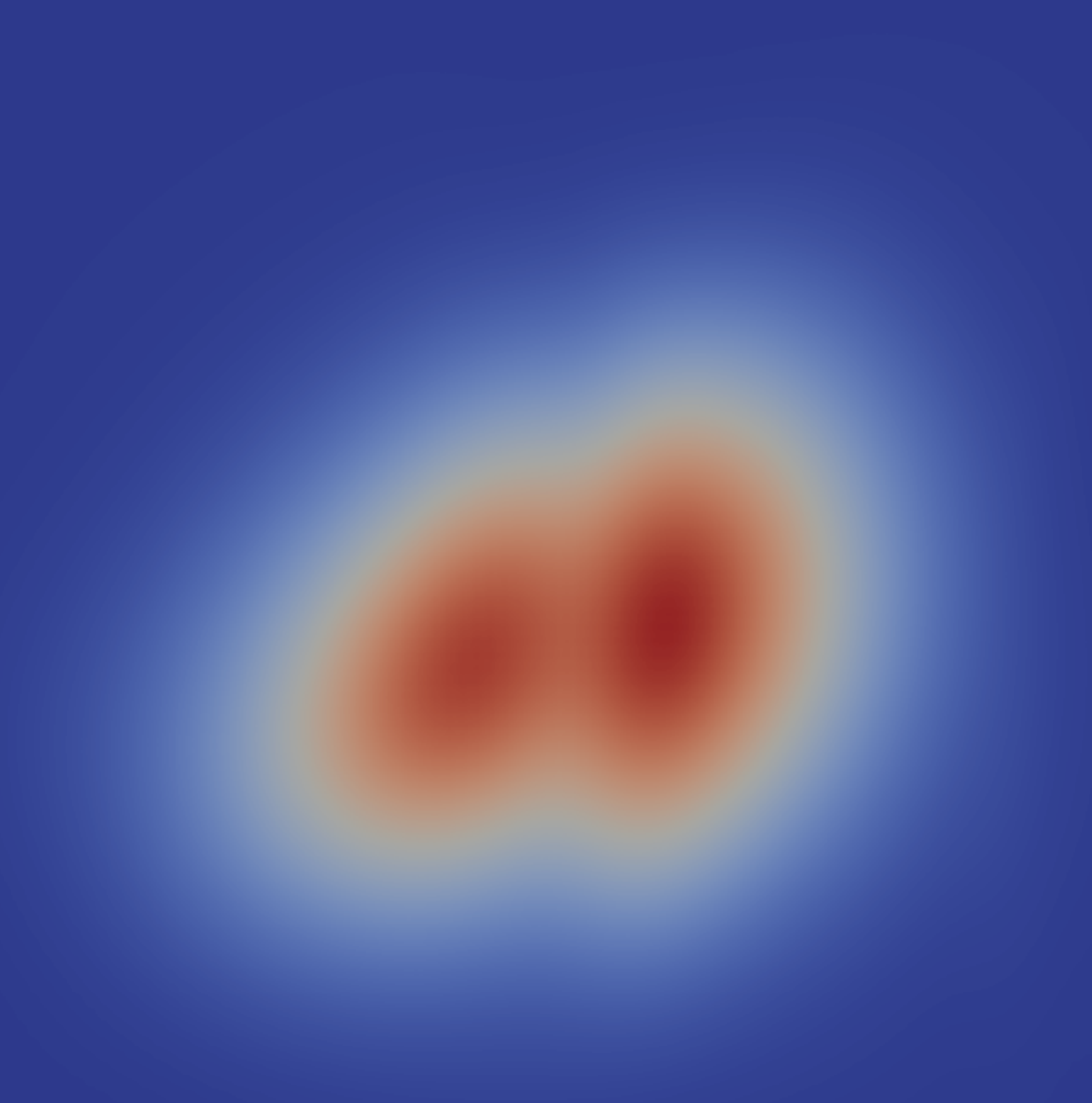}
\includegraphics[width=.19\textwidth]{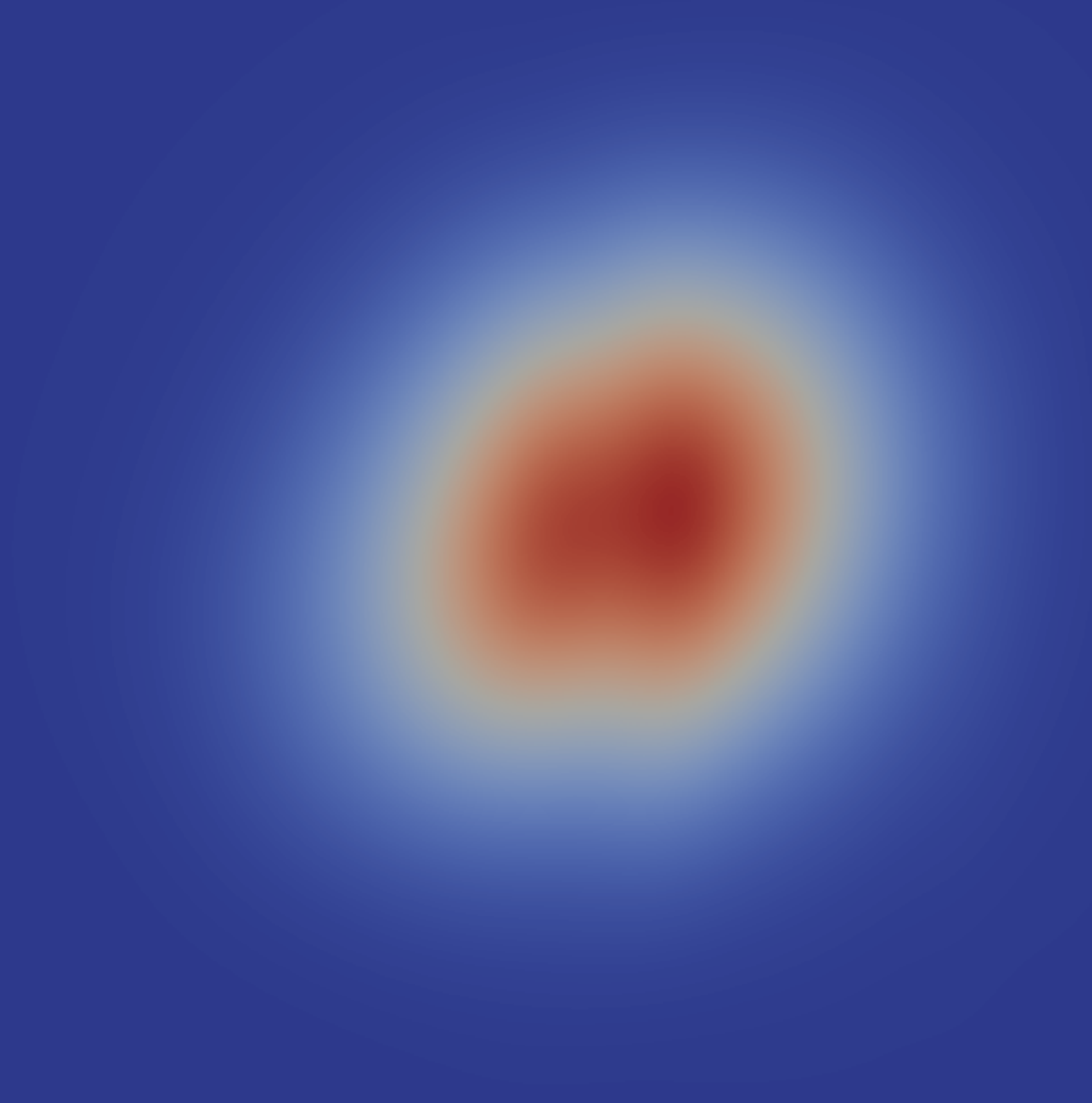}
\includegraphics[width=.19\textwidth]{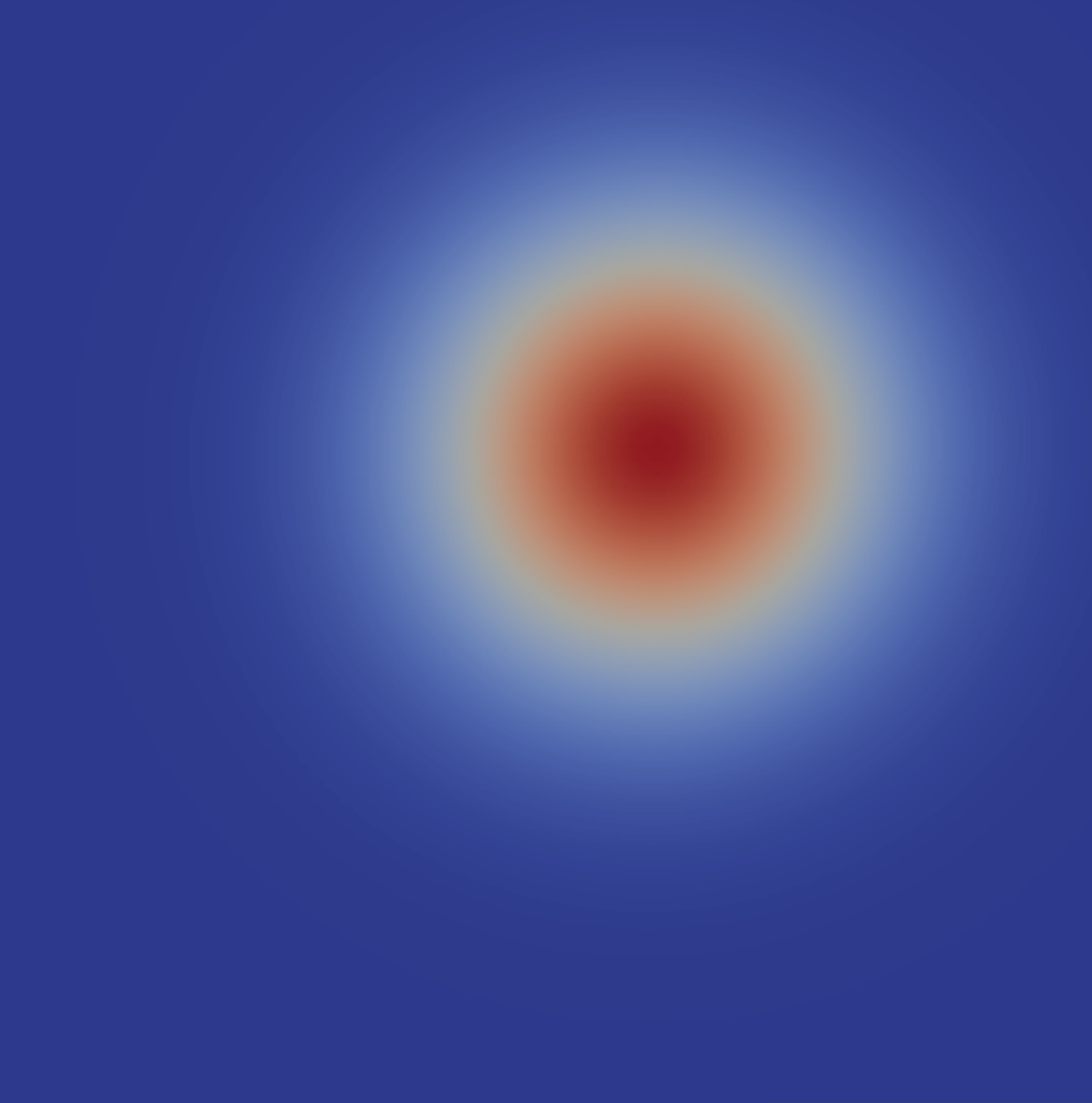}\\[0.07cm]
\includegraphics[width=.19\textwidth]{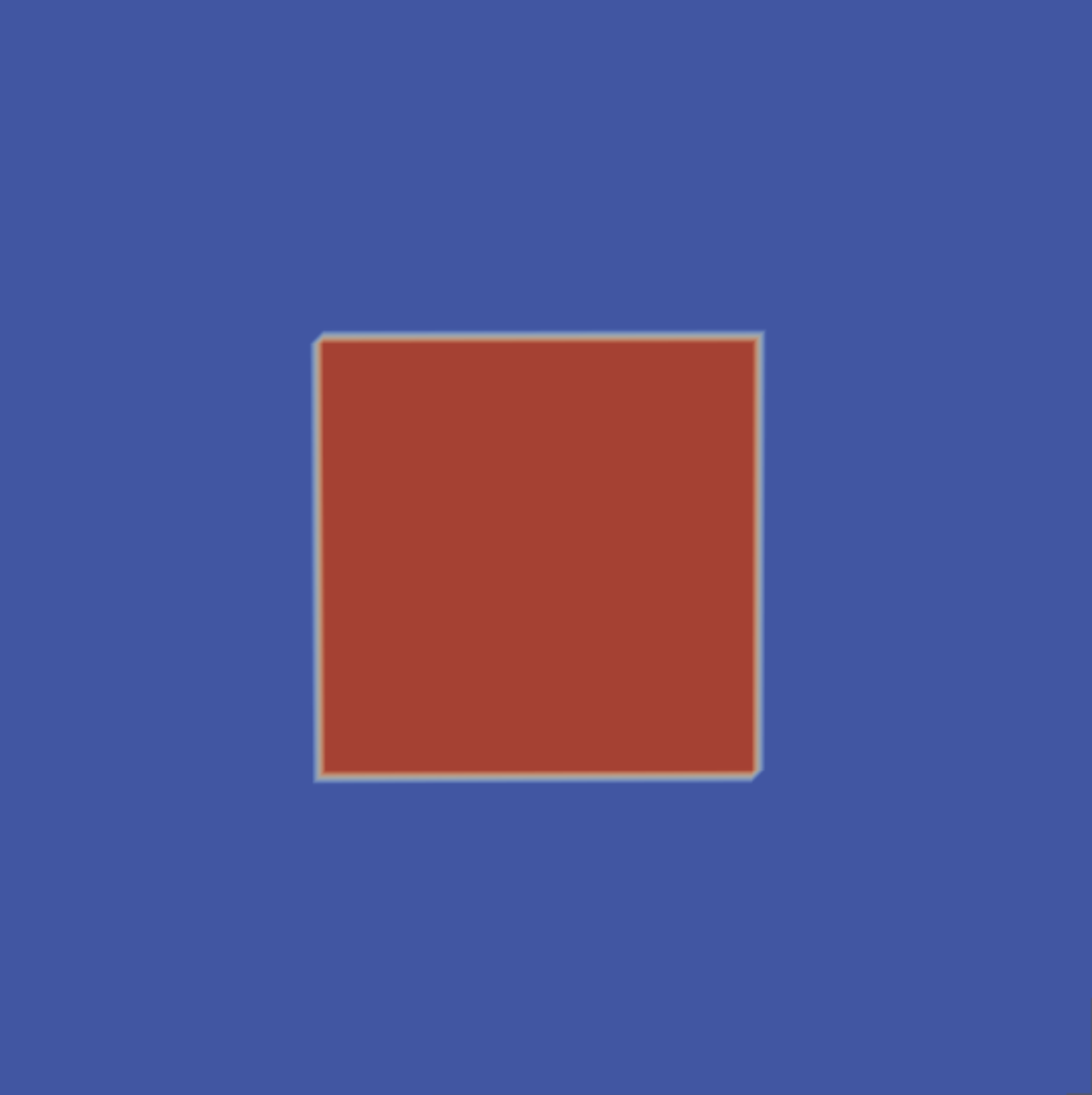}
\includegraphics[width=.19\textwidth]{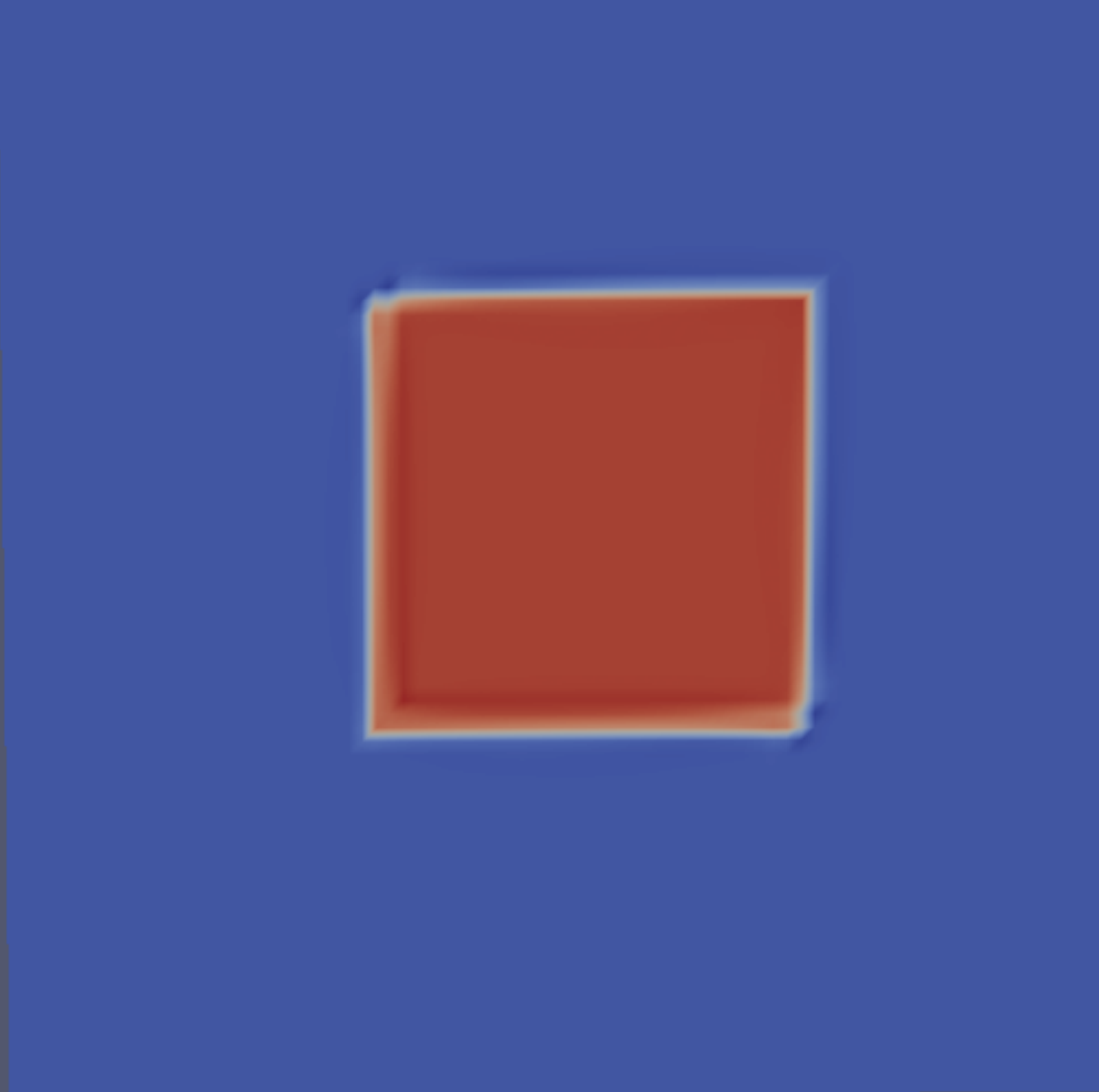}
\includegraphics[width=.19\textwidth]{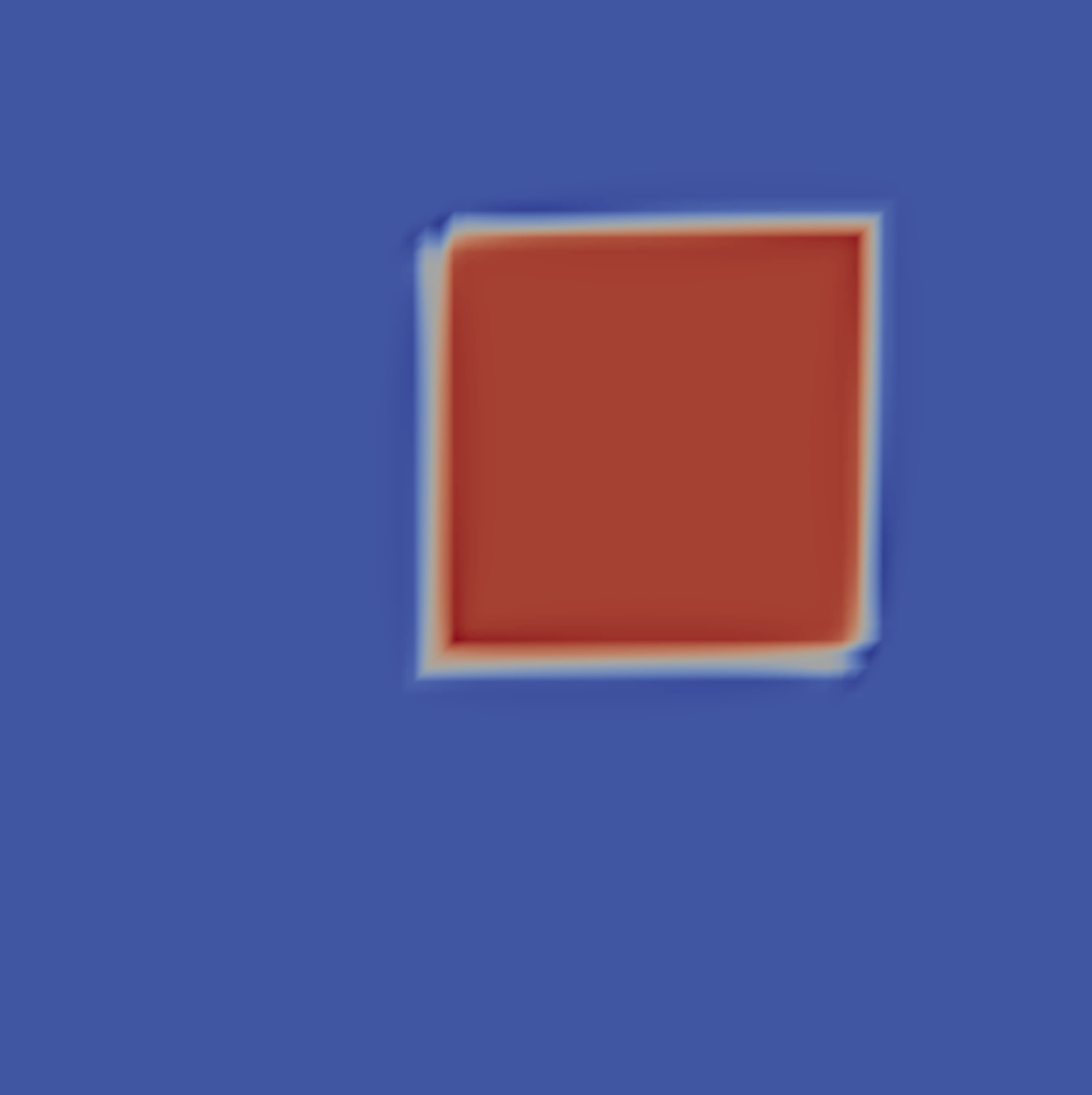}
\includegraphics[width=.19\textwidth]{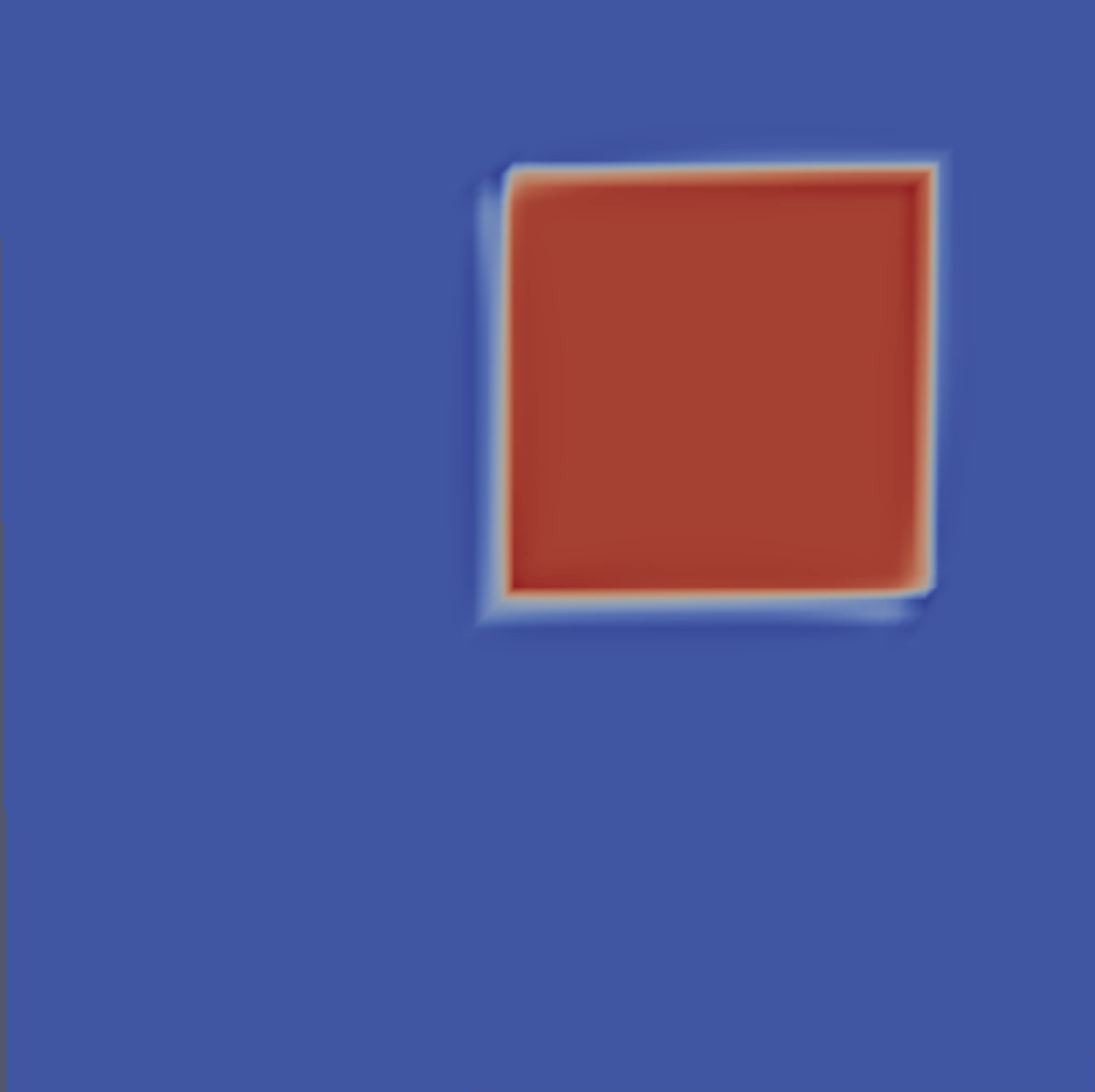}
\includegraphics[width=.19\textwidth]{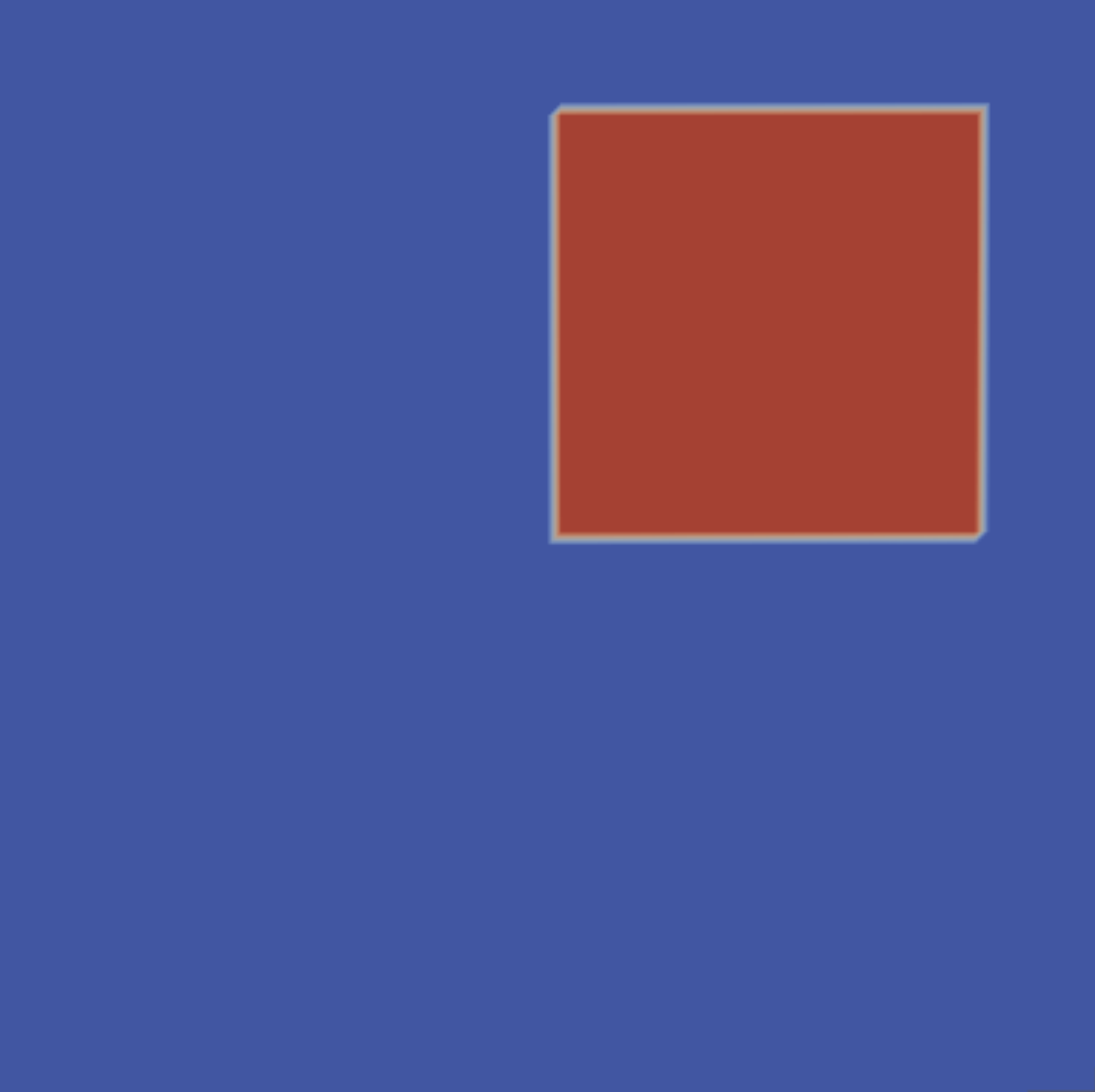}\\[0.07cm]
\includegraphics[width=.19\textwidth]{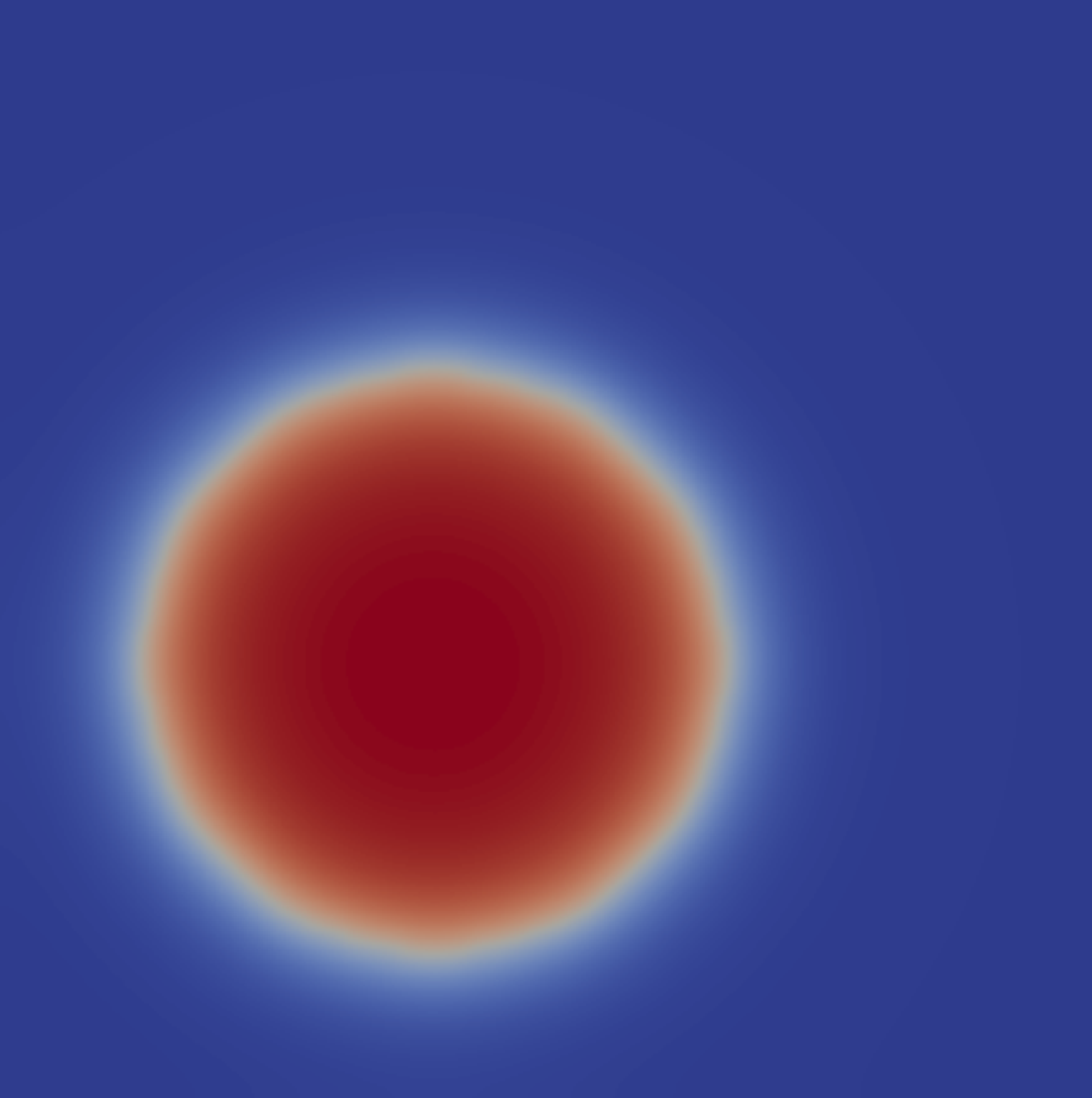}
\includegraphics[width=.19\textwidth]{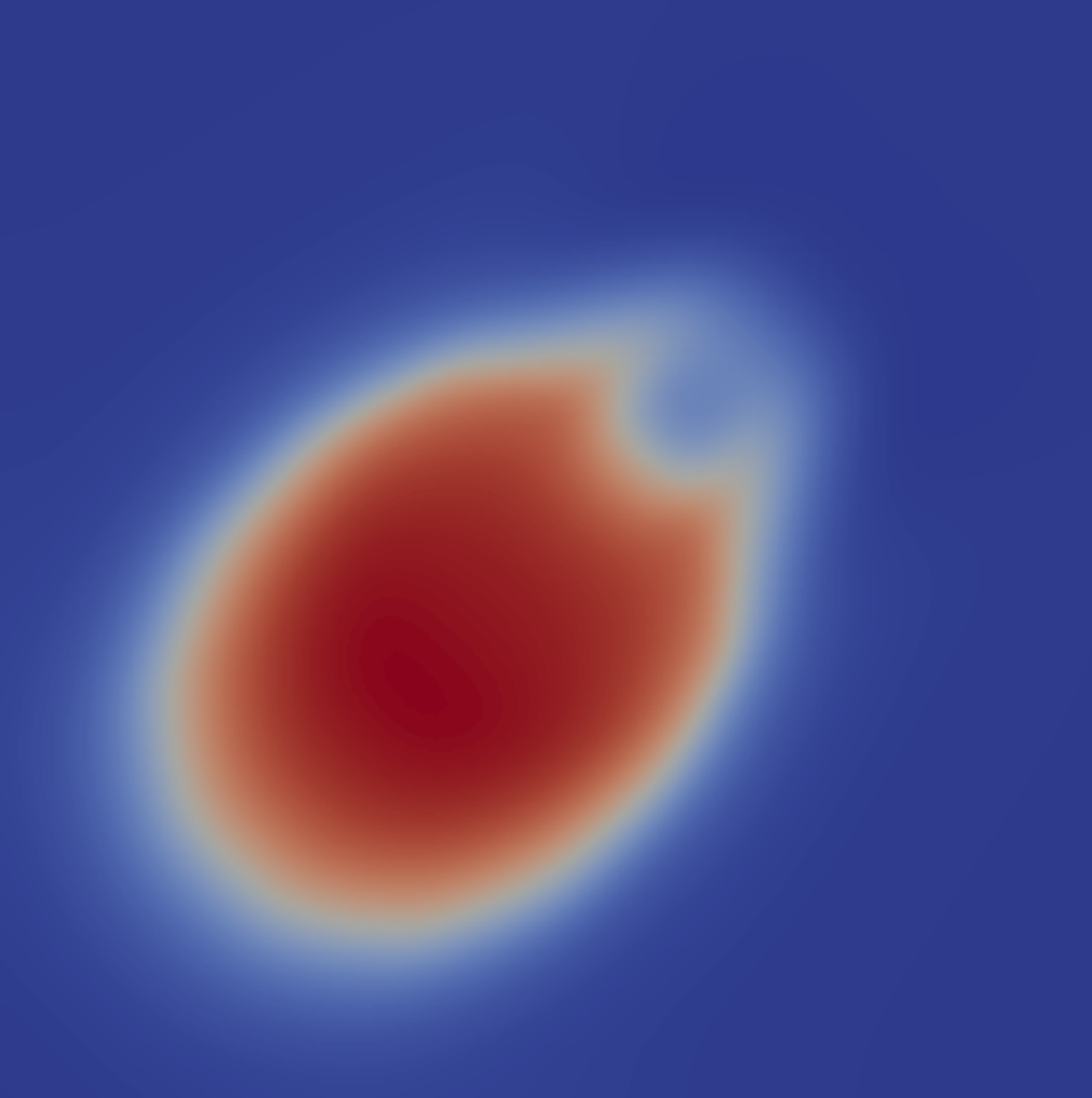}
\includegraphics[width=.19\textwidth]{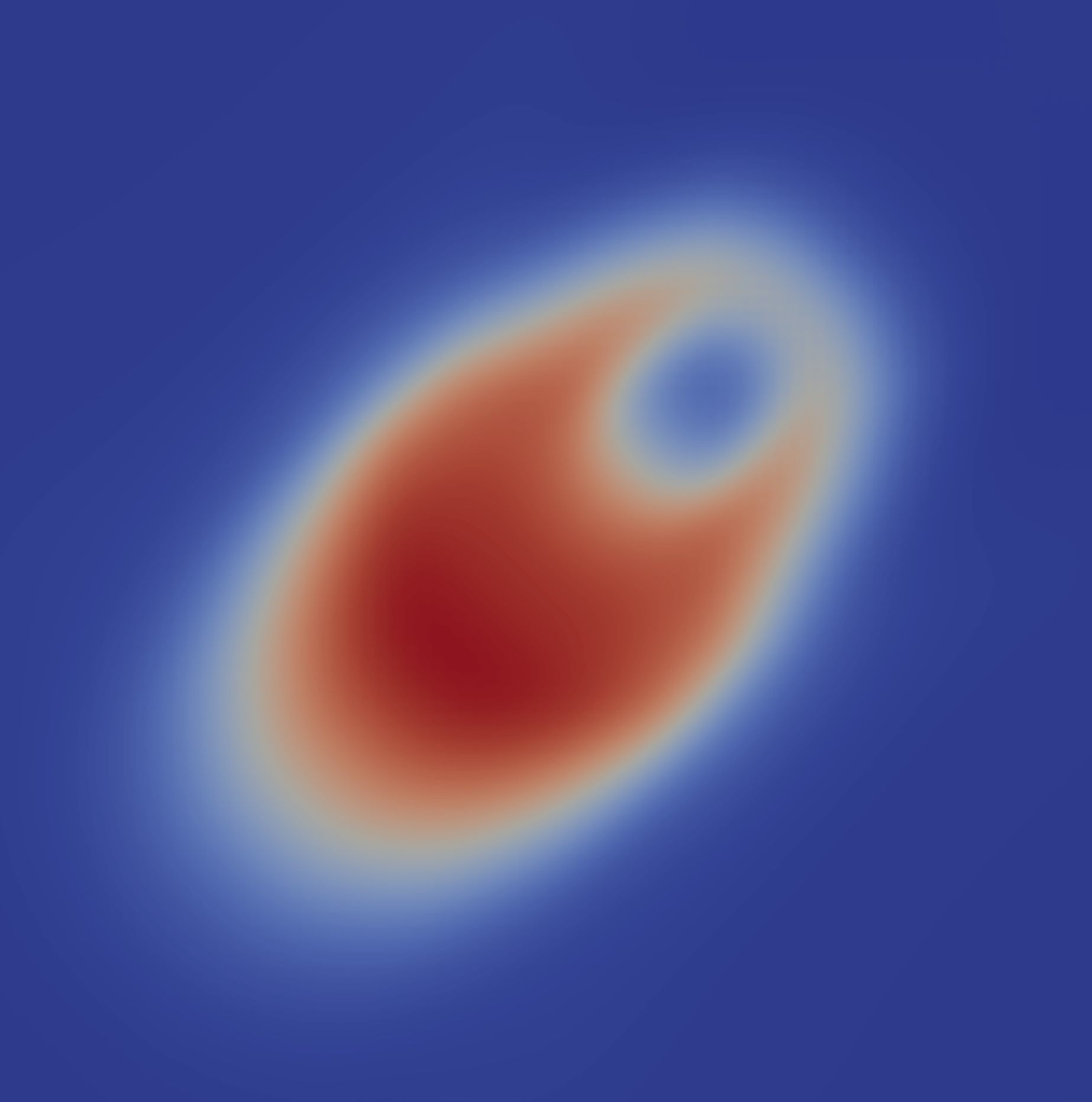}
\includegraphics[width=.19\textwidth]{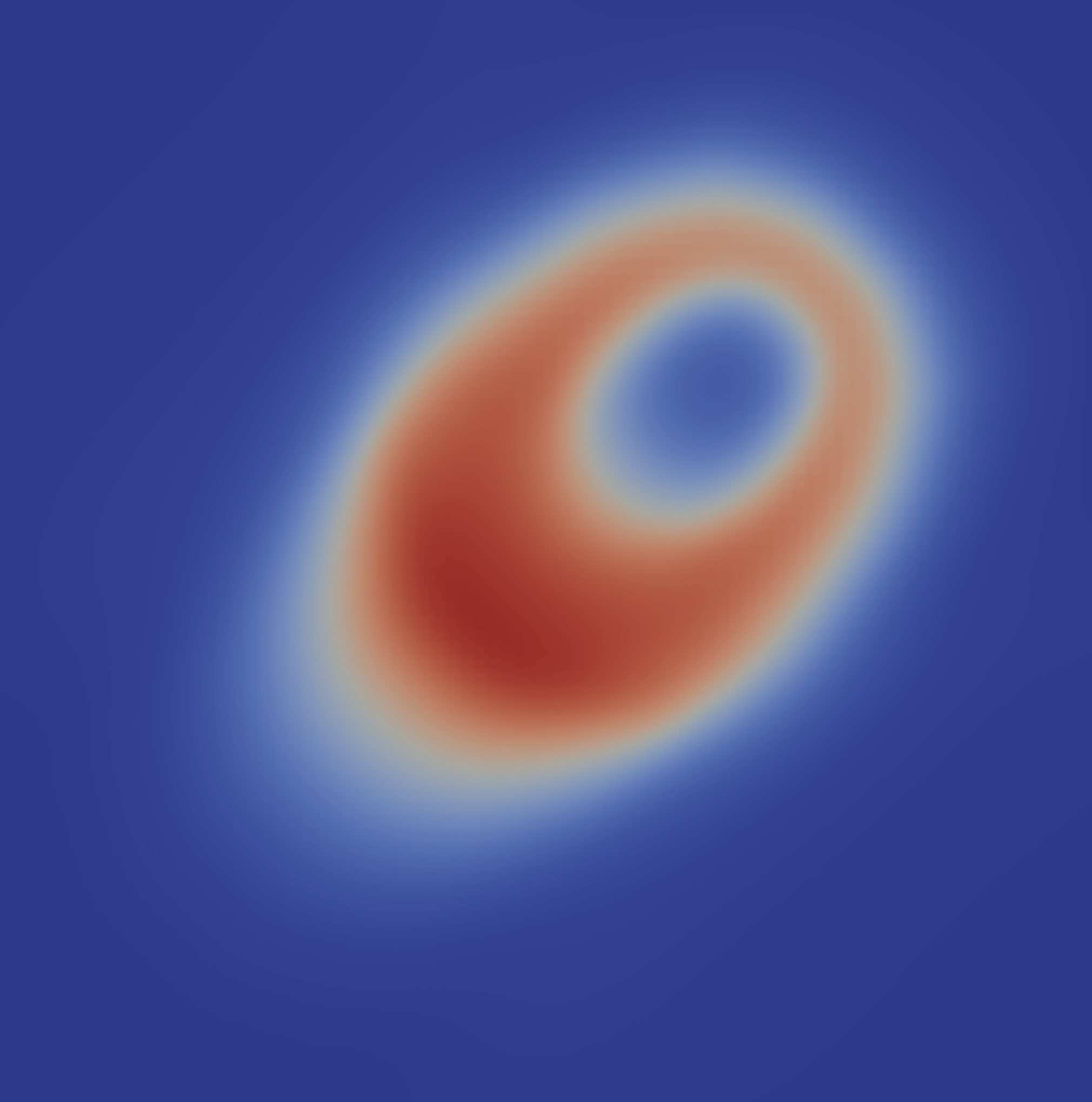}
\includegraphics[width=.19\textwidth]{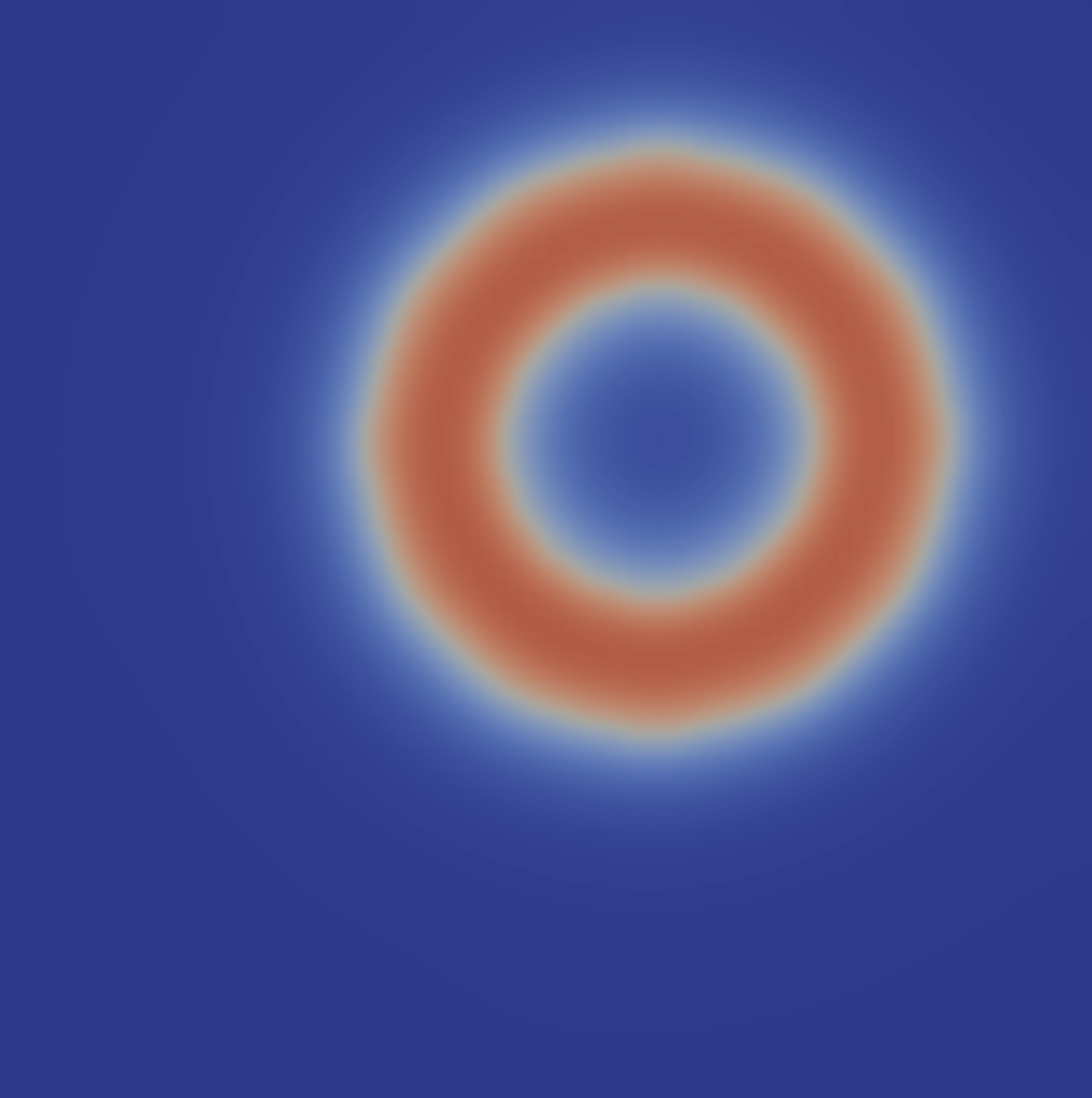}\\[0.07cm]
\includegraphics[width=.19\textwidth]{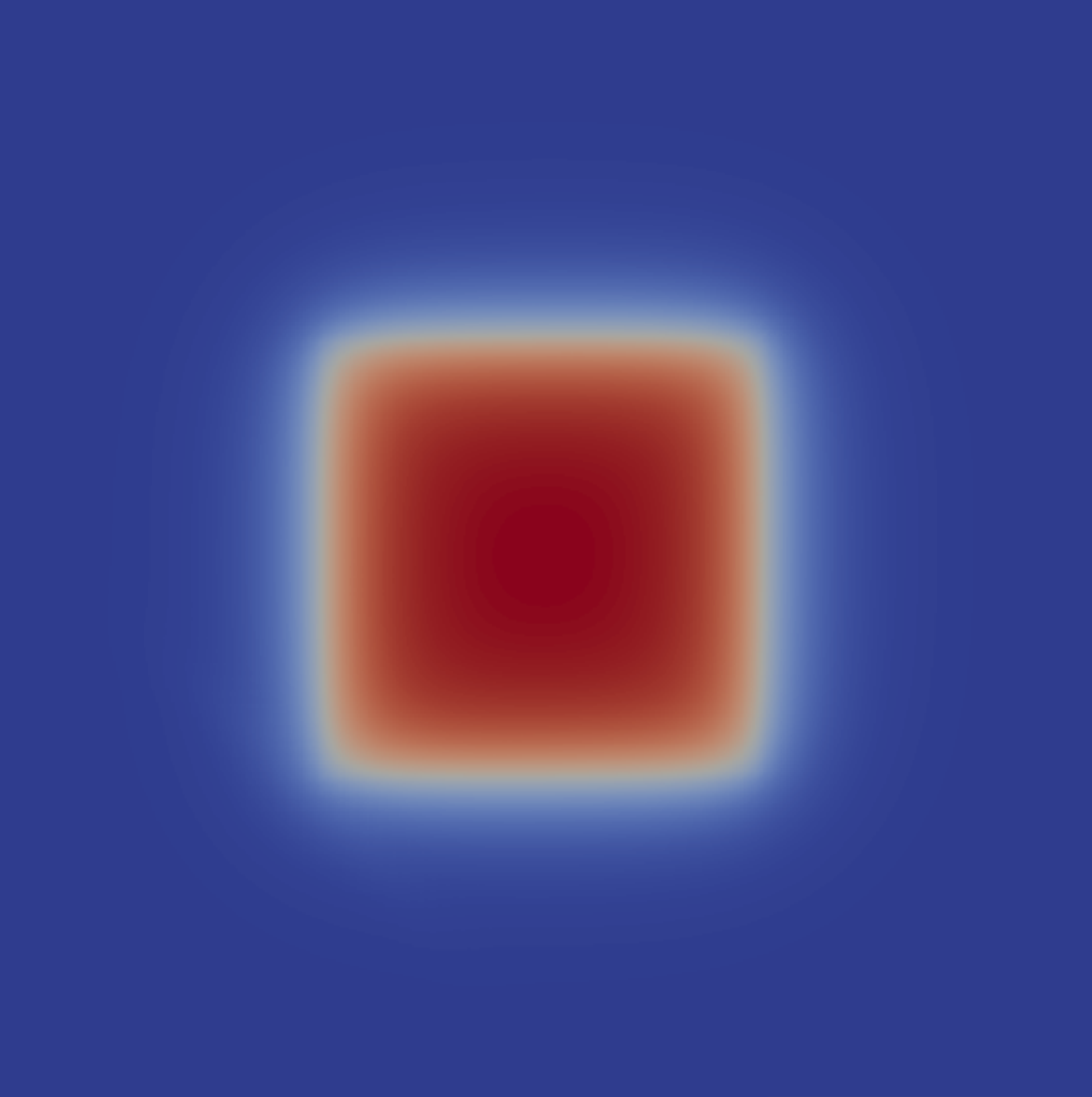}
\includegraphics[width=.19\textwidth]{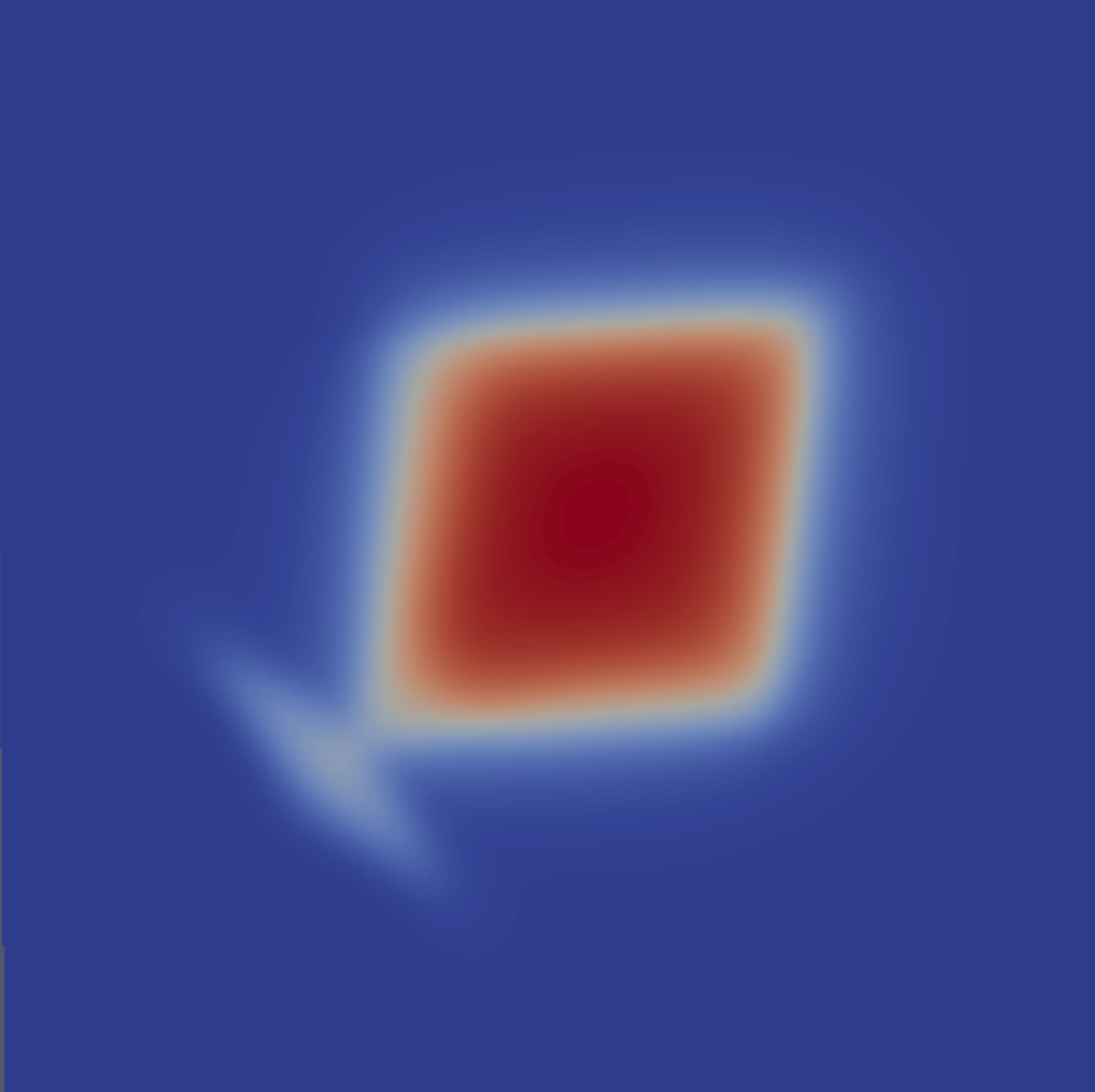}
\includegraphics[width=.19\textwidth]{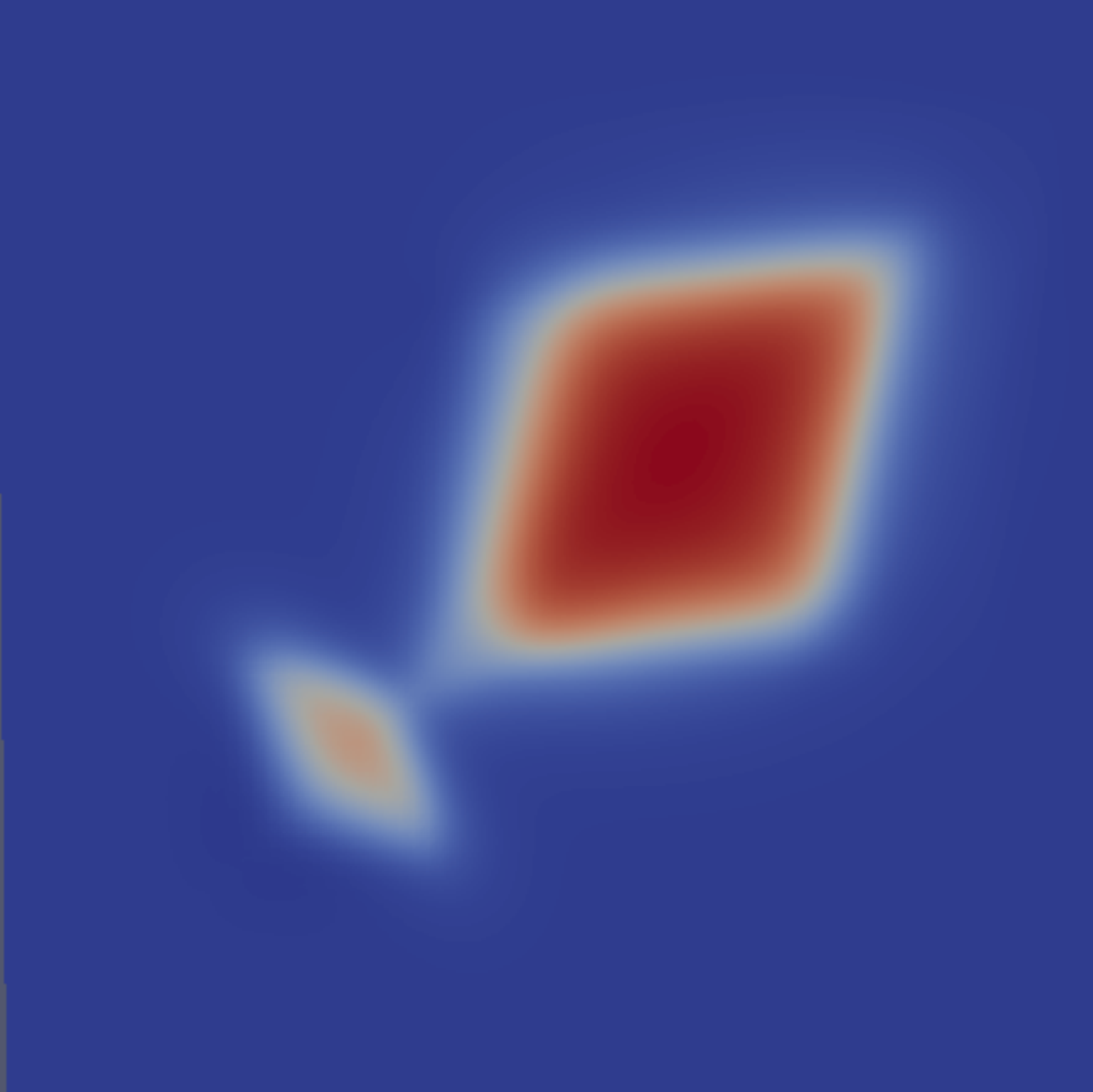}
\includegraphics[width=.19\textwidth]{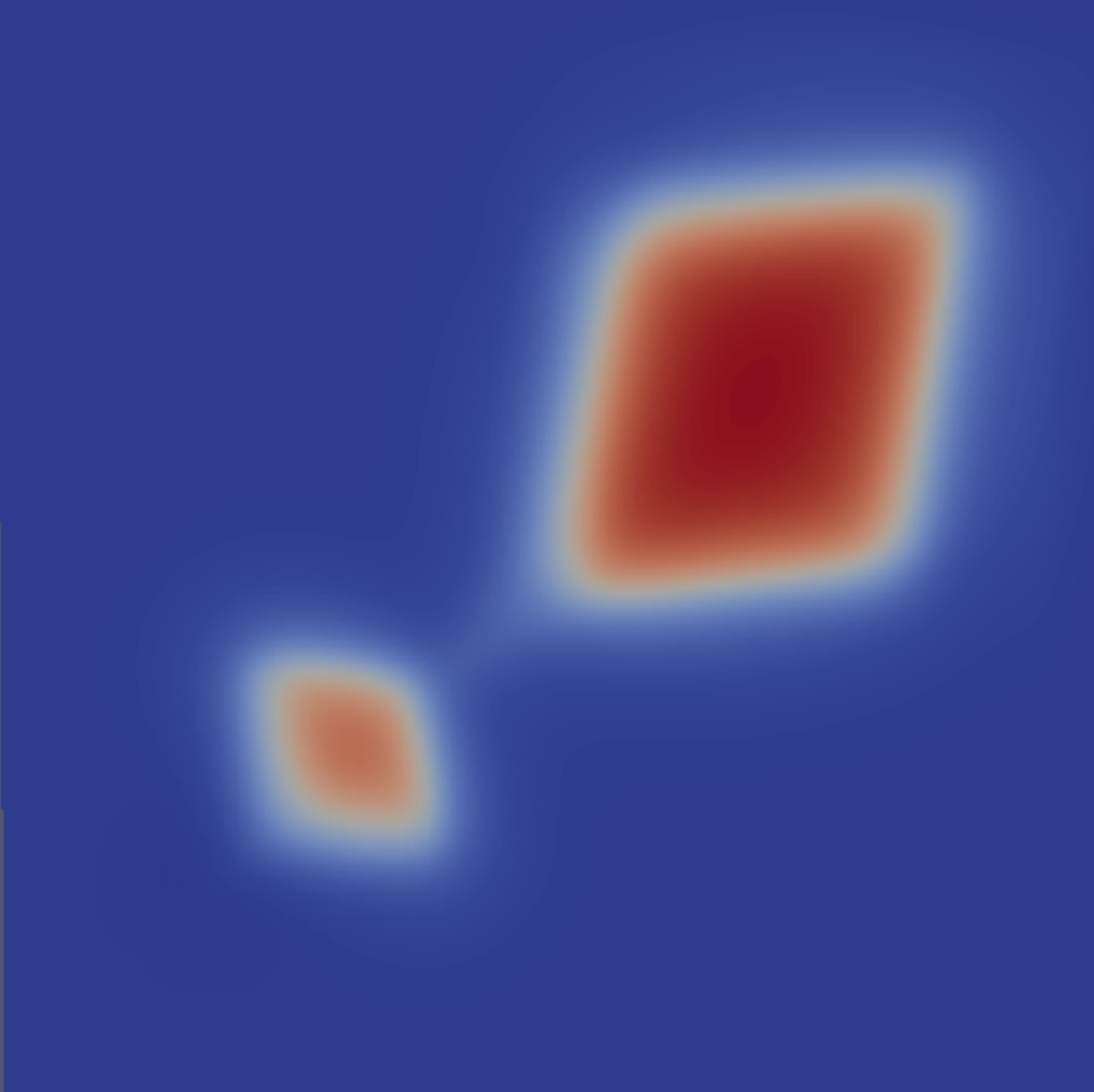}
\includegraphics[width=.19\textwidth]{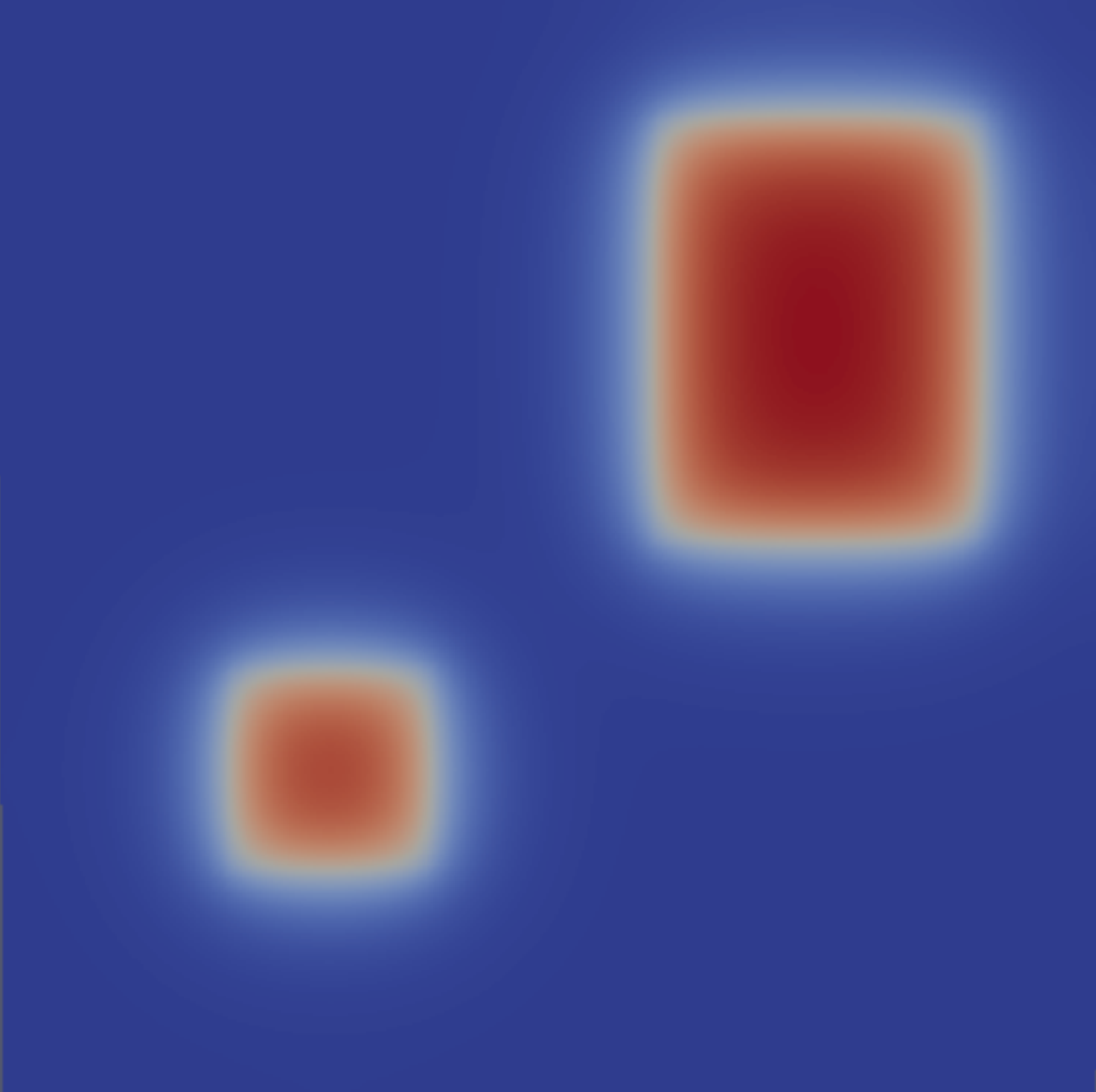}
\end{minipage}
\caption[Results for the outer problem for various template and target pairs in
two dimensions using $\sigma^{-2}=10^5$.]{Results for the outer problem for
various template and target pairs in two dimensions for $\sigma^{-2}=10^5$ to
seek more diffeomorphic geodesics. Each row shows the evolution at time
$t=0,0.25,0.5,0.75,1$.}
\label{fig:adjoint2d}
\end{figure}

We also present some results for $d=2$ in Fig. \ref{fig:adjoint2d_c=1} and
\ref{fig:adjoint2d} for different values of $\sigma$. Here we use the boundary
conditions $u_t|_{\partial\Gamp} = 0$, $\forall t\in [0,1]$, meaning that there
is no inflow or outflow boundaries in space. In the spatial dimension,
$h=0.0125$ and extruded in time with $15$ subdivisions. The observations are
similar to those in the one-dimensional case. We comment on a few of these. For
transport-type problems (the first and third rows of Figs.
\ref{fig:adjoint2d_c=1}--\ref{fig:adjoint2d}), the parameter $\sigma$ does
not change the nature of the solutions, albeit smearing out the smooth data in
the first row when the penalty parameter is low. This behaviour is not
present for discontinuous data which we attribute the flow being almost
grid-aligned and constant in time. $\sigma$ does not appear to alter the
\emph{merging} effect in the example provided in the second row of the figures,
though overall the intensity values are lower whenever $\sigma=1$. Again this
can be explained by advection being enforced to a lesser extent, but not
changing the direction of the velocity field. Lastly are the fourth and fifth
rows where we truly see the influence of $\sigma$. When $\sigma=1$ in Fig.
\ref{fig:adjoint2d_c=1}, very little advection is observed and the motion
between the template and the target can be viewed as a simple "fading in/out"
effect. In Fig. \ref{fig:adjoint2d}, however, we see some interesting results
showing more meaningful (qualitatively speaking), geodesics between the images.
Although these results are still of a toy problem nature they show clear promise
for our approach. These results can be generated in less than 24 hours on
standard consumer hardware (without using any preconditioning whatsoever).
Applying high performance parallel algorithms are the subject of future
research.

\section{Summary}\label{outer:conclusion}

We successfully cast the metamorphosis problem in a variational Hilbertian
space-time setting. The variational setting introduced in section \ref{inner}
provides a novel framework to study weakly space-time $\bgrad$-differentiable
geodesics for metamorphosis. A \emph{primal} least-squared formulation was
studied and we showed several theoretical results for the inner problem. Based
on this we presented a simple conforming finite element method. Section
\ref{outer_problem} showed a practical adjoint-based way to solve the inverse
problem to find a suitable velocity field $u$ requiring minimal implementation
effort by the user.\\

There are many open problems and options yet to be explored for space-time
metamorphosis. The main theoretical challenges for the inner problem is showing
a regularity result akin to the usual elliptic regularity for $H^1$ variational
problems. Preconditioning of the least-squared system in section \ref{inner:lsq}
is another option for future research, as is the investigation of weaker notions
of discretizations e.g. using the discontinuous Petrov-Galerkin
 framework, see \cite{demkowicz2010class}. The work of \cite{anderson2007numerical} offers some
insights into preconditioning for a space-time wave equation and offers
excellent numerical evidence that might be useful in this endeavour. Space-time
adaptivity could also be explored in order to resolve details of images
represented as finite element functions, in which case one would dispense with a
quasi-uniformity assumption on the mesh. In view of our numerical experiments 
for $d=2$ we believe this paper provides the basis for a powerful method given
good preconditioners and more bespoke numerical optimizers.\\ The outer problem
also presents some possibilities for future work. As mentioned, a conforming
finite element discretization of $V$ is desirable, although a
$\mathsf{C}^0$ conforming $H^3$ non-conforming finite element exists and could be
explored, see \cite{wu2019nonconforming}. Since the Lipschitz regularity of the
velocity field is only required in space, such an implementation over the
spatial domain could be combined with extruded quadrilateral meshes
 to impose different temporal and
spatial regularity, see \cite{mcrae2016automated,bercea2016structure}.\\
In this project we have chosen to treat the outer problem with a
\emph{discretize-then-optimize} approach because of the availability of the
seamless adjoint framework in Firedrake.  Another approach that was attempted
was to derive the Euler-Lagrange equations for the "full" metamorphosis problem
(also called EPDiff in the literature) to solve for $u$, $I$ and $z$
\emph{simultaneously}, and then discretize this non-linear system of equations. A
Newton method could then be applied for this system. Some work has been done
previously in this area for LDDMM, see \cite{hernandez2008gauss}.  Also, the work
\cite{zhang2015finite} establishes a fast algorithm for computing
finite-dimensional velocity fields which could be used in conjunction with the
inner problem formulations above.  This may be morally similar to discretizing
the velocity field using finite elements but we simply mention this as a
possible avenue owing to the numerical evidence in support of the authors' work.



\bibliographystyle{abbrv}
\bibliography{ms}

\clearpage

\appendix
\newpage

\section*{Appendix A.\ Inner Problem Convergence Results for $d=1$}\label{app:innerconv1d}

\begin{table}[!ht]
\centering
\caption{Convergence rates for example 0 from Table {\ref{tab:cg_manufactured}}}
{\begin{tabular}{c|cc|cc}
\hline\bfseries $h\inv$ & $L^2$ error & Order & Energy error & Order
\begin{filecontents*}{imgs/inner_problem/test0_rates_1d.csv}
hmesh,ltwoerr,ltwoord,energyerr,energyord
4,2.4669e-01,-,1.1330e+00,-
16,9.3254e-02,0.70,1.9597e+00,-0.40
36,8.2272e-02,0.15,3.0075e+00,-0.53
64,7.1942e-02,0.23,4.0880e+00,-0.53
100,7.1870e-02,0.00,5.1363e+00,-0.51
144,6.0582e-02,0.47,6.1323e+00,-0.49
196,5.7787e-02,0.15,7.1294e+00,-0.49
256,5.4580e-02,0.21,8.1551e+00,-0.50
324,5.2206e-02,0.19,9.1998e+00,-0.51
400,5.1512e-02,0.06,1.0236e+01,-0.51
484,4.8312e-02,0.34,1.1245e+01,-0.49
576,4.6738e-02,0.19,1.2253e+01,-0.49
\end{filecontents*}
\csvreader[head to column names]{imgs/inner_problem/test0_rates_1d.csv}{}
{\\\hline \hmesh & \ltwoerr &\ltwoord & \energyerr & \energyord}\\\hline
\end{tabular}}
\label{tab:cg_conv0}
\end{table}

\begin{table}[!ht]
\centering
\caption{Convergence rates for example 1 from Table {\ref{tab:cg_manufactured}}}
{\begin{tabular}{c|cc|cc}
\hline\bfseries $h\inv$ & $L^2$ error & Order & Energy error & Order
\begin{filecontents*}{imgs/inner_problem/test1_rates_1d.csv}
hmesh,ltwoerr,ltwoord,energyerr,energyord
4,1.2943e-01,-,8.1847e-01,-
16,6.4739e-02,0.50,4.4521e-01,0.44
36,3.0632e-02,0.92,1.5345e-01,1.31
64,1.3834e-02,1.38,5.9659e-02,1.64
100,6.6335e-03,1.65,2.6630e-02,1.81
144,3.4437e-03,1.80,1.3348e-02,1.89
196,1.9280e-03,1.88,7.3401e-03,1.94
256,1.1521e-03,1.93,4.3442e-03,1.96
324,7.2697e-04,1.95,2.7265e-03,1.98
400,4.7997e-04,1.97,1.7944e-03,1.99
484,3.2910e-04,1.98,1.2279e-03,1.99
576,2.3294e-04,1.99,8.6803e-04,1.99
\end{filecontents*}
\csvreader[head to column names]{imgs/inner_problem/test1_rates_1d.csv}{}
{\\\hline \hmesh & \ltwoerr &\ltwoord & \energyerr & \energyord}\\\hline
\end{tabular}}
\label{tab:cg_conv1}
\end{table}

\begin{table}[!ht]
\centering
\caption{Convergence rates for example 2 from Table {\ref{tab:cg_manufactured}}}
{\begin{tabular}{c|cc|cc}
\hline\bfseries $h\inv$ & $L^2$ error & Order & Energy error & Order
\begin{filecontents*}{imgs/inner_problem/test2_rates_1d.csv}
hmesh,ltwoerr,ltwoord,energyerr,energyord
4,1.2494e-01,-,1.0894e+00,-
16,5.9254e-02,0.54,3.8319e-01,0.75
36,2.8248e-02,0.91,1.2476e-01,1.38
64,1.2967e-02,1.35,4.8549e-02,1.64
100,6.3042e-03,1.62,2.1807e-02,1.79
144,3.3031e-03,1.77,1.0981e-02,1.88
196,1.8597e-03,1.86,6.0554e-03,1.93
256,1.1150e-03,1.92,3.5898e-03,1.96
324,7.0499e-04,1.95,2.2552e-03,1.97
400,4.6604e-04,1.96,1.4851e-03,1.98
484,3.1979e-04,1.98,1.0166e-03,1.99
576,2.2647e-04,1.98,7.1886e-04,1.99
\end{filecontents*}
\csvreader[head to column names]{imgs/inner_problem/test2_rates_1d.csv}{}
{\\\hline \hmesh & \ltwoerr &\ltwoord & \energyerr & \energyord}\\\hline
\end{tabular}}
\label{tab:cg_conv2}
\end{table}

\begin{table}[!ht]
\centering
\caption{Convergence rates for example 3 from Table {\ref{tab:cg_manufactured}}}
{\begin{tabular}{c|cc|cc}
\hline\bfseries $h\inv$ & $L^2$ error & Order & Energy error & Order
\begin{filecontents*}{imgs/inner_problem/test3_rates_1d.csv}
hmesh,ltwoerr,ltwoord,energyerr,energyord
4,2.1457e-01,-,1.4589e+00,-
16,5.8776e-02,0.93,1.2915e+00,0.09
36,9.8465e-03,2.20,3.6886e-01,1.55
64,3.6522e-03,1.72,1.2592e-01,1.87
100,1.6125e-03,1.83,5.2731e-02,1.95
144,8.1228e-04,1.88,2.5650e-02,1.98
196,4.5136e-04,1.91,1.3901e-02,1.99
256,2.7014e-04,1.92,8.1655e-03,1.99
324,1.7131e-04,1.93,5.1038e-03,1.99
400,1.1378e-04,1.94,3.3511e-03,2.00
484,7.8465e-05,1.95,2.2900e-03,2.00
576,5.5831e-05,1.96,1.6174e-03,2.00
\end{filecontents*}
\csvreader[head to column names]{imgs/inner_problem/test3_rates_1d.csv}{}
{\\\hline \hmesh & \ltwoerr &\ltwoord & \energyerr & \energyord}\\\hline
\end{tabular}}
\label{tab:cg_conv3}
\end{table}

\newpage
\section*{Appendix B.\ Inner Problem Convergence Results for $d=2$}\label{app:innerconv2d}

\begin{figure}[h!]
\centering
\includegraphics[scale=.4]{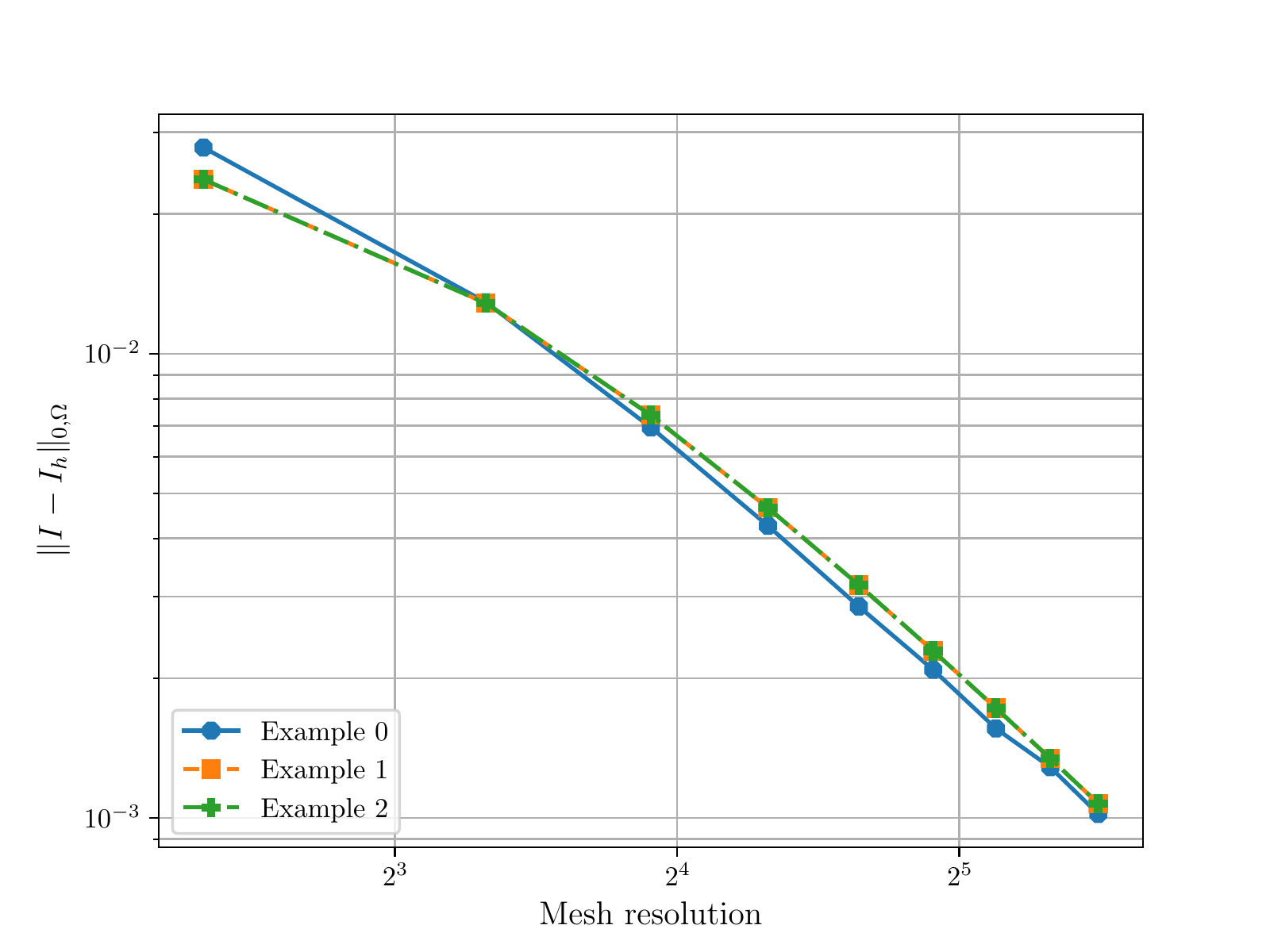}
\includegraphics[scale=.4]{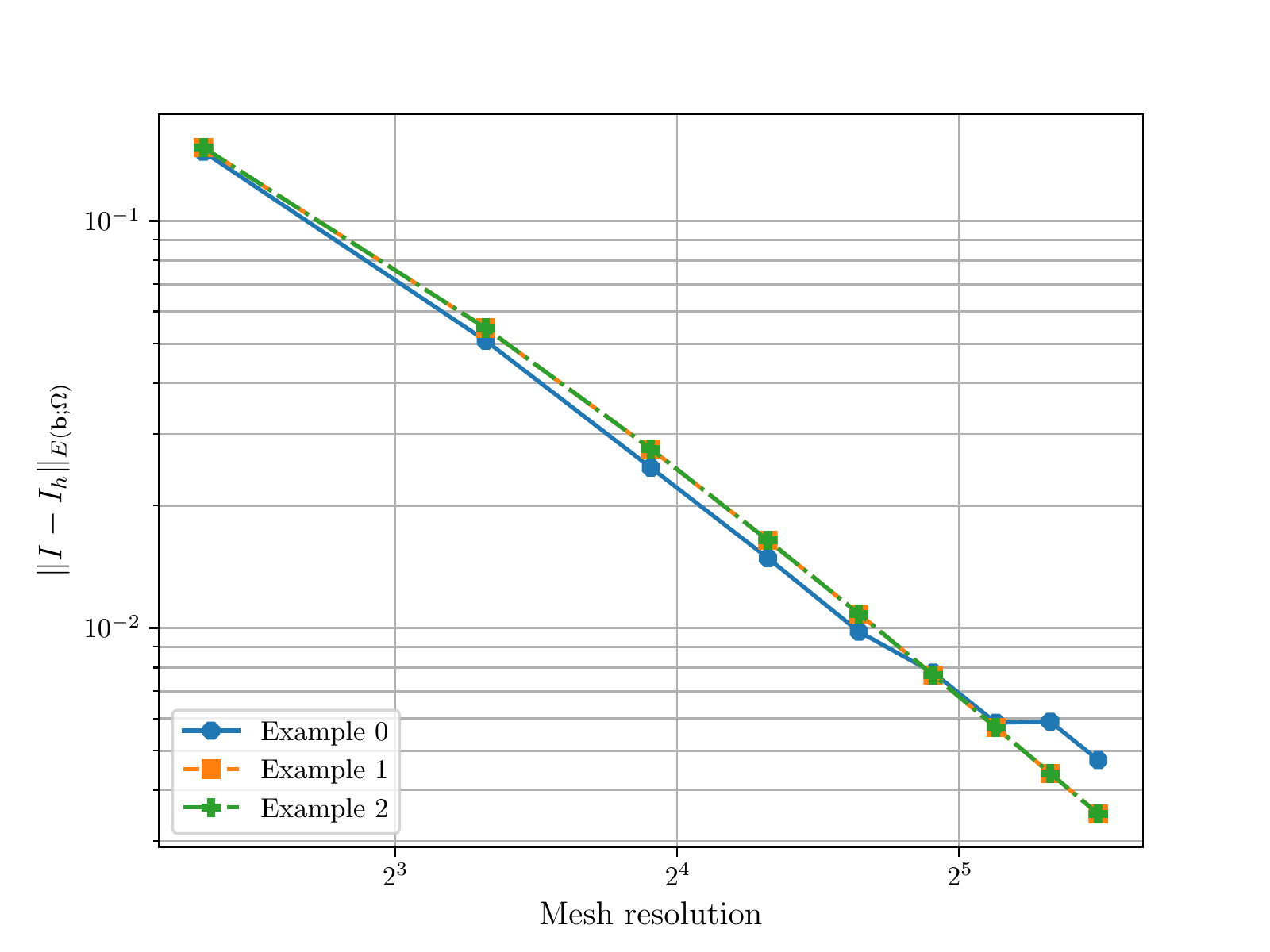}
\caption{Convergence as a function of mesh resolution
$h\inv$ for the examples in Table {\ref{tab:cg_manufactured2d}}. Left:
$L^2$ errors. Right $\Ebo$ errors.}
\label{fig:cg_inner2d}
\end{figure}

We also construct a test suite for $d=2$ for some simple advection problems.
Table \ref{tab:cg_manufactured2d} shows these manufactured solutions
\footnote{Again, only defined up to a linear function for periodicity.} with
convergence rates shown in Table \ref{tab:cg_conv0_2d}, \ref{tab:cg_conv1_2d}
and \ref{tab:cg_conv2_2d}. These rates are also depicted in Fig.
\ref{fig:cg_inner2d}.

\begin{table}[!h]
\centering
\caption{Manufactured solutions to
{\eqref{leastsquareadvection}} where $d=2$}
{%
\begin{tabular}{c|c|c}
\hline No. & $I(x,t)$ & $u(x,t)$\\
\hline 0 & $e^{-25((x-0.3(1-t)-0.6t)^2+(y-0.3(1-t)-0.6t)^2)}$ & $(0.3, 0.3)$\\
\hline 1 & $e^{-25((x-0.3(1-t)-0.7t)^2+(y-0.5)^2)}$ & $(0.4, 0)$\\
\hline 2 & $e^{-25((y-0.3(1-t)-0.7t)^2+(x-0.5)^2)}$ & $(0, 0.4)$\\
\hline 
\end{tabular}
}
\label{tab:cg_manufactured2d}
\end{table}

\begin{table}[!t]
\centering
\caption{Convergence rates for example 0 from Table
{\ref{tab:cg_manufactured2d}}}
{\begin{tabular}{c|cc|cc}
\hline\bfseries $h\inv$ & $L^2$ error & Order & Energy error & Order
\begin{filecontents*}{imgs/inner_problem/test0_rates_2d.csv}
hmesh,ltwoerr,ltwoord,energyerr,energyord
5,2.7822e-02,-,1.4797e-01,-
10,1.2917e-02,1.11,5.0807e-02,1.54
15,6.9421e-03,1.53,2.4787e-02,1.77
20,4.2600e-03,1.70,1.4856e-02,1.78
25,2.8541e-03,1.79,9.8016e-03,1.86
30,2.0860e-03,1.72,7.7759e-03,1.27
35,1.5585e-03,1.89,5.8581e-03,1.84
40,1.2864e-03,1.44,5.8927e-03,-0.04
45,1.0200e-03,1.97,4.7406e-03,1.85
\end{filecontents*}
\csvreader[head to column names]{imgs/inner_problem/test0_rates_2d.csv}{}
{\\\hline \hmesh & \ltwoerr &\ltwoord & \energyerr & \energyord}\\\hline
\end{tabular}}
\label{tab:cg_conv0_2d}
\end{table}

\begin{table}[!t]
\centering
\caption{Convergence rates for example 1 from Table
{\ref{tab:cg_manufactured2d}}}
{\begin{tabular}{c|cc|cc}
\hline\bfseries $h\inv$ & $L^2$ error & Order & Energy error & Order
\begin{filecontents*}{imgs/inner_problem/test1_rates_2d.csv}
hmesh,ltwoerr,ltwoord,energyerr,energyord
5,2.3747e-02,-,1.5156e-01,-
10,1.2861e-02,0.88,5.4566e-02,1.47
15,7.3752e-03,1.37,2.7568e-02,1.68
20,4.6565e-03,1.60,1.6426e-02,1.80
25,3.1745e-03,1.72,1.0842e-02,1.86
30,2.2905e-03,1.79,7.6685e-03,1.90
35,1.7250e-03,1.84,5.7005e-03,1.92
40,1.3432e-03,1.87,4.3992e-03,1.94
45,1.0742e-03,1.90,3.4954e-03,1.95
\end{filecontents*}
\csvreader[head to column names]{imgs/inner_problem/test1_rates_2d.csv}{}
{\\\hline \hmesh & \ltwoerr &\ltwoord & \energyerr & \energyord}\\\hline
\end{tabular}}
\label{tab:cg_conv1_2d}
\end{table}

\begin{table}[!t]
\centering
\caption{Convergence rates for example 2 from Table
{\ref{tab:cg_manufactured2d}}}
{\begin{tabular}{c|cc|cc}
\hline\bfseries $h\inv$ & $L^2$ error & Order & Energy error & Order
\begin{filecontents*}{imgs/inner_problem/test2_rates_2d.csv}
hmesh,ltwoerr,ltwoord,energyerr,energyord
5,2.3747e-02,-,1.5156e-01,-
10,1.2861e-02,0.88,5.4566e-02,1.47
15,7.3752e-03,1.37,2.7568e-02,1.68
20,4.6565e-03,1.60,1.6426e-02,1.80
25,3.1745e-03,1.72,1.0842e-02,1.86
30,2.2905e-03,1.79,7.6685e-03,1.90
35,1.7250e-03,1.84,5.7005e-03,1.92
40,1.3432e-03,1.87,4.3992e-03,1.94
45,1.0742e-03,1.90,3.4954e-03,1.95
\end{filecontents*}
\csvreader[head to column names]{imgs/inner_problem/test2_rates_2d.csv}{}
{\\\hline \hmesh & \ltwoerr &\ltwoord & \energyerr & \energyord}\\\hline
\end{tabular}}
\label{tab:cg_conv2_2d}
\end{table}

\end{document}